\documentclass[a4paper,10pt]{article}

\usepackage{subcaption}
\usepackage{amsmath}
\usepackage{amssymb}
\usepackage{amsfonts}
\usepackage{amsthm}
\usepackage{bbm}
\usepackage[mathscr]{eucal}
\usepackage{a4wide}
\usepackage{authblk}
\usepackage{float}

\usepackage{tikz}
\usetikzlibrary{calc,fadings,decorations.pathreplacing}
\usetikzlibrary{shapes,snakes}

\tikzset{%
  >=latex, 
  inner sep=0pt,%
  outer sep=2pt,%
  rootnode/.style={inner sep=0pt,outer sep=0pt,minimum size=15pt,fill=black,star,star points=7},
  root2/.style={inner sep=0pt,outer sep=0pt,minimum size=8pt,fill=black,circle},
  neighborhood2/.style={inner sep=0pt,outer sep=0pt,minimum size=7pt,draw=black,diamond},
  notneighborhood2/.style={inner sep=0pt,outer sep=0pt,minimum size=5pt,fill=black,circle},
  insidecircle/.style={inner sep=0pt,outer sep=0pt,minimum size=5pt,fill=black,circle}%
}

\usepackage{hyperref}
\hypersetup{
    colorlinks = true,
    urlcolor = {magenta},
    citecolor = {red},
    linkcolor = {blue}
}
\usepackage{doi}
\usepackage[square,numbers]{natbib}

\usepackage{xcolor}
\usepackage{todonotes}

\newcommand{\RNum}[1]{\uppercase\expandafter{\romannumeral #1\relax}}

\newcommand{\pr}{\mathbb{P}}								

\newcommand{\Prob}[1]{\pr\left(#1\right)}					
																
\newcommand{\e}{\mathbb{E}}								

\newcommand{\Exp}[1]{\e\left[#1\right]}					
															
\newcommand{\CExp}[2]{\e\left[\left.#1\right|#2\right]}	

\newcommand{\plim}{\ensuremath{\stackrel{\pr}{\rightarrow}}}	
\newcommand{\dlim}{\ensuremath{\stackrel{d}{\rightarrow}}}		

\newcommand{\bigO}[1]{O\left(#1\right)}				
\newcommand{\smallO}[1]{o\left(#1\right)}			
\newcommand{\bigT}[1]{\Theta\left(#1\right)}		

\newcommand{\1}{\mathbbm{1}}								
\newcommand{\ind}[1]{\1_{\left\{#1\right\}}}	

\newcommand{\ch}[1]{{\color{red}#1\color{black}}}

\newcommand\numberthis{\addtocounter{equation}{1}\tag{\theequation}}

\allowdisplaybreaks

\newtheorem{theorem}{Theorem}[section]
\newtheorem{lemma}[theorem]{Lemma}

\newtheorem{proposition}[theorem]{Proposition}
\newtheorem{corollary}[theorem]{Corollary}
\newtheorem{conjecture}[theorem]{Conjecture}
\newtheorem{definition}{Definition}[section]

\newtheorem{remark}{Remark}[section]

\numberwithin{equation}{section}

\title{Scaling of the clustering function in\\ spatial inhomogeneous random graphs}
\author[1]{Remco van der Hofstad, Pim van der Hoorn, Neeladri Maitra}

\affil[1]{Department of Mathematics and Computer Science,\linebreak Eindhoven University of Technology}

\begin{document}
\maketitle

\begin{abstract}
We consider an infinite spatial inhomogeneous random graph model with an integrable connection kernel that interpolates nicely between existing spatial random graph models. Key examples are versions of the weight-dependent random connection model, the infinite geometric inhomogeneous random graph, and the age-based random connection model. These infinite models arise as the local limit of the corresponding finite models. For these models we identify the asymptotics of the \emph{local clustering} as a function of the degree of the root in different regimes {in a unified way}. We show that the scaling exhibits \emph{phase transitions} as the interpolation parameter moves across different regimes. This allows us to draw conclusions on the geometry of a \emph{typical} triangle contributing to the clustering in the different regimes. 
\end{abstract}





\section{Introduction}\label{sec:intro}

%

The last two decades has seen the rise of models for complex networks that use geometry~\cite{KB_GIRGs_19,garcia2016hidden,age_dep_RCMs,hoff2002latent,krioukov2010hyperbolic}. The main idea is that geometry can encode similarity of nodes and thus serves as a component for deciding whether an edge will be present. Random graph models with underlying geometry have proven to be very successful in reproducing several key features of real-world networks. As such, geometry of complex networks has become an important research domain~\cite{bloznelis2019clustering,boguna2021network,michielan2022detecting,clust_complex_nets_1, kolk2022anomalous}. 
One key feature properly captured by geometric models is the presence of triangles, triples of connected nodes, measured in terms of their \emph{clustering}. In the literature, three different measures for clustering are often considered. The first is the so-called global clustering coefficient, defined as three times the ratio of the number of triangles to the number of paths of length two in the graph. The second is the local clustering coefficient, which for each node counts the fraction of triangles it participates in and then averages this over all nodes. Finally, the \emph{clustering function}, sometimes called the clustering spectrum, is the local clustering function restricted to nodes of a given degree $k$. That is, it is a function that takes as input $k \in \mathbb{N}$, and outputs the average number of triangles among all possible triangles that are present among degree $k$ vertices in the graph, see \eqref{eq:cLust_fn_SIRG_finite} below for a definition.  

This paper studies the behaviour of the \emph{clustering function} of a \emph{typical vertex}, in a model of infinite spatial inhomogeneous random graphs. Spatial inhomogeneous random graphs can be seen as \emph{spatial} versions of inhomogeneous random graphs, as studied in \cite{IRGs_BRJ}, where nodes have a position in a geometric space, see also \cite[Remark 1.5]{LWC_SIRGs_2020} for a discussion. In a typical construction of these graphs, one considers vertices as atoms of some marked Poisson point process, where the edges are included conditionally independently given the vertices {and their marks}, according to a connection function $\kappa$. These graphs also generalize infinite geometric inhomogeneous random graphs \cite{KB_Avg_dist_16, KB_SamplingGIRGs_17, KB_GIRGs_19} and hyperbolic random graphs \cite{KPKVB_HRG}. Graphs of similar flavour have been considered in \cite{age_dep_RCMs, WDRCM_Rec_trans}.


{For an infinite graph, a functional similar to the clustering function can be considered under some conditions.} In particular, if a sequence of finite random graphs converges locally to an infinite graph, {then} the associated sequence of clustering functions also converge to a limiting functional of the infinite graph, see \cite[Exercise 2.32]{RGCN_2}. {This limit is} expressed as the clustering function of a \emph{typical vertex} of this graph, see \eqref{eq:cLust_fn_SIRG} below. We study how the clustering function behaves as the degree of the typical vertex diverges. Problems of similar flavor have been studied in \cite{Clustering_HRGs_2020, Stegehuis_2019} for the hyperbolic random graph. In particular, {a phase transition in the clustering function is established} in \cite[Proposition 1.4]{Clustering_HRGs_2020}, which disproved the conjectured {inverse} linear scaling of the clustering function of these graphs \cite{KPKVB_HRG}. For spatial versions of preferential attachment networks, {convergence of the average and global clustering coefficients are shown in \cite{SPA_EJ_PM}}, {while} scalings of the clustering function for a different spatial preferential attachment model have been studied in \cite{CF_scaling_PAM}. Typically, the clustering function vanishes polynomially, as the degree of a typical vertex {tends} to infinity, for {sparse} network models from the literature. For example, for the models considered in \cite{CF_gen_PAMs, Props_highly_clustered_nets}, the authors see {an inverse} linear decay, while a polynomial decay is observed, for example, in \cite{struc_large_soc_net, book_dyn_large_nets, clust_complex_nets_1, top_dyn_props_int}. The main intuition behind the prediction of an inverse linear decay is that most triangles are edge-disjoint, and there can be at most $k/2$ such edge-disjoint triangles incident to a vertex if it has degree $k$. 

In this paper, we see that under suitable choices of the model parameters, even for sparse random graphs, the clustering function {can} decay logarithmically, {or not decay} at all, see for example Theorem \ref{thm:scaling} below. This last behavior is to be understood as that \emph{there are rare vertices of {arbitrarily} large degrees, and for these vertices, a constant proportion of their neighbor pairs are also neighbors}, i.e., the neighborhood of such {rare} high-degree vertices are \emph{almost} cliques - a phenomenon that has not yet been observed in the literature for sparse random graphs.  



The rest of this section is organized as follows: In Sections \ref{ssec:SIRG} and \ref{ssec:int_kernel} we define our random graph models. In Section \ref{ssec:clustering} we discuss the clustering function. We state our main result in Section \ref{ssec:main_result}, followed by a discussion on different examples covered under our setup in Section \ref{ssec:examples}. Finally in Section \ref{ssec:proof_strategy} we discuss the proof strategy of our main result.

\subsection{The random graph model}\label{ssec:SIRG}
{The random graphs we study} are spatial inhomogeneous random graphs{, or SIRGs} \cite{LWC_SIRGs_2020}. The vertex set is either $[n]=\{1,\dots,n\}$ (finite model) or $\mathbb{N} \cup \{0\}$ (infinite model). Each vertex $i$ has a location $X_i$ in $\mathbb{R}^d$, where for the finite model the locations are i.i.d.\ uniform in the box
\[
 I_n:=\left[-\frac{n^{1/d}}{2},\frac{n^{1/d}}{2}\right]^d.
\]
For the infinite model,
all the locations $ (X_i)_{i \geq 1}$ except that of $0$, form the collection of atoms of a rate $1$ homogeneous Poisson point process {$\Gamma$} on $\mathbb{R}^d$, and we let the vertex $0$ have the location $X_0:=\mathbf{0} \in \mathbb{R}^d$. {In particular, the vertex locations of $\mathbb{G}^{\infty}$ form a  palm version $\Gamma \cup \{\mathbf{0}\}$ of the Poisson process $\Gamma$.} Moreover, each vertex $i$ has a random positive real weight $W_i$ associated to it, {where $(W_i)_{i \geq 1}$ is a sequence of i.i.d.\ copies of a random variable $W$.} Edges between pairs of vertices $i$ and $j$, conditionally on $(X_i,W_i)$ and $(X_j,W_j)$, are placed independently with probability 
\[
 \kappa(\|X_i-X_j\|,W_i,W_j),
\]
where $\kappa: \mathbb{R}_+ \times \mathbb{R} \times \mathbb{R} \to [0,1]$ is assumed to be symmetric in its last two arguments. We denote the finite model as {$\mathbb{G}^{(n)}(\kappa,W)$, and the infinite model as $\mathbb{G}^{\infty}(\kappa,W)$.} The infinite rooted graph $(\mathbb{G}^{\infty},0)$ arises as the local limit in probability of the finite model, see \cite[Theorem 1.2]{LWC_SIRGs_2020}.

For the rest of the paper, we fix the sequence of weights $\mathbb{W}=(W_i)_{i\in S}$, where $S$ is either $[n]$ or $\mathbb{N}\cup\{0\}$, corresponding respectively to the finite and infinite models, to be a sequence of i.i.d.\ copies from a $\text{Pareto}(\beta-1,1)$ distribution. That is, each {$W$} has density function 
\begin{equation}\label{eq:Pareto_weights}
    f_W(t)=\frac{\beta-1}{t^{\beta}}\ind{t>1},
\end{equation}
with $\beta>2$, i.e., the weights have finite expectation. For such choices of power-law weights, with appropriate choice of the connection function $\kappa$, the degree distribution of a {vertex} also obeys a power law, so that they naturally give rise to scale-free spatial random graphs.

\subsection{The interpolation kernel} \label{ssec:int_kernel}
We now discuss the connection function $\kappa$ we are concerned with. For a \emph{long-range} parameter $\alpha \in (0, \infty]$, we consider
\begin{align*}
\kappa(r,s,t)=\kappa_{\alpha}(r,s,t)=\begin{cases}
        \ind{r<g(s,t)^{1/d}}+\left(\frac{g(s,t)}{r^d}\right)^{\alpha}\ind{r\geq g(s,t)^{1/d}} &\mbox{if}\;\;\alpha \in (0,\infty),\\
        \ind{r<g(s,t)^{1/d}} &\mbox{if}\;\; \alpha=\infty,
    \end{cases} \numberthis \label{eq:conn_fn_form}
\end{align*}
for some symmetric function $g:\mathbb{R} \times \mathbb{R} \to \mathbb{R}_+$. We refer to $g$ as the \emph{weight function} of $\kappa$. For vertices $i$ and $j$, this is to be interpreted as follows: when $\alpha<\infty$, {the edge $\{i,j\}$ is present with probability $1$ if $\|X_i-X_j\|<g(W_i,W_j)^{1/d}$,} otherwise it is present with probability $\left(\frac{g(W_i,W_j)}{\|X_i-X_j\|^d}\right)^{\alpha}$; while when $\alpha=\infty$, the edge $\{i,j\}$ is present precisely when $\|X_i-X_j\|<g(W_i,W_j)^{1/d}$.

Let us introduce the following notation which we use throughout the paper: for any two reals $r_1$ and $r_2$,
\[r_1 \vee r_2 := \max\{r_1,r_2\},\;\;\text{and}\;\;r_1 \wedge r_2 = \min\{r_1,r_2\}.\]

In this paper, we work with the weight function
\begin{equation}\label{defn:interpolated_wt_fn}
    g(s,t)=g_a(s,t)=(s\vee t)(s \wedge t)^a,
\end{equation}
for a parameter $a \in [0 ,\infty)$. {Different values of the interpolation parameter $a$ in \eqref{defn:interpolated_wt_fn} capture different spatial random graphs from the literature, see Section \ref{ssec:examples} for details.}


\subsection{The clustering function}\label{ssec:clustering}

Let $G=(V(G),E(G))$ be a {locally} finite graph. For $v \in V(G)$, {we denote its degree in $G$ by $d_v(G)$ and further define}
\begin{equation}\label{eq:defn_Clust_Func}
    \Delta_v(G):=\sum_{v_1,v_2 \in V(G)}\ind{\{v,v_1\}\in E(G)}\ind{\{v,v_2\}\in E(G)}\ind{\{v_1,v_2\}\in E(G)},
\end{equation}
{so that} $\Delta_v(G)$ is the number of triangles in $G$ containing $v$ as a vertex. {In addition, we} let $N_k(G)$ be the number of vertices in $G$ with $d_v(G)=k$.

\begin{definition}[Clustering function]
The clustering function $\mathrm{CC}_G:\mathbb{N}\to\mathbb{R}_+$ of $G$ is defined as 
\begin{equation}\label{eq:cLust_fn_SIRG_finite}
    \mathrm{CC}_G(k):=\begin{cases}
        \frac{1}{N_k(G)}\sum_{v \in V(G): d_v(G)=k} {\frac{\Delta_v(G)}{k(k-1)}} &\mbox{if}\;\;0<N_k(G)< \infty,\\
        0 &\mbox{otherwise.}
        \end{cases}
\end{equation}
\end{definition}

Thus, the quantity $\mathrm{CC}_G(k)$ measures the proportion of two-paths in $G$ that are triangles, where the middle vertex of the two-path has degree $k$. 

It is not immediately clear how to define such a notion for an infinite graph $G$, since $N_k(G)$ can be infinite. However, if $G$ arises as the local limit of a growing sequence $(G_n)_{n \geq 1}$ of finite graphs, then the corresponding sequence of clustering functions {has} a limit, which can then be defined as the clustering function of $G$. {Indeed,} it was shown in \cite[Corollary 2.10]{LWC_SIRGs_2020} that, for any $k \in \mathbb{N}$,
\begin{align*}
    \mathrm{CC}_{\mathbb{G}^{(n)}}(k) \plim \gamma(k),
\end{align*}
where 
\begin{equation}\label{eq:cLust_fn_SIRG}
    \gamma(k)=\frac{1}{\binom{k}{2}}\CExp{\Delta_0(\mathbb{G}^{\infty})}{d_0(\mathbb{G}^{\infty})=k}.
\end{equation}
We call $\gamma(k)$ the \emph{clustering function} of $\mathbb{G}^{\infty}$. In this paper we {analyse} the asymptotics of $\gamma(k)$ for large $k$.

\subsection{Main result}\label{ssec:main_result}

Let us denote the Lebesgue measure of a ball of unit radius in $\mathbb{R}^d$ by $\omega_d$. For any $\alpha \in (1,\infty]$ we then define
\begin{equation}\label{eq:def_xi}
	\xi_\alpha = \xi_{\alpha,\beta,d} := \frac{\alpha \omega_d (\beta-1)}{\alpha-1},
\end{equation}
where $\xi_\infty := \lim_{\alpha \to \infty} \xi_\alpha = \omega_d (\beta-1)$. For any $w_1, w_2 \in (1,\infty)$, we define the functions
\begin{align*}
	I(w_1,w_2) &:= \iint_{(\mathbb{R}^d)^2} \kappa(\|\mathbf{x}\|, 1, w_1) \kappa(\|\mathbf{x}\|, 1, w_2)
			\kappa(\|\mathbf{y}\|, w_1, w_2) d\mathbf{x} d\mathbf{y}, \text{ and}\\
	\mathcal{S}_1(w_1,w_2) &:= \iint_{(\mathbb{R}^d)^2} \kappa(\|\mathbf{x}\|, 1, w_1) \kappa(\|\mathbf{y}\|, 1, w_2)
		\kappa(\|\mathbf{x}-\mathbf{y}\|, w_1, w_2) d\mathbf{x} d\mathbf{y}.
\end{align*}
Further, let $\mathrm{S}_k(a,\beta)$ be defined as
\begin{equation}\label{eq:scaling_fn}
    \mathrm{S}_k(a,\beta)=\begin{cases}
        k^{-1} &\mbox{if}\;\;\beta>(a+\frac{3}{2})\vee 2, a \in (0,\infty),\\
        k^{-1}\log{k} &\mbox{if}\;\; \beta=a+\frac{3}{2}, a\in (\frac{1}{2},\infty),\\
        k^{2a+2-2\beta} &\mbox{if}\;\; \beta\in(a+1,a+\frac{3}{2}), a\in (\frac{1}{2},\infty),\\
        (\log{k})^{-2} &\mbox{if}\;\; \beta=a+1, a\in (1,\infty),\\
        1 &\mbox{if}\;\; \beta \in (2,a+1), a\in (1,\infty).
    \end{cases}
\end{equation}
Next for $\mathbf{w}=(w_1,w_2)\in \mathbb{R}^2$, using the notation 
\begin{equation}\label{eq:notn_two_dim_weight_integral}
    f_W(\mathbf{w})d\mathbf{w}:=f_W(w_1)dw_1f_W(w_2)dw_2
\end{equation}
for the product measure on $\mathbb{R}^2$ with density $f_W(w_1)f_W(w_2)$, let
$\Gamma(a,\alpha,\beta,d)$ be defined as
\begin{equation}\label{eq:all_constants}
\Gamma(a,\alpha,\beta,d) = \begin{cases}
			\frac{\beta-a-1}{\xi_\alpha} 
			\hspace*{46pt}\iint_1^\infty I(w_1,w_2) f_W(\mathbf{w})d\mathbf{w} 
			&\mbox{if}\;\;\beta>\left(a+\frac{3}{2}\right)\vee 2, a > 0,\\
	        \frac{2a+1}{4\xi_\alpha} 
	        \hspace*{54pt} {(\beta-1)^2\int_0^\infty I(1,r) r^{-\beta}dr} 
	        &\mbox{if}\;\; \beta=a+\frac{3}{2}, a > \frac{1}{2},\\
	        \left(\frac{\beta-a-1}{\xi_\alpha}\right)^{4+2a-2\beta}
	        \iint_0^\infty \mathcal{S}_1(w_1,w_2) f_W(\mathbf{w})d\mathbf{w} 
	        &\mbox{if}\;\; \beta \in (a+1,a+\frac{3}{2}), a > \frac{1}{2},\\
	        \xi_\alpha^{-1}, 
	        \hspace*{58pt}\iint_0^\infty \mathcal{S}_1(w_1,w_2) f_W(\mathbf{w})d\mathbf{w} 
	        &\mbox{if}\;\; \beta=a+1, a > 1,\\
	        \left(\frac{(1+a-\beta)(\beta-2)}{(a-1)\xi_\alpha}\right)^2 
	        \iint_0^\infty \mathcal{S}_1(w_1,w_2) f_W(\mathbf{w})d\mathbf{w} 
	        &\mbox{if}\;\; \beta \in (2,a+1), a > 1.
	\end{cases}    
\end{equation}
With these definitions we are ready to state our main result about the scaling of the clustering function:

\begin{theorem} \label{thm:scaling}
Let $\alpha>1$, $\beta>2$ and let $W$ be a random variable with density~\eqref{eq:Pareto_weights}. Furthermore, let $\kappa$ be defined as in~\eqref{eq:conn_fn_form} with $g$ as in \eqref{defn:interpolated_wt_fn}.
Consider the infinite SIRG $\mathbb{G}^{\infty}(\kappa,W)$ and associated clustering function $\gamma(k)$ as defined in~\eqref{eq:cLust_fn_SIRG}. Then, as $k \to \infty$,
\begin{align*}
    \frac{\gamma(k)}{\mathrm{S}_k(a,\beta)} \to \Gamma(a,\alpha,\beta,d) \in (0,\infty),
\end{align*}
where $\mathrm{S}_k(a,\beta)$ is as in \eqref{eq:scaling_fn}, with $\Gamma(a,\alpha,\beta,d)$ is as in \eqref{eq:all_constants}.
\end{theorem}

The main message of Theorem~\ref{thm:scaling} is that the clustering function of the infinite SIRG $\mathbb{G}^\infty$ has \emph{three} different scaling regimes, depending on the parameters $a$ and $\beta$, characterized by the function $\mathrm{S}_k(a,\beta)$. 

\begin{figure}[h]
    \centering
    \resizebox{0.7\textwidth}{!}{
        \begin{tikzpicture}[xscale=2,
  pinknode/.style={shape=rectangle, draw=pink, line width=8},
  purplenode/.style={shape=rectangle, draw=purple, line width=8},
  yellownode/.style={shape=rectangle, draw=yellow, line width=8},
  greennode/.style={shape=rectangle, draw=green, line width=8},
  cyannode/.style={shape=rectangle, draw=cyan, line width=8}
]
\def\figheight{5};
\def\figwidth{3};

\coordinate (O) at (0,0);
\coordinate (xaxis) at (\figwidth,0);
\coordinate (yaxis) at (0,\figheight);
\coordinate (xmark1) at (0.5,0);
\coordinate (xmark2) at (1,0);
\coordinate (ymark1) at (0,1);
\coordinate (ymark2) at (0,1.5);
\coordinate (ymark3) at (0,2);
\coordinate (intersect1) at (0.5,2);
\coordinate (intersect2) at (1,2);
\coordinate (endpoint1) at (\figwidth,\figheight);
\coordinate (endpoint2) at (\figwidth,4.5);
\coordinate (endpoint3) at (\figwidth,4);
\coordinate (endpoint4) at (\figwidth,3);
\coordinate (endpoint5) at (\figwidth, 2);

\path[fill=pink] (ymark3) -- (intersect1) -- (endpoint2) -- (endpoint1) -- (yaxis) -- cycle;
\path[fill=yellow] (intersect1) -- (intersect2) -- (endpoint3) -- (endpoint2) -- cycle;
\path[fill=cyan] (intersect2) -- (endpoint5) -- (endpoint3) -- cycle;
 
\draw[-stealth,thick] (O) -- (xaxis) node[right, align=center] {$a$};
\draw[-stealth,thick] (O) -- (yaxis) node[above, align=center] {$a=0$\\$\beta$} ;

\draw[dashed] (xmark1) node[below, align=center] {$(\frac{1}{2},0)$} -- (intersect1) ;
\draw[dashed] (xmark2) node[below, align=center] {$(1,0)$} -- ++(0,\figheight)node[above, align=center] {$a=1$ };
\draw[dashed] (ymark3) node[left, align=center] {$(0,2)$} -- ++(\figwidth,0);

\draw[thick] (ymark3) -- (endpoint1) node[right,align=center] {$\beta=a+2$};
\draw[purple,thick] (intersect1) -- (endpoint2) node[right,align=center] {$\beta=a+\frac{3}{2}$};
\draw[green,thick] (intersect2) -- (endpoint3) node[right,align=center] {$\beta=a+1$};
\draw[orange,thick] (intersect2) -- (endpoint4) node[right,align=center] {$\beta=\frac{a+3}{2}$};

\matrix [draw,below right] at (current bounding box.east) {
  \node [pinknode,label=right:{$\mathrm{S}_k(a,\beta)=k^{-1}$}] {}; \\
  \node [purplenode,label=right:{$\mathrm{S}_k(a,\beta)=k^{-1}\log{k}$}] {}; \\
  \node [yellownode,label=right:{$\mathrm{S}_k(a,\beta)=k^{2a+2-2\beta}$}] {}; \\
  \node [greennode,label=right:{$\mathrm{S}_k(a,\beta)=\frac{1}{\log{k}^2}$}] {}; \\
  \node [cyannode,label=right:{$\mathrm{S}_k(a,\beta)=\text{constant}$}] {}; \\
};

\end{tikzpicture}
    }
    \caption{Phase diagram of the scaling function $\mathrm{S}_k(a,\beta)$.}
    \label{fig:phase_diagram_scalings}
\end{figure}
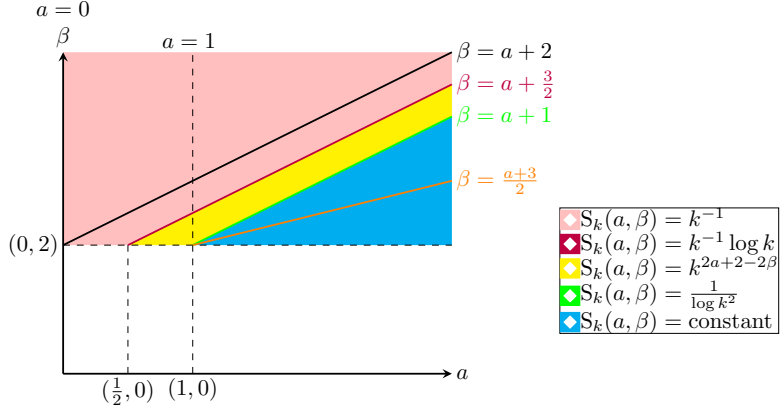

In Figure \ref{fig:phase_diagram_scalings}, we provide the phase diagram of the scaling function $\mathrm{S}_k(a,\beta)$ in terms of $a$ and $\beta$. The lines $a=0$, $a=1$ and $\beta=a+2$ all correspond to known random graph models, see Section \ref{ssec:examples} for a discussion on the different examples of random graph models covered under our set up.

\begin{remark}[Range of parameters]\label{rmk:parameters}
Let us discuss our choice of parameters $\beta>2$ and $\alpha>1$. 
When $\alpha\leq 1$, the conditional degree of $0$ given $\{W_0=w\}$, is almost surely infinite for any $w>0$ by Lemma \ref{lem:nbr_PP}, which implies that $d_0(\mathbb{G}^{\infty})$ is almost surely infinite. Similarly, when $\alpha>1$ and $\beta\leq 2$, $d_0(\mathbb{G}^{\infty})$ is almost surely infinite, again by Lemma \ref{lem:nbr_PP}. In particular our choice of the range of parameters $\beta>2$, $\alpha>1$ are so that the random variable $d_0(\mathbb{G}^{\infty})$ is finite almost surely, and the quantity of interest $\gamma(k)$ makes sense.
\end{remark}




\begin{remark}[Phase transitions]
We note that there are three different phases for the scaling of the clustering function $\gamma(k)$: $\beta>(a+\frac{3}{2})\vee 2$, $\beta \in (a+1,a+\frac{3}{2})$ and $\beta \in (2,a+1)$, with two phase transition boundaries at $\beta=a+\frac{3}{2}$ and $\beta=a+1$. The scalings of $\gamma(k)$ at $\beta=a+1$ and $\beta \in (2,a+1)$ are new, especially the constant scaling of $\gamma(k)$ in the regime $\beta \in (2,a+1)$ is surprising. The phase transition at $\beta=a+\frac{3}{2}$ has already been observed for a particular case: namely $a=1$ and $d=1$, with $\alpha=\infty$. {In this case the} corresponding infinite SIRG $\mathbb{G}^{\infty}$ is a good approximation to the local limit of the threshold hyperbolic random graph model with appropriate parameters (see e.g. \cite[Section 2.1.3]{LWC_SIRGs_2020}), using arguments as appearing in \cite[Section 9]{JK_BL_explosions_HRGs_19}. For these random graphs, the phase transition at the point $(1,5/2)$ of Figure \ref{fig:phase_diagram_scalings} was observed in \cite[Proposition 1.4]{Clustering_HRGs_2020}.    
\end{remark}

\begin{remark}[Choice of Pareto weights]
Our choice of the weight distribution \eqref{eq:Pareto_weights} is purely for technical reasons. One can also consider regularly-varying distributions with tail parameter $\beta-1$, but Pareto weights make the arguments {and constants} cleaner. We expect Theorem \ref{thm:scaling} to go through when one has regularly varying weights {with similar proofs, albeit more technical, except in the critical regimes $\beta=a+\frac{3}{2}$ and $\beta=a+1$. Here} we expect the particular choice of the regularly varying distribution to introduce additional slowly-varying corrections, which change the scaling. The choice of Pareto weights also helps us draw a cleaner link between the models of SIRG we consider and existing random graph models from the literature, as discussed in Section \ref{ssec:examples}. {All these random graph models have a mixed-Poisson degree distribution, and hence have the same power-law tail as the mixing random variable (e.g. see \cite[Section 2.2]{LWC_SIRGs_2020}) }   
\end{remark}

\subsection{Examples}\label{ssec:examples}
We now discuss different spatial random graph models from the literature that different values of $a$ in \eqref{defn:interpolated_wt_fn} give rise to. 

\paragraph{Product kernels.} When one takes $a=1$ in \eqref{defn:interpolated_wt_fn}, $g(s,t)$ becomes the product $st$. This corresponds to the models having a product structure in their weights, most notably Geometric Inhomogeneous Random Graphs (GIRGs) \cite{KB_Avg_dist_16, KB_SamplingGIRGs_17,KB_GIRGs_19}, Weight Dependent Random Connection Models (WDRCMs) \cite{WDRCM_Rec_trans}, and Hyperbolic Random Graphs (HRGs)  \cite{KPKVB_HRG}. Let us discuss these examples briefly.

The infinite GIRG (see \cite[Section 2.1.2]{LWC_SIRGs_2020}) has as its vertex set a homogeneous Poisson process on $\mathbb{R}^d$. Each vertex has a power-law weight associated to it, where the collection of weights form an i.i.d.\ collection. Conditionally on everything, the edges are placed independently with a connection function of the form \eqref{eq:conn_fn_form} with $g(s,t)=st$. When restricted to $d=1$, these graphs give a  version of HRGs \cite{KPKVB_HRG} - for a discussion of this, see \cite[Section 9]{JK_BL_explosions_HRGs_19}. Here we comment that typically in GIRGs or HRGs, the connection function does not exactly have the form \eqref{eq:conn_fn_form}, but can easily be seen to be bounded from above and below up to constant factors by \eqref{eq:conn_fn_form}. In particular, assuming the form \eqref{eq:conn_fn_form} for the connection function does not affect the conclusions of Theorem \ref{thm:scaling}, but comes with less unnecessary technicalities, which is why we prefer \eqref{eq:conn_fn_form}.

The $a=1$ case can also be seen as a special case of the \emph{weight-dependent random connection model} \cite{WDRCM_Rec_trans}. Here the vertices are again a homogeneous Poisson point process on $\mathbb{R}^d$, but each vertex now has a uniformly distributed $[0,1]$ mark associated to it, as opposed to power-law weights in the case of GIRGs, and each edge between two vertices $(\mathbf{x},s)$ and $(\mathbf{y},t)$ is included conditionally independently with probability $\rho(g(s,t)\|\mathbf{x}-\mathbf{y}\|^d)$, where the function $\rho(\cdot)$ is non-increasing and integrable.  In particular, let the \emph{profile function} $\rho(g(s,t)\|\mathbf{x}-\mathbf{y}\|^d)$ in \cite[(1.2)]{WDRCM_Rec_trans} be equal to $\kappa_{\alpha}(\|\mathbf{x}-\mathbf{y}\|,s,t)$, with weight function $g$ as in \eqref{defn:interpolated_wt_fn} with $a=1$. This corresponds to taking $\beta=1$ in the product kernel \cite[(1.6)]{WDRCM_Rec_trans} of these models, and taking $1/\gamma=\beta-1$, where now $\beta$ is the tail parameter of the Pareto distribution coming from \eqref{eq:Pareto_weights}. This can be seen using the standard fact that if $U$ is uniform on $(0,1)$, then $U^{-\gamma}$ has the density \eqref{eq:Pareto_weights} with $1/\gamma=\beta-1$.

\paragraph{The max kernel.} When one takes $a=0$ in \eqref{defn:interpolated_wt_fn}, the function $g(s,t)$ {becomes} $s \vee t$. This also corresponds to a special case of the \emph{ weight-dependent random connection model} \cite{WDRCM_Rec_trans}. To realize this case, we take $\beta=1$ for the \emph{min kernel} $g^{\mathrm{min}}(s,t)=(s \wedge t)/\beta$ as in \cite[p3]{WDRCM_Rec_trans}, and let $1/\gamma=\beta-1$ as before. This transforms the connection function $\rho(g(s,t)\|\mathbf{x}-\mathbf{y}\|)$ into the connection function $\kappa_{\alpha}$ as in \eqref{eq:conn_fn_form}, with $a=0$ and $1/\gamma=\beta-1$, with weight function $g$ as in \eqref{defn:interpolated_wt_fn} with $a=0$.

\paragraph{Age-dependent random connection models.} When one takes $a=\beta-2$, the corresponding infinite model $\mathbb{G}^{\infty}$ gives rise to the \emph{age-dependent random connection model}, which comes out as the local limit of the \emph{age-based spatial preferential attachment model}, as studied in \cite{age_dep_RCMs}. This model also has a homogeneous Poisson point process as its vertex set, with i.i.d.\ marks associated to them that are uniformly distributed on $[0,1]$. The connection function (see \cite[Section 3]{age_dep_RCMs}) for this model has the form $\varphi(\beta^{-1}(s \vee u)^{1-\gamma}(s \wedge u)^{\gamma}\|\mathbf{x}-\mathbf{y}\|^d)$ (see e.g. the bottom of \cite[p315]{age_dep_RCMs}), for parameters $\beta \in (0, \infty)$, $\gamma \in (0,1)$, $s,t$ being i.i.d.\ uniforms on $(0,1)$, and a non-increasing, integrable profile function $\varphi$. Let $\beta=1$, and consider the profile function $\varphi(r)=\ind{r<1}+\left(\frac{1}{r}\right)^{\alpha}\ind{r \geq 1}$. For these choices, we can write the above connection function as $\kappa_{\alpha}(r,s',t')$, as appearing in \eqref{eq:conn_fn_form}, with the weight function $g_a(s',t')$ being $(s'\vee t')(s' \wedge t')^a$, where $s'=s^{-\gamma}$, $t'=t^{-\gamma}$, and $a=1/\gamma-1$. As before, using the standard fact that if $U$ is uniform on $(0,1)$, then $W=U^{-\gamma}$ is a Pareto random variable with scale parameter $1$ and shape parameter $1/\gamma$, if we let $1/\gamma=\beta-1$, then using $\gamma \in (0,1)$, $\beta>2$, both $s'$ and $t'$ have the density \eqref{eq:Pareto_weights}, and the weight function $g_a$ is \eqref{defn:interpolated_wt_fn} with $a=\beta-2$. This completes the link between the infinite SIRGs we consider, and the age-based spatial preferential attachment networks.

\medskip
A similar kernel was also considered in \cite{finiteness_one_dim_Lukas}, where the authors look at percolation questions for one dimensional models, and in {the recent} work \cite{JJ_JK_DM_CSD_1}, the authors look at a similar interpolating kernel, and study cluster-size decay of these models.

\subsection{Proof strategy}\label{ssec:proof_strategy}

In this section we discuss the main ideas to prove Theorem \ref{thm:scaling}. {First, let us recall some standard notations. Throughout the paper, for real sequences $(f_k)_{k \geq 1}$ and $(g_k)_{k \geq 1}$, we use the standard Landau notations $f_k=\bigO{g_k}$ and $f_k=\Omega\left( g_k\right)$, as $k \to \infty$, to mean respectively that there exist real constants $C,c$ such that
\begin{align*}
    \limsup_{k \to \infty}\left|\frac{f_k}{g_k}\right|\leq C,\;\;\text{and}\;\;\liminf_{k \to \infty}\left|\frac{f_k}{g_k}\right|\geq c.
\end{align*}
We write $f_k=\Theta(g_k)$ when both $f_k=\bigO{g_k}$ and $f_k=\Omega\left( g_k\right)$. Further, we write $f_k=\smallO{g_k}$ and $f_k=\omega(g_k)$ to mean respectively that
\begin{align*}
    \limsup_{k \to \infty} \frac{f_k}{g_k}=0\;\;\text{and}\;\;\liminf_{k \to \infty}\frac{f_k}{g_k}=\infty.
\end{align*}
}

The first step is to reduce the expression of $\gamma(k)$ using properties of the locations and weights of the vertices of the graph $\mathbb{G}^{\infty}$. Let us for now denote the conditional density of the weight $W_0$ of the root $0$ given that its degree $d_0$ equals $k$ by $w\mapsto f_k(w)$. 
{Using this conditional density}, we obtain the following representation of $\gamma(k)$,
\[
\gamma(k)=\int_{\mathbb{R}_+}\frac{T(w)}{\mathcal{M}(w)^2}f_k(w)dw,
\]
where the quantity ${T(w)}/{\mathcal{M}(w)^2}$ is a bounded function of $w$. Here $T(w)$ can be understood as the expected number of triangles containing $0$ given $W_0=w$, and $\mathcal{M}(w)$ is the expected degree of $0$ given $W_0=w$. This representation is obtained in Proposition \ref{prop:simplifying_gamma}. We then use concentration properties of Poisson random variables to show that asymptotically, as $k \to \infty$, the density $f_k(w)$ behaves like a Dirac density, with a point mass at $\mathcal{M}^{-1}(k)$. This is the content of Proposition \ref{prop:conc_weights}. As a consequence,
\begin{equation}\label{eq:proof_strategy_approx_gamma_T}
	\gamma(k) = (1+o(1)) \frac{T(\mathcal{M}^{-1}(k))}{k^2},
\end{equation}
as $k \to \infty$, and thus it suffices to show that
\[
	\frac{T(\mathcal{M}^{-1}(k))}{k^2 \mathrm{S}_k(a,\beta)} \to \Gamma(a,\alpha,\beta,d).
\]

%
We then proceed in two phases. First in Section~\ref{sec:m_inv_scaling} we establish the exact scaling of $\mathcal{M}^{-1}(k)$. More precisely, if we define
\[
    \sigma(w)=\begin{cases}
       	w &\mbox{if } \beta> \left(a + \frac{3}{2}\right) \vee 2, a > 0,\\
        w\log{\left(w \right)} &\mbox{if } \beta=a + \frac{3}{2}, a > \frac{1}{2},\\  
        w^{4+2a-2\beta} &\mbox{if } \beta \in (2,a+ \frac{3}{2}), a > \frac{1}{2}.
        \end{cases}
\]
then our result implies that $\sigma(\mathcal{M}^{-1}(k))/(k^2 \mathrm{S}_k(a,\beta))$ converges to a constant, see Corollary~\ref{cor:limits_M_inv_k}.

For the final part we write
\[
	\frac{T(\mathcal{M}^{-1}(k))}{k^2 \mathrm{S}_k(a,\beta)} = \frac{\sigma(\mathcal{M}^{-1}(k))}{k^2 \mathrm{S}_k(a,\beta)}
	\frac{T(\mathcal{M}^{-1}(k))}{\sigma(\mathcal{M}^{-1}(k))}.
\]
The first term converges by Corollary~\ref{cor:limits_M_inv_k}, so that all that is left is to study the limit $T(w)/\sigma(w)$ as $w \to \infty$ in more detail. To prepare for this, we establish general upper and lower bounds on $T(w)$ and study their asymptotic behavior in Section~\ref{sec:T_bounds}. Finally, in Section~\ref{sec:main_thm_proof}, we establish the limits $T(w)/\sigma(w)$ for the three different cases and wrap up the proof of Theorem~\ref{thm:scaling}.

\subsection{More on the phase transitions}

{

As Theorem~\ref{thm:scaling} shows, the clustering function undergoes two phase transitions. While the exact scaling of the clustering function in the different regimes, as well as the exact location of the transition, depends on the choice of the function $g$ in $\kappa$, there is something general we can say about the driving force behind both phase transitions. 

We start with the first transition. To discuss this, for $w,{w}_1,{w}_2>0$ we introduce the notation
\begin{equation}\label{eq:notn_h_wst}
    h(w,{w}_1,{w}_2):=g(w_1,{w}_2)[g(w,{w}_1)\wedge g(w,{w}_2)].
\end{equation}
Recall from~\eqref{eq:proof_strategy_approx_gamma_T} that the scaling of $\gamma(k)$ is determined by that of $T(\mathcal{M}^{-1}(k))$. In Section~\ref{ssec:proof_strategy} we mentioned that we will establish general upper and lower bounds on $T(w)$. 

In Lemma~\ref{lem:N_UB}, we show that
\begin{align*}
    T(w)
    &\leq \xi_\alpha^2 \int_{1}^{\infty}\int_{1}^{\infty} h(w,{w}_1,{w}_2)\ind{{w}_1\wedge {w}_2\leq w}({w}_1{w}_2)^{-\beta}d{w}_1d{w_2}\\
    &\hspace{10pt} + \xi_\alpha^2\int_{1}^{\infty}\int_{1}^{\infty}g(w,{w}_1)g(w,{w}_2) \ind{{w}_1 \wedge{w}_2>w}({w}_1{w}_2)^{-\beta}d{w}_1d{w}_2,
\end{align*}
while in Lemma~\ref{lem:N_LB}, we obtain
\begin{align*}
    &T(w) \geq \omega_d^2\int_{1}^{\infty} \int_{1}^{\infty} \ind{g({w}_1,{w}_2)<g(w,{w}_2)} h(w,{w}_1,{w}_2)\\
    &\hspace{100 pt} \times \left[1-\left(\frac{g({w}_1,{w}_2)}{(g(w,{w}_1) \wedge g(w,{w}_2)}\right)^{1/d} \right]^d ({w}_1{w}_2)^{-\beta}d{w}_1d{w}_2.
\end{align*}

We remark that when $w$ is sufficiently large, up to constant factors, the first term of the upper bound, as well as the lower bound, are both truncated integrals of the function $h(w,{w}_1,{w}_2)$. For the first term of the upper bound, this is immediate, while for the lower bound we use that, for $w$ sufficiently large, the term
\[
	\left[1-\left(\frac{g({w}_1,{w}_2)}{(g(w,{w}_1) \wedge g(w,{w}_2)}\right)^{1/d} \right]^d
\]
is of constant order. In particular, as we see in the proof of Theorem \ref{thm:scaling}, the upper and lower bounds are of the same order for large $w$, which means that the scaling of $T(w)$, for large $w$, is guided by truncated integrals of the function $h(w,{w}_1,{w}_2)$. As we see in the proof of Theorem \ref{thm:scaling}, it is the phase transition in the integrability of the function $h(w,{w}_1,{w}_2)$ that causes the first phase transition in the scaling of $\gamma(k)$ at $\beta=a+\frac{3}{2}$, see Theorem \ref{thm:scaling}. In particular, when $\beta>a+\frac{3}{2}$, the function $h(w,{w}_1,{w}_2)$, when multiplied by the densities of ${w}_1$ and ${w}_2$ having the form \eqref{eq:Pareto_weights}, is integrable in ${w}_1$ and ${w}_2$, and the obtained function of $w$ decays linearly in $w$, which is the main cause for the {inverse} linear scaling of $\gamma(k)$ in this regime. Instead when $\beta \leq a+\frac{3}{2}$, the function $h(w,{w}_1,{w}_2)$ is not integrable anymore, and its truncated integral diverges logarithmically when $\beta=a+\frac{3}{2}$, and polynomially when $\beta<a+\frac{3}{2}$.
 
Finally, we discuss the occurence of the second phase transition at $\beta=a+1$. Here we understand the phase transition in terms of how the scaling of $\mathcal{M}(w)$ changes. In particular, as $w\to \infty$, one observes the following transition in behaviour of the scaling of the mean degree $\mathcal{M}(w)$ of $0$ given $W_0=w$, about $\beta=a+1$: 
\begin{align*}
     \mathcal{M}(w)=\begin{cases}
        \Theta(w) &\mbox{if}\;\; \beta>(a+1)\vee 2, a\in [0,\infty),\\
        \Theta(w \log{w}) &\mbox{if}\;\; \beta=a+1, a\in (1,\infty),\\
        \Theta(w^{2+a-\beta}) &\mbox{if}\;\; \beta \in (2,a+1), a\in (1,\infty),
    \end{cases}
\end{align*}
see Figure \ref{fig:scaling_M} and Section \ref{sec:m_inv_scaling} for more details on the behaviour of $\mathcal{M}(w)$.

In particular, whenever $\beta<a+\frac{3}{2}$, ${T(\mathcal{M}^{-1}(k))}/{k^2}$ always scales like ${\mathcal{M}^{-1}(k)^{4+2a-2\beta}}/{k^2},$ as we will prove in Section \ref{sec:m_inv_scaling}. Thus, the transition in behaviour of the scaling of $\mathcal{M}^{-1}(k)$ guides the transition in scaling of ${T(\mathcal{M}^{-1}(k))}/{k^2}$ around $\beta=a+1$.}

\paragraph{Organization of the paper.} Section \ref{sec:prelims} gives preliminary simplifications to analyse $\gamma(k)$. {In Section \ref{sec:m_inv_scaling}, we analyse the scaling of $\mathcal{M}^{-1}(k)$.} In Section \ref{sec:T_bounds}, we develop some bounding techniques which aid us in analysing upper and lower asymptotics of $\gamma(k)$. Using results from Sections \ref{sec:m_inv_scaling} and \ref{sec:T_bounds}, we prove Theorem \ref{thm:scaling} in Section \ref{sec:main_thm_proof}. Section \ref{sec:disc} is devoted to discussion. {Appendix A contains some standard results which are used to derive a result of Section \ref{sec:prelims}. Appendix B contains the proof of a technical result used in Section \ref{sec:m_inv_scaling}}.

\section{Preliminaries}\label{sec:prelims}

This section is devoted to two simplifications of the clustering function $\gamma(k)$ from~\eqref{eq:cLust_fn_SIRG}. First, in Proposition~\ref{prop:simplifying_gamma} we reduce its expression to an integral form as given in~\eqref{eq:gamma_integral_form}. Next, in Proposition~\ref{prop:conc_weights} we prove a concentration result which enables us to argue that asymptotically, this integral is really an integral with respect to a Dirac measure. This simplifies the analysis for its asymptotics.

\subsection{Simplifying $\gamma(k)$}

{Let us denote
\begin{align*}
    \mathbb{U}:=\mathbb{R}^d \times \mathbb{R}_+
\end{align*} for the parameter space of the vertices: the first coordinate denotes the location and the second the weight of a vertex.}
{We first prove that conditionally on the weight of a typical vertex, the neighbors of the typical vertex in any infinite SIRG form an inhomogeneous Poisson point process on $\mathbb{U}$}:


\begin{lemma}[Neighbors of $0$]\label{lem:nbr_PP}
{Consider an infinite SIRG $\mathbb{G}^{\infty}$, where the vertex locations are given by a Palm version $\Gamma \cup \{\mathbf{0}\}$ of a homogeneous Poisson process $\Gamma$ on $\mathbb{R}^d$, the weights are i.i.d.\ with density $f_W$, and the connection function $\kappa: \mathbb{R}_+\times \mathbb{R} \times \mathbb{R} \to [0,1]$ is symmetric in its last two arguments. Let us denote the weight of $0$ by $W_0$. Then,} conditionally on $\{W_0=w\}$, the locations and weights of the neighbors of $0$ in the graph $\mathbb{G}^{\infty}$ are distributed as the atoms of an inhomogeneous Poisson point process $\mathcal{N}_0^{(w)}$ on $\mathbb{U}$, with intensity function, for any $\mathbf{p}=(\mathbf{x},x) \in \mathbb{U}$, given by{
\begin{equation}\label{eq:density_bias}
    \rho_w(\mathbf{p})=\rho_w(\mathbf{x},x):=\kappa(\|\mathbf{x}\|,w,x)f_W(x).
\end{equation}
}
 
\end{lemma}
{The proof of Lemma \ref{lem:nbr_PP} is a consequence of standard theory of markings and thinnings of Poisson point processes. The required necessary concepts and the proof is postponed to Appendix A.}

{Recall $\gamma(k)$ from \eqref{eq:cLust_fn_SIRG}. Let us use Lemma \ref{lem:nbr_PP} to simplify the expression of $\gamma(k)$, when $\mathbb{G}^{\infty}$ has independent weights distributed according to the density \eqref{eq:Pareto_weights}, and with connection function \eqref{eq:conn_fn_form}, for some symmetric weight function $g$. Denoting in this case the integral of the above intensity function by $\mathcal{M}(w)$, we have,
\begin{equation}\label{eq:mean_degree_given_wt}
	\mathcal{M}(w) = \int_{\mathbb{U}}\rho_w(\mathbf{p})d\mathbf{p} 
	= \int_{1}^{\infty} \int_{\mathbb{R}^d}\kappa_{\alpha}(\|\mathbf{x}\|,w,x)d\mathbf{x}f_W(x)dx.
\end{equation}
%
Moreover, Lemma \ref{lem:nbr_PP} implies that the number of neighbors of $0$ conditionally on the event $\{W_0=w\}$ has distribution
\begin{equation}\label{eq:deg_law_given_wt}
 d_0(\mathbb{G}^{\infty}) \big|_{W_0=w} \stackrel{d}{=} \text{Poi}\left(\mathcal{M}(w) \right).
\end{equation}
}

Using \eqref{eq:deg_law_given_wt}, a simple application of the law of total probability shows that the conditional density function $f_k(w)$ of the weight $W_0$ of $0$, conditionally on the event $\{d_0(\mathbb{G}^{\infty})=k\}$ has the form
\begin{equation}\label{eq:cond_density}
    f_k(w) :=\frac{\Prob{\text{Poi}(\mathcal{M}(w))=k}f_W(w)}{\int_{1}^{\infty}\Prob{\text{Poi}(\mathcal{M}(r))=k}f_W(r)dr} 
    =\frac{e^{-\mathcal{M}(w)}\mathcal{M}(w)^k f_W(w)}{\int_{1}^{\infty}e^{-\mathcal{M}(r)}\mathcal{M}(r)^k f_W(r)dr}.
\end{equation}
From now on, to simplify notation, we write $d\underline{\mathbf{x}}=d\mathbf{x}_1d\mathbf{x}_2$ for {the} Lebesgue measure on $(\mathbb{R}^d)^2$. Also recall the notation $f_W(\underline{\mathbf{w}})d\underline{\mathbf{w}}$ from \eqref{eq:notn_two_dim_weight_integral}. 
 
\begin{proposition}[Simplifying $\gamma(k)$]\label{prop:simplifying_gamma}
{Consider an infinite SIRG $\mathbb{G}^{\infty}$ with vertex locations distributed as $\Gamma \cup \{\mathbf{0}\}$ where $\Gamma$ is a homogeneous Poisson point process on $\mathbb{R}^d$, with independent vertex weights distributed according to the density \eqref{eq:Pareto_weights}, and with connection function \eqref{eq:conn_fn_form} for some symmetric weight function $g$. Assume that the integral on the RHS of \eqref{eq:mean_degree_given_wt} is finite. Then} the clustering function $\gamma(k)$ of $\mathbb{G}^{\infty}$ has the form
\begin{equation}\label{eq:gamma_integral_form}
    \gamma(k)=\int_{1}^{\infty}\frac{T(w)}{\mathcal{M}(w)^2}f_k(w)dw,
\end{equation}
where {$f_k(w)$ is as in \eqref{eq:cond_density}, and}

\begin{equation}\label{eq:def_T_w}
	T(w) := \iint_0^\infty \iint_{(\mathbb{R}^d)^2} \kappa_{\alpha}(\|\mathbf{x}_1\|,{w}_1,w) \kappa_{\alpha}(\|\mathbf{x}_2\|,{w}_2,w) \kappa_{\alpha}(\|\mathbf{x}_1-\mathbf{x}_2\|,{w}_1,{w}_2) 
	d\underline{\mathbf{x}}	\, f_W(\underline{\mathbf{w}})d\underline{\mathbf{w}}.
\end{equation}

\end{proposition}

\begin{proof}

Recall the Poisson process $\mathcal{N}^{(w)}_0$ from Lemma \ref{lem:nbr_PP}. Let $\mathcal{N}_{0}^{(w,k)}$ be the conditional point process on $\mathbb{U}$ defined as
\begin{align*}
    \mathcal{N}_{0}^{(w,k)}:=\mathcal{N}^{(w)}_0 \big|_{\mathcal{N}^{(w)}_0(\mathbb{U})=k}.
\end{align*}

Let $(U_{\{i,j\}})_{i,j=1}^k$ be a symmetric array of {uniformly distributed random variables on $[0,1]$,} independent of everything else.

Denoting the atoms of $\mathcal{N}_0^{(w,k)}$ by $\{\mathbf{p}_i=(\mathbf{x}_i,{x}_i)\in \mathbb{U}: 1\leq i \leq k\}$, we note that 
\begin{align*}
     \CExp{\Delta_0(\mathbb{G}^{\infty})}{d_0(\mathbb{G}^{\infty})=k, W_0=w}
     = \Exp{\sum_{\mathbf{p}_i,\mathbf{p}_j \in \mathcal{N}_0^{(w,k)}}^{\neq} \ind{U_{\{i,j\}} < \kappa_{\alpha}(\|\mathbf{x}_i-\mathbf{x}_j\|,{x}_i,{x}_j)}},
\end{align*}
where $\sum_{\mathbf{p}_i,\mathbf{p}_j \in \mathcal{N}_0^{(w,k)}}^{\neq}$ denotes the sum over distinct pairs $\mathbf{p}_i,\mathbf{p}_j \in \mathcal{N}_0^{(w,k)}$ and the expectation is taken {both} with respect to the point process and the uniform random variables. This implies {that}
\begin{equation}\label{eq:gamma_int_exp_prelim}
    \CExp{\Delta_0(\mathbb{G}^{\infty})}{d_0(\mathbb{G}^{\infty})=k}
    =\int_{1}^{\infty}\Exp{\sum_{\mathbf{p}_i,\mathbf{p}_j \in \mathcal{N}_0^{(w,k)}}^{\neq} \ind{U_{\{i,j\}} < \kappa_{\alpha}(\|\mathbf{x}_i-\mathbf{x}_j\|,{x}_i,{x}_j)}}f_k(w)dw.
\end{equation}

Using the notation $\mathbf{p}_i=(\mathbf{x}_i,x_i)$, we note that the collection of $\mathbb{U}$-valued random variables $\{\mathbf{p}_i\}_{\mathbf{p}_i \in \mathcal{N}_0^{(w,k)}}$ form an i.i.d.\ collection of random variables with density given by{
\begin{equation}\label{eq:nbr_density}
    \mathbf{p}\mapsto h_w(\mathbf{p})=h_w(\mathbf{x},x)
    = \frac{\kappa(\|\mathbf{x}\|,w,x) f_W(x)}
    {\int_1^\infty \int_{\mathbb{R}^d} \kappa(\|\mathbf{x}\|,w,x) f_W(x) d \mathbf{x} \, dx}
    = \frac{\kappa(\|\mathbf{x}\|,w,x) f_W(x)}{\mathcal{M}(w)}.
\end{equation}
In particular, the expectation of the sum in the integrand in \eqref{eq:gamma_int_exp_prelim} evaluates to
\begin{align*}
    &\binom{k}{2}\int_{\mathbb{U}} \int_{\mathbb{U}} \kappa_{\alpha}(\|\mathbf{x}_1-\mathbf{x}_2\|,x_1,x_2)h_w(\mathbf{p}_1) h_w(\mathbf{p}_2)d(\mathbf{p}_1)d(\mathbf{p}_2),\\
    &= \frac{\binom{k}{2}}{\mathcal{M}(w)^2} \int_1^\infty \int_1^\infty\iint_{(\mathbb{R}^d)^2}\kappa_{\alpha}(\|\mathbf{x}_1\|,w_1,w)\kappa_{\alpha}(\|\mathbf{x}_2\|,w_2,w)\kappa_{\alpha}(\|\mathbf{x}_1-\mathbf{x}_2\|,w_1,w_2) d\underline{\mathbf{x}} \, f_W(\underline{\mathbf{w}})d\underline{\mathbf{w}}\\
    &=\frac{\binom{k}{2}T(w)}{\mathcal{M}(w)^2},
\end{align*}
where $\mathbf{p}_1=(\mathbf{x}_1,x_1),\mathbf{p}_2=(\mathbf{x}_2,x_2)$, and to obtain the first equality we have simply changed variables $(x_1,x_2)\mapsto (w_1,w_2)$.} Since $\gamma(k)=\binom{k}{2}^{-1}\CExp{\Delta_0(\mathbb{G}^{\infty})}{d_0(\mathbb{G}^{\infty})=k}$ from \eqref{eq:cLust_fn_SIRG}, the result follows. 
\end{proof}

\subsection{Concentration of the weight of $0$ when $d_0(\mathbb{G}^{\infty})$ is large}

{The main goal of this section is to argue that the conditional density $w \mapsto f_k(w)$ of \eqref{eq:cond_density} is asymptotically like a Dirac mass at $\mathcal{M}^{-1}(k)$.} Informally, the argument goes as follows: Conditionally on the Poisson variable \eqref{eq:deg_law_given_wt} being large, {standard Poisson concentration arguments} imply that its mean $\mathcal{M}(w)$ also has to be large. In particular, if $\mathcal{M}(w)$ is a nice monotone function of $w$, this means that $w$ also has to be large. That is, the degree of $0$ in $\mathbb{G}^{\infty}$ cannot be large unless its weight is also large. {In this section, we formalise this by arguing that} when the degree of $0$ is $k$, and $k \to \infty$, the weight of $0$ is highly concentrated, and diverges {as $\mathcal{M}^{-1}(k)$}. As a consequence, the conditional density function $f_k(w)$ appearing in the integral \eqref{eq:gamma_integral_form} behaves like a Dirac mass with {an} atom at {$\mathcal{M}^{-1}(k)$}. {Consequently, using \eqref{eq:gamma_integral_form}, the $k \to \infty$ asymptotic behaviour of $\gamma(k)$ is the same as that of ${T(\mathcal{M}^{-1}(k))}/{k^2}$.}

Recall the notations $T(w)$ and $\mathcal{M}(w)$ from \eqref{eq:mean_degree_given_wt} and \eqref{eq:def_T_w}. First, we need a regularity condition for diverging sequences, that ensures that they are robust under {certain lower order} perturbations:

\begin{definition}[Fluctuation-robust sequences]\label{defn:log_robust}
A sequence $(\psi(k))_{k \geq 1}$ diverging to infinity is called \emph{fluctuation-robust}, if for any $C>0$, as $k \to \infty$,
\begin{align*}
    \frac{\left(\psi(k)+C\sqrt{\psi(k)\log{\psi(k)}}\right)\vee \psi\left(k+C\sqrt{k \log{k}}\right)}{\left(\psi(k)-C\sqrt{\psi(k)\log{\psi(k)}}\right) \wedge \psi\left(k-C\sqrt{k\log{k}}\right)} \to 1.
\end{align*}
\end{definition}
Examples of fluctuation-robust sequences are positive powers of $k$, or some positive power of $k$ times some logarithmic factor: $k \log{k}$, $k^{\epsilon} (\log{k})^{\delta}$ for any $\epsilon>0$ and $\delta \in \mathbb{R}$, etc. It is easy to verify that for any constant $C>0$ and $L \in \mathbb{N}$, all of $C\psi(k)$, $\psi(Ck)$, $\psi^{(L)}(k)$ are fluctuation-robust, if $\psi(k)$ is so, where $\psi^{(L)}(\cdot)$ denotes the $L-$fold composition of $\psi(\cdot)$. Note that the sequence $\psi(k)=e^{k}$ is not fluctuation-robust. Observe that if $\psi(k)$ is fluctuation-robust, then any sequence $\phi(k)$ with ${\psi(k)}/{\phi(k)}\to 1$ is also fluctuation-robust.

Let us now state the main result of this section. Recall the connection function \eqref{eq:conn_fn_form}. For the following proposition, we do not need the form \eqref{defn:interpolated_wt_fn} of $g$, and it is true as soon as the weight function $g$ satisfies certain assumptions:
\begin{proposition}\label{prop:conc_weights}
{Let $\mathbb{G}^{\infty}$ be an infinite SIRG with vertex locations distributed as $\Gamma \cup \{\mathbf{0}\}$ for a homogeneous Poisson point process $\Gamma$ on $\mathbb{R}^d$, independent vertex weights having density \eqref{eq:Pareto_weights}, and with the connection function \eqref{eq:conn_fn_form} for some symmetric function $g$. Recall also $\mathcal{M}(w)$ from \eqref{eq:mean_degree_given_wt}.} 
Assume the following:
\begin{itemize}
    \item[1.] For fixed $s \in \mathbb{R}$, the function {$w \mapsto g(s,w)$ is monotone increasing in $w$.}
    \item[2.] {$\mathcal{M}(w)< \infty$ for all $w >0$, with $\mathcal{M}(w) \to \infty$ as $w \to \infty$. Further, there exists $ K_0> 0$ and $\zeta \in \mathbb{R}$ such that $\mathcal{M}(w)$ is strictly increasing and differentiable in $(K_0,\infty)$, with $\mathcal{M}'(w)\leq w^{\zeta}$ for all $w>K_0$.} 
    {\item[3.] There exists $\xi \in \mathbb{R}$ such that $\mathcal{M}^{-1}(k)\leq k^{\xi}$, as $k \to \infty$.} 
    {\item[4.]} $T(\mathcal{M}^{-1}(k))=\Omega(k^{-\eta})$ for some $\eta>0$, as $k \to \infty$.
\end{itemize}
Under the above mentioned assumptions,

\begin{itemize}
    \item[(a.)] if further $T(\mathcal{M}^{-1}(k))=\Theta(\psi(k))$ for some fluctuation-robust sequence $(\psi(k))_{k \in \mathbb{N}}$, then as $k \to \infty$,
    \begin{align*}
        \gamma(k)=\Theta(T(\mathcal{M}^{-1}(k))k^{-2}).
    \end{align*}
    \item[(b.)] if further the sequence $(T(\mathcal{M}^{-1}(k)))_{k \in \mathbb{N}}$ is itself fluctuation-robust, then, as $k \to \infty$,
    \begin{align*}
        \frac{\gamma(k)}{T(\mathcal{M}^{-1}(k))k^{-2}} \to 1.
    \end{align*}
\end{itemize}

\end{proposition}

Thus, under the conditions of Proposition \ref{prop:conc_weights}, asymptotically the density $w \mapsto f_k(w)$ appearing in the integral \eqref{eq:gamma_integral_form} behaves like a Dirac mass at $\mathcal{M}^{-1}(k)$, as $k \to \infty$, provided the sequence $(T(\mathcal{M}^{-1}(k)))_{k \in \mathbb{N}}$ scales like a fluctuation-robust sequence.

{Proposition \ref{prop:conc_weights} is used in the proof of Theorem \ref{thm:scaling} in the following way: Note that the sequence $\mathrm{S}_k(a,\beta)$ in the statement of Theorem \ref{thm:scaling} is fluctuation robust, for given $a$ and $\beta$. Consequently, so is the sequence $k^2\mathrm{S}_k(a,\beta)$. Hence it is enough to show that $\frac{T(\mathcal{M}^{-1}(k))}{k^2 \mathrm{S}_k(a,\beta)} \to \Gamma(a,\alpha, \beta,d)$, since this automatically implies that $T(\mathcal{M}^{-1}(k))$ is fluctuation-robust, and hence the appropriately rescaled limit of both $T(\mathcal{M}^{-1}(k))/k^2$ and $\gamma(k)$ is the same, due to Proposition \ref{prop:conc_weights} (b).}

Before going into the proof, {we first state some useful results.} {We begin with} the standard  Chernoff concentration bound for Poisson random variables, e.g. see \cite[(2.12)]{Clustering_HRGs_2020}, as
\begin{equation}\label{eq:basic_Chernoff}
    \Prob{|\text{Poi}(\lambda)-\lambda| \geq x} \leq 2 \exp{\left(-\frac{x^2}{2(\lambda+x)}\right)},
\end{equation}
for any $x > 0$. In particular, for any $C>0$, and $\lambda=\lambda_n \to \infty$,
\begin{equation}\label{eq:div_Chernoff}
    \Prob{|\text{Poi}(\lambda_n)-\lambda_n| \geq C\sqrt{\lambda_n\log{\lambda_n}}} = \bigO{\lambda_n^{-\frac{C^2}{2}}}.
\end{equation}

We also recall the standard stochastic domination property of Poisson random variables: for $\lambda_1>\lambda_2$ and any $x>0$,
\begin{equation}\label{eq:poisson_stoc_dom} 
    \Prob{\text{Poi}(\lambda_2)>x} \leq \Prob{\text{Poi}(\lambda_1)>x}.
\end{equation}

{Next, we state a lemma which will be useful in the proof of Proposition \ref{prop:conc_weights}:
\begin{lemma}\label{lem:conc_lem_LB_degreetail}
    Under the assumptions of Proposition \ref{prop:conc_weights}, there exists $K_1>0$ such that for all $k>K_1$,\[\int_1^{\infty}\Prob{\text{Poi}(\mathcal{M}(w))=k}w^{-\beta}dw \geq k^{-\xi\beta-\zeta \xi-2}.\] 
\end{lemma}

We are now ready to prove Proposition \ref{prop:conc_weights}.

}

\begin{proof}[Proof of Proposition \ref{prop:conc_weights}]
Let $I_k^{-}:=k-C\sqrt{k \log{k}}$, and define $I_k^{+}$ analogously, where $C$ is a sufficiently large constant to be chosen later. 

Note that $J_k^{-}:=\mathcal{M}^{-1}(I_k^{-})$ and $J_k^{+}:=\mathcal{M}^{-1}(I_k^{+})$ are well defined for all large $k$, since $\mathcal{M}(w)$ is continuous and strictly increasing on $(K_0,\infty)$.

Further, note that if $w<J_k^{-}$, i.e., $\mathcal{M}(w)<I_k^{-}$, then
\begin{align*}
    \Prob{\mathrm{Poi}(\mathcal{M}(w))=k}
    & \leq \Prob{\mathrm{Poi}(\mathcal{M}(w))> I_k^{-}+C\sqrt{I_k^{-}\log{I_k^{-}}}}\\
    & \leq \Prob{\mathrm{Poi}(I_k^{-})> I_k^{-}+C\sqrt{I_k^{-}\log{I_k^{-}}}}=\bigO{k^{-C^2/2}},
\end{align*}
as $k \to \infty$, where to obtain the first inequality, we have used that $I_k^{-}+C\sqrt{I_k^{-}\log{I_k^{-}}}<k$ for all large $k$, and to obtain the last inequality, we have used $\mathcal{M}(w)<I_k^{-}$ along with \eqref{eq:div_Chernoff}. A similar argument shows that, for $w>J_k^{+}$, 
\begin{align*}
    \Prob{\mathrm{Poi}(\mathcal{M}(w))=k}=\bigO{k^{-C^2/2}}.
\end{align*}
We conclude that if $w \notin [J_k^-,J_k^+]$, 
\begin{equation}\label{eq:actual_up_lo_tail_gen}
    \Prob{\mathrm{Poi}(\mathcal{M}(w))=k}=\bigO{k^{-C^2/2}}
\end{equation}
Hence, we can write
\begin{equation}\label{gamma_ovgamma_comp_error}
    \frac{\gamma(k)}{T(\mathcal{M}^{-1}(k))k^{-2}}=\frac{\overline{\gamma}(k)}{T(\mathcal{M}^{-1}(k))k^{-2}}+\varepsilon_k,
\end{equation}
where
\begin{equation}\label{eq:ovgamma}
    \overline{\gamma}(k):=\int_{J_k^{-}}^{J_k^{+}}\frac{T(w)}{\mathcal{M}(w)^2}f_k(w)dw,
\end{equation}
and
\begin{align*}
    \varepsilon_k&:=\frac{1}{T(\mathcal{M}^{-1}(k))k^{-2}}\int_{(1,\infty)\setminus [J_k^-,J_k^+]}\frac{T(w)}{\mathcal{M}(w)^2}f_k(w)dw\\
    &\leq \frac{1}{T(\mathcal{M}^{-1}(k))k^{-2}} \frac{1}{\int_{1}^{\infty}\Prob{\text{Poi}(\mathcal{M}(r))=k}r^{-\beta}dr} \int_{(1,\infty)\setminus [J_k^-,J_k^+]}\Prob{\text{Poi}(\mathcal{M}(w))=k}w^{-\beta}dw. 
\end{align*}
Here to obtain the first inequality, we have used the easily checked upper bound of $1$ on $T(w)/{\mathcal{M}(w)^2}$.
{Observe that using Assumption (4), Lemma \ref{lem:conc_lem_LB_degreetail}, and \eqref{eq:actual_up_lo_tail_gen}, for all large $k$,
\begin{align*}
    \varepsilon_k=\bigO{k^{2-\eta+\xi\beta+\zeta\xi+2-C^2/2}}=\smallO{1},\numberthis \label{eq:comp_error_neg_gen}
\end{align*}
by choosing $C>0$ sufficiently large.
}

{
{Next we aim to bound $\overline{\gamma}(k)$ from above and below.} Recall that we assume that $w \mapsto g_s(w)$ is monotone increasing in $w$. As a consequence, for any $x \in \mathbb{R}^d$ and $s \in \mathbb{R}_+$, $w \mapsto \kappa_{\alpha}(\|\mathbf{x}\|,s,w)$ is also monotone increasing in $w$. In particular, this implies that $w \mapsto T(w)$ is monotone increasing in $w$. Also, $w \mapsto \mathcal{M}(w)$ is strictly increasing in $w$ whenever $w$ is sufficiently large by assumption. Using the monotonicity of $T(w)$ and $\mathcal{M}(w)$, and recalling $\overline{\gamma}(k)$ from \eqref{eq:ovgamma}, we get can an upper bound on it by simply bounding by the supremum value of the integrand,
\[
\overline{\gamma}(k)\leq \frac{T(J_k^+)}{\mathcal{M}(J_k^-)^2}\int_{J_k^-}^{J_k^+}f_k(w)dw=\frac{T(\mathcal{M}^{-1}(I_k^+))}{(I_k^-)^2}(1-\delta_k),
\]
where
\begin{align*}
    1-\delta_k:= \int_{J_k^-}^{J_k^+}f_k(w)dw
    &=\frac{\int_{J_k^-}^{J_k^+}\text{Poi}(\mathcal{M}(w)=k)w^{-\beta}dw}{\int_{1}^{\infty}\text{Poi}(\mathcal{M}(r)=k)r^{-\beta}dr}\\
    &=1-\frac{\int_{(1,\infty)\setminus[J_k^-,J_k^+]}\text{Poi}(\mathcal{M}(r)=k)r^{-\beta}dr}{\int_{1}^{\infty}\text{Poi}(\mathcal{M}(w)=k)w^{-\beta}dw}\\\
    &=1+\bigO{k^{(-C^2/2)+{\xi\beta+\zeta \xi+2}}} \to 1,\numberthis \label{eq:del_k_to_1_gen}
\end{align*}
as $k \to \infty$, {by choosing $C>0$ large,} where we have used \eqref{eq:actual_up_lo_tail_gen} {and Lemma \ref{lem:conc_lem_LB_degreetail}} to obtain \eqref{eq:del_k_to_1_gen}.

Similarly, we can also obtain a lower bound on $\overline{\gamma}(k)$, by bounding it from below by the infimum value of the integrand. In particular, we obtain the following chain of inequalities,
\begin{equation}\label{eq:UB_LB_conclusion_gen}
    \frac{k^2}{T(\mathcal{M}^{-1}(k))}\frac{T(\mathcal{M}^{-1}(I_k^-))}{(I_k^+)^2}(1-\delta_k) \leq \frac{\overline{\gamma}(k)}{T(\mathcal{M}^{-1}(k))k^{-2}} \leq \frac{k^2}{T(\mathcal{M}^{-1}(k))}\frac{T(\mathcal{M}^{-1}(I_k^+))}{(I_k^-)^2}(1-\delta_k).
\end{equation}

}

Now let us prove (a.). Since $T(\mathcal{M}^{-1}(k))=\Theta(\psi(k))$, there are constants $C_1,c_1>0$ such that whenever $k$ is sufficiently large,
\begin{align*}
    c_1\leq \frac{T(\mathcal{M}^{-1}(k))}{\psi(k)}\leq C_1.
\end{align*}
Hence,
\begin{align*}
 \frac{T(\mathcal{M}^{-1}(I_k^-))}{T(\mathcal{M}^{-1}(k))}\leq \frac{C_1\psi(I_k^-)}{c_1\psi(k)},   
\end{align*}
which implies, as $\psi(k)$ is fluctuation-robust and consequently ${\psi(I_k^-)}/{\psi(k)}\to 1$, that\\ ${T(\mathcal{M}^{-1}(I_k^-))}/{T(\mathcal{M}^{-1}(k))} = \bigO{1}$, as $k \to \infty$. A similar argument shows that\\ ${T(\mathcal{M}^{-1}(I_k^-))}/{T(\mathcal{M}^{-1}(k))} = \Omega(1)$, as $k \to \infty$. Hence, as $k \to \infty$, ${T(\mathcal{M}^{-1}(I_k^-))}/{T(\mathcal{M}^{-1}(k))} = \Theta(1)$, and similarly, ${T(\mathcal{M}^{-1}(I_k^+))}/{T(\mathcal{M}^{-1}(k))} = \Theta(1)$. Further, note from the definitions of $I_k^+$ and $I_k^-$ that both ${k}/{I_k^+}, {k}/{I_k^-} \to 1$ as $k \to \infty$. Also, \eqref{eq:del_k_to_1_gen} shows that $\delta_k \to 0$. Hence, both the left-most and the right-most expressions of \eqref{eq:UB_LB_conclusion_gen} are $\Theta(1)$, as $k \to \infty$. Hence, by \eqref{eq:UB_LB_conclusion_gen},
\begin{equation}\label{eq:ov_gamma_const_order}
    \frac{\overline{\gamma}(k)}{T(\mathcal{M}^{-1}(k))k^{-2}} = \Theta(1).
\end{equation}
Finally, recall \eqref{gamma_ovgamma_comp_error}. Since \eqref{eq:comp_error_neg_gen} shows that $\varepsilon_k \to 0$, we have both that ${\overline{\gamma}(k)}/{(T(\mathcal{M}^{-1}(k))k^{-2})}$ and ${\gamma(k)}/{(T(\mathcal{M}^{-1}(k))k^{-2})}$ are of the same order, which concludes the proof of (a), using \eqref{eq:ov_gamma_const_order}.

Now, we move on to proving (b). The philosophy is quite similar to that of the proof of (a). In this case, recall that $T(\mathcal{M}^{-1}(k))$ is itself log-robust. This implies that ${T(\mathcal{M}^{-1}(I_k^-))}/{T(\mathcal{M}^{-1}(k))}\to 1$ and $ {T(\mathcal{M}^{-1}(I_k^+))}/{T(\mathcal{M}^{-1}(k))} \to 1$. Thus arguing similarly as before, both the left-most and the right-most expressions of \eqref{eq:UB_LB_conclusion_gen} converge to $1$, as $k \to \infty$. So that also ${\overline{\gamma}(k)}/{(T(\mathcal{M}^{-1}(k))k^{-2})} \to 1$ as $k \to \infty$. Recalling \eqref{gamma_ovgamma_comp_error} and \eqref{eq:comp_error_neg_gen}, we obtain that ${\gamma(k)}/{(T(\mathcal{M}^{-1}(k))k^{-2})} \to 1$, as $k \to \infty$, which concludes the proof of (b). 
\end{proof} 
{Finally, let us give the proof of Lemma \ref{lem:conc_lem_LB_degreetail}:
\begin{proof}[Proof of Lemma \ref{lem:conc_lem_LB_degreetail}]
Let $K_1>0$ be sufficiently large so that $\mathcal{M}(w)$ is strictly increasing and differentiable on $(K_1,\infty)$, and further, $\mathcal{M}'(w)\leq w^{\zeta}$ and $\mathcal{M}^{-1}(k)\leq 2k^{\xi}$ for all $w \in (K_1,\infty)$. Note that the existence of such a $K_1$ is guaranteed by Assumptions (2.) and (3.) of Proposition \ref{prop:conc_weights}. Restricting the domain of integration to $(K_1,\infty)$ and then changing variables as $t\mapsto \mathcal{M}(w)$, we get a lower bound on $\int_1^{\infty}\Prob{\text{Poi}(\mathcal{M}(w)=k)}w^{-\beta}dw$ as,
\begin{align*}
    \int_{\mathcal{M}(K_1)}^{\infty}\frac{t^k}{k!}e^{-t}(\mathcal{M}^{-1}(t))^{-\beta}\frac{dt}{\mathcal{M}'(\mathcal{M}^{-1}(t))}&\geq \frac{2^{-\beta-\zeta}}{k!}\int_{\mathcal{M}(K_1)}^{\infty} e^{-t}t^{(k-\xi\beta-\xi\zeta-1)+1}dt\\&=\Omega\left(\frac{\Gamma(k-\xi\beta-\xi\zeta-1,\mathcal{M}(K_1))}{\Gamma(k+1)} \right),
\end{align*}
where $\Gamma(t,x)=\int_x^{\infty}e^{-r}r^{t-1}dr$ is the upper incomplete gamma function, and $\Gamma(t)=\Gamma(t,0)$ is the gamma function. Next we recall standard scaling relations between gamma functions: for any fixed $x>0$ and $a \in \mathbb{R}$, as $t \to \infty$,
\begin{align*}
    \frac{\Gamma(t+a,x)}{\Gamma(t)}=\frac{\Gamma(t+a)}{\Gamma(t)}-\frac{\int_{0}^xe^{-r}r^{t+a-1}dr}{\Gamma(t)}=\frac{\Gamma(t+a)}{\Gamma(t)}+\smallO{1}=t^{a}+\smallO{1}.
\end{align*}
Taking $t=k+1$ and $a=-\xi\beta-\xi\zeta-2$ finishes the proof of the lemma.
\end{proof}

}

\section{Scaling of $\mathcal{M}^{-1}(k)$}\label{sec:m_inv_scaling}

In the previous section, we have reduced the analysis of the clustering function $\gamma(k)$ to analyzing $T(\mathcal{M}^{-1}(k))/k^2$. The first step is to study $\mathcal{M}^{-1}(k)$, which is the content of this section. Recall the definition from~\eqref{eq:mean_degree_given_wt}. Then by \eqref{eq:Pareto_weights},
\begin{align*}
    \mathcal{M}(w) 
    &=\int_{1}^{\infty} \int_{\mathbb{R}^d}\kappa_{\alpha}(\|\mathbf{x}\|,w,x)d\mathbf{x}f_W(x)dx\\ 
    &=\begin{cases}
		\int_{1}^{\infty} \int_{\mathbb{R}^d} \left[\ind{\|\mathbf{x}\|<g(x,w)^{1/d}}+\left(\frac{g(x,w)}{\|\mathbf{x}\|^d}\right)^{\alpha}\ind{\|\mathbf{x}\|\geq g(x,w)^{1/d}}\right]d\mathbf{x}\,x^{-\beta}dx, 
		&\text{if } \alpha \in (1,\infty),\\
		\int_{1}^{\infty} \int_{\mathbb{R}^d}\ind{\|\mathbf{x}\|<g(x,w)^{1/d}}d\mathbf{x}\,x^{-\beta}dx, 
		&\text{if } \alpha=\infty.
	\end{cases}\\
	&=\begin{cases}
		\frac{\alpha \omega_d(\beta-1)}{\alpha-1}\int_{1}^{\infty} g(x,w)x^{-\beta}dx, &\text{if } \alpha \in (1,\infty),\\ \omega_d(\beta-1)\int_{1}^{\infty} g(x,w)x^{-\beta}dx, &\text{if } \alpha=\infty.
	\end{cases}
\end{align*}

Next, recalling $\xi_\alpha = \xi_{\alpha,\beta,d} := \frac{\alpha \omega_d (\beta - 1)}{\alpha - 1}$ from~\eqref{eq:def_xi}, with $\xi_\infty = \lim_{\alpha \to \infty} \xi_\alpha$, and from \eqref{defn:interpolated_wt_fn} that $g(w_1,w_2) =(w_1 \vee w_2)(w_1 \wedge w_1)^a$, we conclude that
\begin{equation}\label{eq:M_w_integral_expression}
	\mathcal{M}(w) = \xi_\alpha \left(\int_{1}^w wx^{a-\beta}dx+\int_{w}^{\infty}x^{1-\beta}w^adx \right).
\end{equation}
Direct computations of this integral yield the following result about the asymptotic behavior of {$\mathcal{M}(w)$}:

\begin{lemma}\label{lem:scaling_M_w}
Let $\mathcal{M}(w)$ be defined as in~\eqref{eq:mean_degree_given_wt} and define
\[
	m(w) = \begin{cases}
		w			&\text{if } \beta > (a + 1)\vee 2, a \ge 0,\\
		w \log w	&\text{if } \beta = a+1, a > 1,\\
		w^{2 + a -\beta}	&\text{if } \beta \in (2,a+1), a > 1.
	\end{cases}
\]
Then,
\[
	\lim_{w \to \infty} \frac{\mathcal{M}(w)}{m(w)} = \begin{cases}
			\frac{\xi_\alpha}{\beta - a - 1}	&\text{if } \beta > (a+1)\vee 2, a \ge 0,\\
			\xi_\alpha 				&\text{if } \beta = a + 1, a > 1,\\
			\frac{(a-1)\xi_\alpha}{(\beta - 2)(a + 1 - \beta)} &\text{if } \beta \in (2,a+1), a>1.
		\end{cases}
\]
\end{lemma}

\begin{remark}[Infinite mean degree]\label{rmk:infinite_degree}
Observe that if $\beta \in (2,\frac{a+3}{2})$ then $2 + a - 2\beta < -1$. Hence, since in this regime $\mathcal{M}(w) = \bigT{w^{2+a-\beta}}$, it follows that $\mathcal{M}(w)$ is no longer integrable with respect to the density $f_W(w) = \bigT{w^{-\beta}}$. This implies that $\Exp{d_0(\mathbb{G}^\infty)} = \infty$ in this regime.
\end{remark}

The above result shows that, asymptotically, $\mathcal{M}(w)$ is some polynomial of combinations of $w$ and $\log(w)$.
The next step is then to use the asymptotic behavior of $\mathcal{M}(w)$ to derive the behavior of $\mathcal{M}^{-1}(k)$. Let us define certain functions which have the form of $\mathcal{M}(w)$. To begin with, let us call a real function \emph{eventually strictly increasing} if it is strictly increasing on an interval of the form $(K_0,\infty)$ for some $K_0 \in \mathbb{R}$. Then we define:
\begin{definition}[Almost power functions]
A function $f: \mathbb{R}_+ \to \mathbb{R}$ is called an \emph{almost power function}, if it is an eventually strictly increasing continuous function, such that
\begin{equation}\label{eq:almost_power}
    \lim_{x \to \infty} \frac{f(x)}{x^a \log(x)^b} = c \in (0,\infty),
\end{equation}
for some $a > 0$ and $b\ge 0$.
\end{definition}
We use a lemma that relates the scaling of almost power functions to the scaling of their inverses:
\begin{lemma}\label{lem:log_poly_inv_scaling}
Let $f: \mathbb{R}_+ \to \mathbb{R}$ be an almost power function of the form \eqref{eq:almost_power}. Then, for all large $y$, $f^{-1}(y)$ is well defined, and as $y \to \infty$,
\[
	f^{-1}(y) = \left(\frac{a^b y}{c(\log{y})^{b}}\right)^{1/a}(1+\smallO{1}).
\]
\end{lemma}

\begin{figure}[t!]
    \centering
    \resizebox{0.4\textwidth}{!}{
        \begin{tikzpicture}[
  yellownode/.style={shape=rectangle, draw=yellow, line width=8},
  bluenode/.style={shape=rectangle, draw=blue, line width=8},
  rednode/.style={shape=rectangle, draw=red, line width=8},
]

\coordinate (O) at (0,0);
\coordinate (xaxis) at (5,0);
\coordinate (yaxis) at (0,8);
\coordinate (xmark1) at (1,0);
\coordinate (ymark1) at (0,2);
\coordinate (intersect) at (1,2);
\coordinate (intersect1) at (5,2);
\coordinate (endpoint1) at (5,8);
\coordinate (endpoint2) at (5,6);
\coordinate (endpoint3) at (5,4);

\path [fill=yellow] (ymark1) -- (intersect) -- (endpoint2) -- (endpoint1) -- (yaxis) -- cycle;
\path [fill=red] (intersect) -- (endpoint2) -- (endpoint3) -- cycle;
\path [fill=red] (intersect) -- (endpoint3) -- (intersect1) -- cycle;
 
\draw[-stealth,thick] (O) -- (xaxis) node[right, align=center] {$a$};
\draw[-stealth,thick] (O) -- (yaxis) node[above, align=center] {$\beta$};

\draw[dashed] (xmark1) node[below, align=center] {$(1,0)$} -- ++(0,8);
\draw[dashed] (ymark1) node[left, align=center] {$(0,2)$} -- ++(5,0);

\draw[blue, thick] (intersect) -- (endpoint2) node[right,align=center] {$\beta=a+1$};
\draw (intersect) -- (endpoint3) node[right,align=center] {$\beta=\frac{a+3}{2}$};

\matrix [draw,below left] at (current bounding box.south east) {
  \node [yellownode,label=right:{$\mathcal{M}^{-1}(k)=\Theta(k)$}] {}; \\
  \node [bluenode,label=right:{$\mathcal{M}^{-1}(k)=\Theta(k/\log{k})$}] {}; \\
  \node [rednode,label=right:{$\mathcal{M}^{-1}(k)=\Theta(w^{1/(2+a-\beta)})$}] {}; \\
};

\end{tikzpicture}
    }
    \caption{Phase diagram of the scaling of $\mathcal{M}^{-1}(k)$ as $k \to \infty$.}
    \label{fig:scaling_M}
\end{figure}
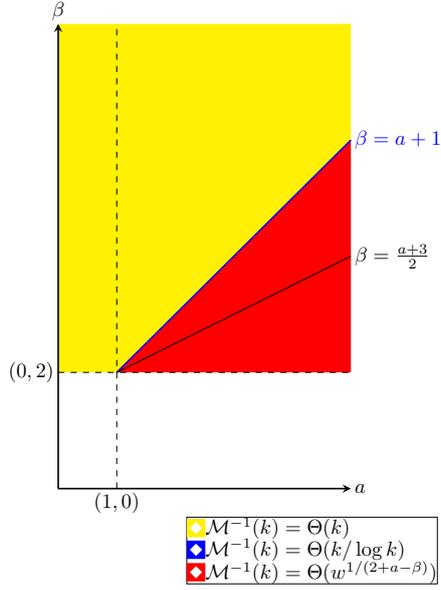

Lemma \ref{lem:log_poly_inv_scaling} is proved in Appendix B. For now we observe that, by~\eqref{eq:M_w_integral_expression} and Lemma~\ref{lem:scaling_M_w}, the function $w \mapsto \mathcal{M}(w)$ satisfies the conditions of the lemma. Hence, we obtain the scaling behavior of the inverse $\mathcal{M}^{-1}(k)$. Figure~\ref{fig:scaling_M} shows the two scaling regimes of $\mathcal{M}^{-1}(k)$ in the $(\beta,a)$ space.

\begin{proposition}\label{prop:scaling_Minv_k}
Let $\mathcal{M}(w)$ be defined as in~\eqref{eq:mean_degree_given_wt} and define
\[
	m'(k) = \begin{cases}
		k			&\text{if } \beta > (a + 1)\vee 2, a \ge 0,\\
		k/\log k	&\text{if } \beta = a+1, a > 1,\\
		k^{1/(2 + a -\beta)}	&\text{if } \beta \in (2,a+1), a > 1.
	\end{cases}
\]
Then,
\[
	\lim_{k \to \infty} \frac{\mathcal{M}^{-1}(k)}{m'(k)} = \begin{cases}
			\frac{\beta - a - 1}{\xi_\alpha}	&\text{if } \beta > (a+1)\vee 2, a \ge 0,\\
			\xi_\alpha^{-1} 				&\text{if } \beta = a + 1, a > 1,\\
			\left( \frac{(\beta - 2)(a + 1 - \beta)}{(a-1)\xi_\alpha}\right)^{1/(2+a-\beta)} &\text{if } \beta \in (2,a+1), a>1.
		\end{cases}
\]
\end{proposition}

Next define
\begin{equation}\label{eq:def_scaling_T_w}
    \sigma(w)=\begin{cases}
       	w &\mbox{if}\;\; \beta>(a+\frac{3}{2}) \vee 2, a \in (0,\infty),\\
        w\log{\left(w \right)} &\mbox{if}\;\; \beta=a+\frac{3}{2}, a\in(\frac{1}{2},\infty),\\  
        w^{4+2a-2\beta} &\mbox{if}\;\; \beta \in (2,a+\frac{3}{2}), a \in(\frac{1}{2},\infty).
        \end{cases}
\end{equation}
Then, combining Proposition~\ref{prop:scaling_Minv_k} with the definition of $\mathrm{S}_k(a,\beta)$ yields the following corollary, which is the first step to proving Theorem~\ref{thm:scaling}:

\begin{corollary}\label{cor:limits_M_inv_k}
As $k \to \infty$,
\[
	 \frac{\sigma(\mathcal{M}^{-1}(k))}{k^2 \mathrm{S}_k(a,\beta)} 
	\to \begin{cases}
		\frac{\beta - a - 1}{\xi_\alpha} &\text{if } \beta > \left(a + \frac{3}{2}\right)\vee 2, a > 0,\\
		\frac{1}{2\xi_\alpha} &\text{if } \beta = a + \frac{3}{2}, a > \frac{1}{2},\\
		\left(\frac{(\beta-a-1)}{\xi_\alpha}\right)^{4+2a-2\beta} &\mbox{if } \beta \in (a+1, a+\frac{3}{2}), a > \frac{1}{2},\\
		\xi_\alpha^{-2} &\mbox{if } \beta=a+1, a > 1,\\  
		\left(\frac{(1+a-\beta)(\beta-2)}{(a-1) \xi_\alpha}\right)^{2} &\mbox{if } {\beta \in (2,a+1), a > 1.}
	\end{cases}
\]
\end{corollary}

\section{Scaling of $T(w)$}\label{sec:T_bounds} 

%

Having established the asymptotic behavior of $\mathcal{M}^{-1}(k)$, in this section we study the scaling of the integral $T(w)$ as $w \to \infty$. In particular, we will show that $T(w) = \bigT{\sigma(w)}$, where $\sigma(w)$ is defined in~\eqref{eq:def_scaling_T_w}. For this, we first obtain upper and lower bounds on $T(w)$ that are sharp, up to constants, for sufficiently large $w$. It is important to note that the bounds we obtain in Lemmas~\ref{lem:N_UB} {and} \ref{lem:N_LB} are valid for any symmetric function $g$ and can thus be used, together with Proposition~\ref{prop:conc_weights}, to study the scaling of the clustering function for different types of $g$'s.

To study the asymptotics of these bounds we will use the specific form of $g$ as in~\eqref{defn:interpolated_wt_fn} and show that for this function $T(w) = \bigO{\sigma(w)}$ and $T(w) = \Omega(\sigma(w))$. We start with the upper bound, as the lower bound requires a bit more computations.

To aid in our analysis we define the spatial integral
\begin{equation}\label{eq:spatial_int}
    \mathcal{S}_w({w}_1,{w}_2):=\int_{(\mathbb{R}^d)^2} \kappa_{\alpha}(\|\mathbf{x}\|,w,{w}_1)\kappa_{\alpha}(\|\mathbf{y}\|,w,{w}_2)\kappa_{\alpha}(\|\mathbf{x}-\mathbf{y}\|,{w}_1,{w}_2) d\mathbf{x}d\mathbf{y},
\end{equation}
so that
\[
	T(w) = \int_{1}^{\infty}\int_{1}^\infty \mathcal{S}_w({w}_1,{w}_2) f_W({w}_1)d{w}_1 f_W({w}_2)d{w}_2.
\]

\subsection{The upper bound}



We obtain an upper bound on $T(w)$ by splitting the integral depending on whether $w_1 \wedge w_2$ is smaller {or larger} than $w$. For each regime we then use a different upper bound for the integrand $\mathcal{S}_w(w_1,w_2)$. Then we study each term separately:

\begin{lemma} \label{lem:N_UB}
Let $T(w)$ be defined as in~\eqref{eq:def_T_w} with $\kappa_{\alpha}$ from~\eqref{eq:conn_fn_form} for some symmetric weight function $g$. Then,
\begin{align*}
    T(w)
    &\leq \xi_\alpha^2 \iint_{1}^{\infty} h(w,w_1,w_2)\ind{w_1\wedge w_2\leq w}(w_1w_2)^{-\beta}dw_1d{w_2} \numberthis \label{eq:T_UB_1} \\
    &\hspace{10pt} + \xi_\alpha^2\iint_{1}^{\infty}g(w,w_1)g(w,w_2) \ind{w_1 \wedge w_2>w}(w_2w_2)^{-\beta}dw_1dw_2, \numberthis \label{eq:T_UB_2}
\end{align*}
where $h(w,w_1, w_2)$ is defined in~\eqref{eq:notn_h_wst}.
\end{lemma}
\begin{proof}
Note that for an upper bound on $T(w)$, it is sufficient to focus on the case $\alpha \in (1,\infty)$, since for any $\alpha \in (1,\infty)$, from \eqref{eq:conn_fn_form} one has $\kappa_{\alpha}(\cdot)\geq \kappa_{\infty}(\cdot)$ point-wise. Hence for the rest of the proof we assume $\alpha \in (1,\infty)$.

Next we begin by providing two upper bounds on $\mathcal{S}_w(w_1,w_2)$. Bounding $\kappa_{\alpha}(\|\mathbf{x}\|,w,w_1)$ from above by $1$, we get that
\begin{align*}
    \mathcal{S}_w(w_1,w_2)\leq \left(\frac{\alpha \omega_d}{\alpha-1}\right)^2  g(w,w_2)g(w_1,w_2).  
\end{align*}
Similarly,
\begin{equation}\label{eq:UB_spatial_asymmetric}
    \mathcal{S}_w(w_1,w_2)\leq \left(\frac{\alpha \omega_d}{\alpha-1}\right)^2  g(w,w_1)g(w_1,w_2).
\end{equation}
Thus using the non-negativity of $g$,
\begin{align*}
    \mathcal{S}_w(w_1,w_2)&\leq \left(\frac{\alpha \omega_d}{\alpha-1}\right)^2  \left[g(w,w_1)\wedge g(w,w_2)\right] g(w_1,w_2)=\left(\frac{\alpha \omega_d}{\alpha-1}\right)^2 h(w,w_1,w_2). \numberthis \label{eq:min_UB_spatial}
\end{align*}

Also note that by bounding $\kappa_{\alpha}(\|\mathbf{x}-\mathbf{y}\|,w_1,w_2) \leq 1$,
\begin{equation}\label{eq:UB_spatial_mixed_leq_1}
    \mathcal{S}_w(w_1,w_2)\leq \left(\frac{\alpha \omega_d}{\alpha-1}\right)^2  g(w,w_2)g(w,w_1).
\end{equation}

Finally, we split the integral
\begin{align*}
    T(w)
    &=(\beta-1)^2\int_{1}^{\infty}\int_{1}^{\infty}\mathcal{S}_w(w_1,w_2)\ind{w_1\wedge w_2 \leq w}(w_1w_2)^{-\beta}dw_1dw_2\\&\hspace{10 pt}+(\beta-1)^2\int_{1}^{\infty}\int_{1}^{\infty}\mathcal{S}_w(w_1,w_2)\ind{w_1\wedge w_2 > w}(w_1w_2)^{-\beta}dw_1dw_2. \numberthis \label{eq:T_decomp_min_l_minv_add_min_g_minv}
\end{align*}
Using the bounds \eqref{eq:min_UB_spatial} and \eqref{eq:UB_spatial_mixed_leq_1}, respectively, for $\mathcal{S}_w(w_1,w_2)$ in the first and second terms of the right hand side (RHS) of \eqref{eq:T_decomp_min_l_minv_add_min_g_minv} gives the result. 
\end{proof}


\begin{remark}[Bounding $\mathcal{S}_w(w_1,w_2)$ differently]
Note that depending on whether $w_1 \wedge w_2$ is larger or smaller than $w$, we bound $\mathcal{S}_w(w_1,w_2)$ differently. As already discussed in Section \ref{ssec:proof_strategy}, {when $w_1$ and $w_2$ have density \eqref{eq:Pareto_weights}, the function $h(w,w_1,w_2)$ is not integrable in $(w_1,w_2)$ if $\beta \leq a+\frac{3}{2}$.} Consequently, using the bound \eqref{eq:min_UB_spatial} for $w_1 \wedge w_2 > w$, for $\beta \leq a+\frac{3}{2}$, would only give us a trivial upper bound of $\infty$ on $T(w)$, which is not useful. Keeping this in mind, we were required to bound $\mathcal{S}_w(w_1,w_2)$ differently. When $\beta>a+\frac{3}{2}$, $h(w,w_1,w_2)$ is integrable when $w_1$ and $w_2$ have density \eqref{eq:Pareto_weights}, and indeed for this case this distinction does not make a difference: as we will see, when $\beta>a+\frac{3}{2}$, the term \eqref{eq:T_UB_2} is of negligible order compared to the term \eqref{eq:T_UB_1}. \end{remark}

With this bound in hands we can now use the definition of our interpolation function $g$ to study the scaling of $T(w)$.

\begin{lemma}\label{lem:scaling_T_w_UB}
Let $T(w)$ be defined as in~\eqref{eq:def_T_w} with $\kappa$ and $g$ as in~\eqref{eq:conn_fn_form} and~\eqref{defn:interpolated_wt_fn}, respectively. Then, as $w \to \infty$,
\begin{equation}\label{eq:scaling_T_w_UB}
	T(w) =
	\begin{cases}
		\bigO{w} &\mbox{if}\;\; \beta> \left(a + \frac{3}{2}\right) \vee 2, a >0,\\
		\bigO{w\log{\left(w \right)}} &\mbox{if}\;\; \beta=a+\frac{3}{2}, a > \frac{1}{2},\\ 
		\bigO{(w)^{4+2a-2\beta}}  &\mbox{if}\;\; \beta \in (2,a+\frac{3}{2}), a > \frac{1}{2}. 
	\end{cases}
\end{equation}
\end{lemma}

\begin{proof}
We show that the first term~\eqref{eq:T_UB_1} from Lemma~\ref{lem:N_UB} has the exact same scaling as in the statement of the lemma. For the second term~\eqref{eq:T_UB_2} we then use the definition of $g$ from \eqref{defn:interpolated_wt_fn} to get, as $w \to \infty$,
\[
    \int_w^{\infty}\int_{w}^{\infty} g(w,w_1)g(w,w_2)f_W(\underline{\mathbf{w}})d\mathbf{w}
    =\left(w \right)^{2a}\left((\beta-1)\int_{w}^{\infty} w_1^{1-\beta}dw_1 \right)^2=\bigO{\left(w\right)^{4+2a-2\beta}}.
\]

Noting that $w^{4+2a-2\beta}$ is bounded by each term on the RHS of~\eqref{eq:scaling_T_w_UB}, in its respective case, then finishes the proof.

For~\eqref{eq:T_UB_1}, we again use the definition of $g$ from \eqref{defn:interpolated_wt_fn} to get
\begin{align*}
    &\hspace{-20pt}\int_1^{\infty}\int_{1}^{\infty} h(w,w_1,w_2)
    	\ind{w_1\wedge w_2\leq w}(w_1w_2)^{-\beta}dw_1dw_2\\
    &=\int_1^{\infty}\int_{1}^{\infty} \left[g(w,w_1)\wedge g(w,w_2)\right] g(w_1,w_2)\ind{w_1\wedge w_2\leq w}(w_1w_2)^{-\beta}dw_1dw_2\\
    &=\int_1^{\infty}\int_1^{\infty} w(w_1 \wedge w_2)^a (w_1 \vee w_2)(w_1 \wedge w_2)^a\ind{w_1 \wedge w_2\leq w}(w_1w_2)^{-\beta}dw_1dw_2\\
    &=2w\int_1^{w} \int_{w_1}^{\infty}{w_1}^{2a-\beta}{w_2}^{1-\beta}dw_2dw_1\\
    &=\frac{2}{\beta-2}w\int_1^{w} w_1^{2+2a-2\beta}dw_1,
\end{align*}
where to obtain the second last inequality we have used the symmetry in $w_1$ and $w_2$. Now, the intended result follows by noting that as $w \to \infty$,
\begin{align*}
    \int_1^{w} w_1^{2+2a-2\beta}dw_1=
    	\begin{cases}
        \bigO{1} &\mbox{if}\;\; \beta> \left(a + \frac{3}{2}\right) \vee 2, a >0,\\
        \bigO{\log{\left(w \right)}} &\mbox{if}\;\; \beta=a+\frac{3}{2}, a > \frac{1}{2},\\ 
        \bigO{(w)^{3+2a-2\beta}}  &\mbox{if}\;\; \beta \in(2,a+\frac{3}{2}), a > \frac{1}{2}.
        \end{cases}
\end{align*}
\end{proof}

\subsection{The lower bound}

Next, we prove a lower bound on $T(w)$. For this we rely on geometric techniques. Recall that $\omega_d$ denotes the Lebesgue measure of a ball of unit radius in $\mathbb{R}^d$.

\begin{lemma}\label{lem:N_LB}
Consider the connection function \eqref{eq:conn_fn_form}, with some symmetric weight function $g$. Recall $h(w,w_1,w_2)$ from \eqref{eq:notn_h_wst}. Then,
\begin{align*}
    &T(w) \geq \omega_d^2\int_{1}^{\infty} \int_{1}^{\infty} \ind{g(w_1,w_2)<g(w,w_2)} h(w,w_1,w_2)\\
    &\hspace{100 pt} \times \left[1-\left(\frac{g(w_1,w_2)}{(g(w,w_1) \wedge g(w,w_2)}\right)^{1/d} \right]^d f_W(w_1)dw_1 f_W(w_2)dw_2. \numberthis \label{eq:T_LB}
\end{align*}
\end{lemma}

\begin{proof}
Note that {for a lower bound on $T(w)$,} we need only consider the case $\alpha= \infty$ since for any $\alpha \in (0,\infty)$, $\kappa_\alpha > \kappa_{\infty}$ point-wise. Thus, we assume $\alpha=\infty$ for the rest of the proof. For $\mathbf{x} \in \mathbb{R}^d$ and $r>0$, we also use the notation 
\[
B(\mathbf{x},r):= \{\mathbf{y} \in \mathbb{R}^d: \|\mathbf{x}-\mathbf{y}\|<r \}
\]
to denote the open ball of radius $r$ around $\mathbf{x}$ in $\mathbb{R}^d$.

{For $\alpha=\infty$,}
\begin{align*}
    & T(w) =  \int_{1}^{\infty} \int_{1}^{\infty} \int_{(\mathbb{R}^d)^2} \ind{\|\mathbf{x}\|<g(w,w_1)^{1/d}}\ind{\|\mathbf{y}\|<g(w,w_2)^{1/d}}\\&\hspace{118 pt}\times \ind{\|\mathbf{x}-\mathbf{y}\|<g(w_1,w_2)^{1/d}} d\mathbf{x}d\mathbf{y}\;f_W(w_1)dw_1 f_W(w_2)dw_2,
\end{align*}
where $g$ is the weight function \eqref{defn:interpolated_wt_fn}.
A simple change of variable $\mathbf{z}=\mathbf{x}-\mathbf{y}$ yields
\begin{align*}
    & T(w) =  \int_{1}^{\infty} \int_{1}^{\infty} \int_{B(\mathbf{0},g(w,w_1)^{1/d})}\int_{B(\mathbf{x},g(w,w_2)^{1/d})} \ind{\|\mathbf{z}\|<g(w_1,w_2)^{1/d}} d\mathbf{z}d\mathbf{x}\;f_W(w_1)dw_1 f_W(w_2)dw_2.
\end{align*}

Note that if $\|\mathbf{x}\|>g(w,w_1)^{1/d}+g(w_1,w_2)^{1/d}$, then the set $B(\mathbf{x},g(w,w_2)^{1/d})$ is disjoint from the open ball of radius $g(w_1,w_2)^{1/d}$ around the origin $\mathbf{0} \in \mathbb{R}^d$, see Figure \ref{LB_1_fig}.

This means in particular that if $\|\mathbf{z}\|<g(w_1,w_2)^{1/d}$, then $\mathbf{z} \notin B(\mathbf{x},g(w,w_2)^{1/d})$.  Consequently in this case,
\begin{equation}\label{eq:two_circle_inner_integral}
    \int_{B(\mathbf{x},g(w,w_2)^{1/d})} \ind{\|\mathbf{z}\|<g(w_1,w_2)^{1/d}} d\mathbf{z}=0.
\end{equation} 
By this observation,
\begin{align*}
    & T(w) =  \int_{1}^{\infty} \int_{1}^{\infty} \int_{B\left(\mathbf{0},(g(w,w_1)^{1/d}\wedge [g(w,w_2)^{1/d}+g(w_1,w_2)^{1/d}])\right)} \\ & \hspace{110 pt}\times \int_{B(\mathbf{x},g(w,w_2)^{1/d})} \ind{\|\mathbf{z}\|<g(w_1,w_2)^{1/d}} d\mathbf{z}d\mathbf{x}\;f_W(w_1)dw_1 f_W(w_2)dw_2.
\end{align*}

Next, we obtain a generic lower bound on the left hand side (LHS) of \eqref{eq:two_circle_inner_integral}, in the complementary range, that is, for $\|\mathbf{x}\|<g(w,w_2)^{1/d}+g(w_1,w_2)^{1/d}$.

First observe that the integral \eqref{eq:two_circle_inner_integral} is the volume of intersection of the two balls $B(\mathbf{x},g(w,w_2)^{1/d})$ and $B(\mathbf{0},g(w_1,w_2)^{1/d})$. A lower bound on this volume is the volume of any ball that is completely contained in this intersection. It is not too hard to figure out what {the maximum diameter $\ell$} of such a ball, completely contained in this intersection, can be (see Figure \ref{LB_fig}):

\begin{figure}[h]
     \centering
     \begin{subfigure}[t]{0.50\textwidth}
         \centering
   \resizebox{1\textwidth}{!}{
        \begin{tikzpicture}[scale=0.8]

\coordinate (O) at (0,0);
\coordinate (x) at (7,3);
\coordinate (x2) at (12,3);
\coordinate (l) at (0.8,0.6);
\coordinate (l2) at (1.6,1.2);
\coordinate (O2) at (-1,0);

\draw[thick] (O) circle (2);
\draw[thick] (x) circle (5);

\node [right,align=center] at (O) {$\mathbf{0}$} ;

\draw[thick] (O) -- (x) node[above,align=center] {$\mathbf{x}$};
\draw[stealth-stealth,blue,thick] (x) -- node[sloped,above,align=center] {$g(w,{w}_2)^{1/d}$} (x2);
\draw[stealth-stealth,green,thick] (O) -- node[sloped,below, align=center] {$g({w}_1,{w}_2)^{1/d}$} ++(-2,0);
\end{tikzpicture}
    }

    \caption{When $\|\mathbf{x}\|>g(w,w_1)^{1/d}+g(w_1,w_2)^{1/d}$, $B(\mathbf{0},g(w,w_2)^{1/d})$ is disjoint from $B(\mathbf{x},g(w,w_2)^{1/d})$, so that the integral in the LHS of \eqref{eq:two_circle_inner_integral} equals $0$.}
    \label{LB_1_fig}
     \end{subfigure}
     \hfill
     \begin{subfigure}[t]{0.45\textwidth}
         \centering
    \resizebox{1\textwidth}{!}{
        \begin{tikzpicture}[scale=0.8]

\coordinate (O) at (0,0);
\coordinate (x) at (4,3);
\coordinate (x2) at (9,3);
\coordinate (l) at (0.8,0.6);
\coordinate (l2) at (1.6,1.2);
\coordinate (O2) at (-1,0);

\draw[thick] (O) circle (2);
\draw[thick] (x) circle (5);
\draw[thick,red] (l) circle (1);

\node [right,align=center] at (O) {$\mathbf{0}$} ;

\draw[thick] (O) -- (x) node[above,align=center] {$\mathbf{x}$};
\draw[stealth-stealth,blue,thick] (x) -- node[sloped,above,align=center] {$g(w,{w}_2)^{1/d}$} (x2);
\draw[stealth-stealth,red,thick] (O) -- node[sloped,below, align=center] {$\ell$} (l2);
\draw[stealth-stealth,green,thick] (O) -- node[sloped,below, align=center] {$g({w}_1,{w}_2)^{1/d}$} ++(-2,0);
\end{tikzpicture}
    }

    \caption{When $\|\mathbf{x}\|<g(w,w_1)^{1/d}+g(w_1,w_2)^{1/d}$, $B(\mathbf{0},g(w,w_2)^{1/d})\cap B(\mathbf{x},g(w,w_2)^{1/d})\neq \emptyset$. The maximum possible diameter $\ell$ of a ball completely contained in $B(\mathbf{x},g(w,w_2)^{1/d})\cap B(\mathbf{0},g(w_1,w_2)^{1/d})$ is shown in red.}
    \label{LB_fig}
     \end{subfigure}
        \caption{}
        \label{fig:three graphs}
\end{figure}

It is
\begin{align*}
    \ell = \begin{cases}
		0 &\text{if } B(\mathbf{0},g(w_1,w_2)^{1/d}) \cap B(\mathbf{x},g(w,w_2)^{1/d}) = \emptyset,\\
		2g(w_1,w_2)^{1/d} &\text{if } B(\mathbf{0},g(w_1,w_2)^{1/d}) \subset B(\mathbf{x},g(w,w_2)^{1/d}), \\
		2g(w,w_2)^{1/d} &\text{if } B(\mathbf{0},g(w_1,w_2)^{1/d}) \supset B(\mathbf{x},g(w,w_2)^{1/d}),\\
		g(w,w_2)^{1/d}+g(w_1,w_2)^{1/d}-\|\mathbf{x}\| &\text{otherwise.}
	\end{cases}
\end{align*}

In particular, when $g(w,w_2)>g(w_1,w_2)$, since it can never hold that $B(\mathbf{0},g(w_1,w_2)^{1/d}) \supset B(\mathbf{x},g(w,w_2)^{1/d})$, 
\begin{align*}
   &\int_{B(\mathbf{x},g(w,w_2)^{1/d})} \ind{\|\mathbf{z}\|<g(w_1,w_2)^{1/d}} d\mathbf{z} \\& \geq \omega_d \left[\frac{g(w,w_2)^{1/d}+g(w_1,w_2)^{1/d}-\|\mathbf{x}\|}{2} \wedge  g(w_1,w_2)^{1/d}\right]^d \ind{g(w,w_2)>g(w_1,w_2)}.
\end{align*}

Further bounding from below by multiplying by the indicator $\ind{\|\mathbf{x}\|<g(w,w_2)^{1/d}-g(w_1,w_2)^{1/d}}$,
\begin{align*}
    &\int_{B(\mathbf{x},g(w,w_2)^{1/d})} \ind{\|\mathbf{z}\|<g(w_1,w_2)^{1/d}} d\mathbf{z}\\& \geq \omega_d  \ind{g(w_1,w_2)<g(w,w_2)} \ind{\|\mathbf{x}\|<g(w,w_2)^{1/d}-g(w_1,w_2)^{1/d}} g(w_1,w_2). \numberthis \label{eq:final_LB_inner_two_circle_int}
\end{align*}

At this point we introduce some short-hand notations to keep things neat.
Let us denote
\begin{equation}\label{eq:R_1_lb_lem}
    \mathcal{R}_\pm = g(w,w_1)^{1/d}\wedge[g(w,w_2)^{1/d} \pm g(w_1,w_2)^{1/d}].
\end{equation}

Using the lower bound \eqref{eq:final_LB_inner_two_circle_int} and using the notation \eqref{eq:R_1_lb_lem} {(recall the notation $f_{W}(\mathbf{w})d\mathbf{w}$ from \eqref{eq:notn_two_dim_weight_integral})}, we have 
{\begin{align*}
& T(w)
    \\& \geq   \int \int\int_{B(\mathbf{0},\mathcal{R}_+)}\left(\omega_d  \ind{g(w_1,w_2)<g(w,w_2)} \ind{\|\mathbf{x}\|<g(w,w_2)^{1/d}-g(w_1,w_2)^{1/d}} g(w_1,w_2)\right)  d\mathbf{x} f_{W}(\mathbf{w})d\mathbf{w}
    \\&= \omega_d  \int \int g(w_1,w_2) \ind{g(w_1,w_2)<g(w,w_2)}
     \left(\int_{B(\mathbf{0},\mathcal{R}_-)} d\mathbf{x}\right) f_{W}(\mathbf{w})d\mathbf{w} 
    \\&= \omega_d^2\int \int g(w_1,w_2) \ind{g(w_1,w_2)<g(w,w_2)} \left[g(w,w_1)^{1/d} \wedge [g(w,w_2)^{1/d}-g(w_1,w_2)^{1/d}] \right]^d f_{W}(\mathbf{w})d\mathbf{w}.\numberthis \label{eq:gen_geom_lb_almost}
\end{align*}}
Finally, noting that $g(w,w_1)^{1/d}> g(w,w_1)^{1/d}-g(w_1,w_2)^{1/d}$, we have from \eqref{eq:gen_geom_lb_almost},
{\begin{align*}
        & T(w)
        \\& \geq \omega_d^2\int \int g(w_1,w_2) \ind{g(w_1,w_2)<g(w,w_2)} \left[g(w,w_1)^{1/d} \wedge [g(w,w_2)^{1/d}-g(w_1,w_2)^{1/d}] \right]^d f_{W}(\mathbf{w})d\mathbf{w}\\
        & \geq \omega_d^2\int \int g(w_1,w_2) \ind{g(w_1,w_2)<g(w,w_2)} \\& \hspace{40 pt}\times\left[[g(w,w_1)^{1/d}-g(w_1,w_2)^{1/d}] \wedge [g(w,w_2)^{1/d}-g(w_1,w_2)^{1/d}] \right]^d f_{W}(\mathbf{w})d\mathbf{w}\\
        &=\omega_d^2\int \int g(w_1,w_2) \ind{g(w_1,w_2)<g(w,w_2)} \left[[g(w,w_1)^{1/d} \wedge g(w,w_2)^{1/d}]-g(w_1,w_2)^{1/d} \right]^d f_{W}(\mathbf{w})d\mathbf{w}\\
        &=\omega_d^2\int \int g(w_1,w_2) \ind{g(w_1,w_2)<g(w,w_2)} \\& \hspace{40 pt} \times\left(g(w,w_1) \wedge g(w,w_2) \right) \left[1-\left(\frac{g(w_1,w_2)}{(g(w,w_1) \wedge g(w,w_2)}\right)^{1/d} \right]^d f_W(w_1)dw_1f_W(w_2)dw_2\\
        &=\omega_d^2\int \int \ind{g(w_1,w_2)<g(w,w_2)} h(w,w_1,w_2) \left[1-\left(\frac{g(w_1,w_2)}{(g(w,w_1) \wedge g(w,w_2)}\right)^{1/d} \right]^d f_{W}(\mathbf{w})d\mathbf{w},
\end{align*}}
which is \eqref{eq:T_LB} recalling \eqref{eq:Pareto_weights}.
\end{proof}

\begin{lemma}\label{lem:scaling_T_w_LB}
Let $T(w)$ be defined as in~\eqref{eq:def_T_w} with $\kappa$ and $g$ as in~\eqref{eq:conn_fn_form} and~\eqref{defn:interpolated_wt_fn}, respectively. Then, as $w \to \infty$,
\[
	T(w) =
	\begin{cases}
		\Omega\left(w\right) &\mbox{if}\;\; \beta> \left(a + \frac{3}{2}\right) \vee 2, a >0,\\
		\Omega\left(w\log{\left(w \right)}\right) &\mbox{if}\;\; \beta=a+\frac{3}{2}, a > \frac{1}{2},\\ 
		\Omega\left(w^{4+2a-2\beta}\right)  &\mbox{if}\;\; \beta \in (2,a+\frac{3}{2}), a > \frac{1}{2}. 
	\end{cases}
\]
\end{lemma}

\begin{proof}
Using the lower bound from Lemma~\ref{lem:N_LB}, it suffices to show that as $w \to \infty$,
\begin{align*}
    &\hspace{-20pt}\omega_d^2\int_1^{\infty}\int_{1}^{\infty} \ind{g(w_1,w_2)<g(w,w_2)}h(w,w_1,w_2)
    	\left[1-\left(\frac{g(w_1,w_2)}{(g(w,w_1) \wedge g(w,w_2)}\right)^{1/d} \right]^d f_W(\underline{\mathbf{w}})d\underline{\mathbf{w}}\\
    &=\begin{cases}
        \Omega\left(w\right) &\mbox{if}\;\; \beta>(a+\frac{3}{2}) \vee 2, a \in (0,\infty),\\
        \Omega\left(w\log{\left(w \right)}\right) &\mbox{if}\;\; \beta=a+\frac{3}{2}, a\in(1/2,\infty),\\  \Omega\left(w^{4+2a-2\beta}\right) &\mbox{if}\;\; \beta \in(2,a+\frac{3}{2}), a \in(1/2,\infty). \numberthis \label{eq:scaling_LB} 
        \end{cases}
\end{align*}
Observe that if both $w_1,w_2 < {\mathcal{M}^{-1}(k)}/{2}$, then {using the definition of $g$ from \eqref{defn:interpolated_wt_fn}, one has $\ind{g(w_1,w_2)<g(\mathcal{M}^{-1}(k),w_2)}=1$. Furthermore, \begin{align*}
    \left[1-\left(\frac{g(w_1,w_2)}{(g(\mathcal{M}^{-1}(k),w_1) \wedge g(\mathcal{M}^{-1}(k),w_2)}\right)^{1/d} \right]^d
    &=\left[1-\left(\frac{(w_1 \vee w_2)(w_1 \wedge w_2)^a}{\mathcal{M}^{-1}(k) (w_1\wedge \ch{w_2})^a}\right)^{1/d} \right]^d \\
    &\geq \left[1-\left(\frac{1}{2}\right)^{1/d}\right]^d =:\underline{K}.
\end{align*}}

Consequently, we can bound the LHS of \eqref{eq:scaling_LB} from below by restricting the integration domain from $(1,\infty)\times (1,\infty)$ to $\left(1,{\mathcal{M}^{-1}(k)}/{2}\right)\times \left(1,{\mathcal{M}^{-1}(k)}/{2}\right)$, to get the following lower bound on it:
\begin{equation}\label{eq:double_int_T_LB_scaling}
    \underline{K}\omega_d^2(\beta-1)^2\int_{1}^{\frac{\mathcal{M}^{-1}(k)}{2}} \int_{1}^{\frac{\mathcal{M}^{-1}(k)}{2}}  h(\mathcal{M}^{-1}(k),w_1,w_2) (w_1w_2)^{-\beta}dw_1dw_2.
\end{equation}
Hence to conclude Lemma~\ref{lem:scaling_T_w_LB}, it is enough to check that the double integral in \eqref{eq:double_int_T_LB_scaling} is at least the RHS of \eqref{eq:scaling_LB}, which we proceed to do. Note that
\begin{align*}
    &\int_{1}^{\frac{\mathcal{M}^{-1}(k)}{2}} \int_{1}^{\frac{\mathcal{M}^{-1}(k)}{2}}  h(\mathcal{M}^{-1}(k),w_1,w_2) (w_1w_2)^{-\beta}dw_1dw_2\\
    &=\int_{1}^{\frac{\mathcal{M}^{-1}(k)}{2}} \int_{1}^{\frac{\mathcal{M}^{-1}(k)}{2}}  \left[g(\mathcal{M}^{-1}(k),w_1)\wedge g(\mathcal{M}^{-1}(k),w_2)\right] g(w_1,w_2) (w_1w_2)^{-\beta}dw_1dw_2\\
    &= \mathcal{M}^{-1}(k)\int_{1}^{\frac{\mathcal{M}^{-1}(k)}{2}} \int_{1}^{\frac{\mathcal{M}^{-1}(k)}{2}}(w_1 \wedge w_2)^a (w_1 \vee w_2)(w_1 \wedge w_2)^a(w_1w_2)^{-\beta}dw_1dw_2\\
    &=2\mathcal{M}^{-1}(k)\int_{1}^{\frac{\mathcal{M}^{-1}(k)}{2}} \int_{w_1}^{\frac{\mathcal{M}^{-1}(k)}{2}}w_1^{2a-\beta}w_2^{1-\beta}{dw_2dw_1}\\
    &=2{K}\underline{K}\mathcal{M}^{-1}(k)\int_1^{\frac{\mathcal{M}^{-1}(k)}{2}}w_1^{2a+2-2\beta}dw_1\\&\hspace{10 pt}-2{K}\underline{K}\mathcal{M}^{-1}(k)\int_1^{\frac{\mathcal{M}^{-1}(k)}{2}}w_1^{2a-\beta}(\mathcal{M}^{-1}(k)/2)^{2-\beta}dw_1, \numberthis \label{eq:diff_of_two_terms}
\end{align*}
where $K={1}/{(\beta-2)}>0$, and in the second last step we have used the symmetry in $w_1$ and $w_2$.

{To conclude Lemma~\ref{lem:scaling_T_w_LB}, it is enough to check that \eqref{eq:diff_of_two_terms} is at least the RHS of \eqref{eq:scaling_LB}. {To keep notation clear, let us write $R$ for the constant $K\underline{K}$.} We do a case-by-case analysis:
\paragraph{Case 1: $\beta \geq (a+\frac{3}{2})\vee 2$.} In this case, observe that the first term of \eqref{eq:diff_of_two_terms} equals
   \begin{align*}
    2{R}\mathcal{M}^{-1}(k)\times\begin{cases}
        \Omega\left( 1\right) &\mbox{if}\;\; \beta>(a+3/2) \vee 2, a \in (0,\infty),\\
        \Omega\left(\log{\left(\mathcal{M}^{-1}(k) \right)}\right) &\mbox{if}\;\; \beta=a+3/2, a\in(1/2,\infty),\end{cases}
\end{align*} 
while the subtracted second term of \eqref{eq:diff_of_two_terms} equals 
\begin{align*}
    2{R}\mathcal{M}^{-1}(k)\int_1^{\frac{\mathcal{M}^{-1}(k)}{2}}w_1^{2a-\beta}(\mathcal{M}^{-1}(k)/2)^{2-\beta}dw_1=\Theta((\mathcal{M}^{-1}(k))^{4+2a-2\beta}).
\end{align*}
In particular, since $\beta>a+\frac{3}{2}$, the first term of \eqref{eq:diff_of_two_terms} is the dominating term, and is indeed at least the RHS of \eqref{eq:scaling_LB}.

\paragraph{Case 2: $2<\beta<a+\frac{3}{2}$.} In this case, we note that the first term of \eqref{eq:diff_of_two_terms} equals
\begin{align*} \frac{2{R}}{3+2a-2\beta}\mathcal{M}^{-1}(k)\left[\left(\frac{\mathcal{M}^{-1}(k)}{2} \right)^{3+2a-2\beta}-1\right]. \numberthis \label{eq:first_term_diff}
\end{align*}
Simplifying the subtracted second term in \eqref{eq:diff_of_two_terms} gives
{\begin{align*}
    &2^{\beta-1}R\mathcal{M}^{-1}(k)^{3-\beta}\int_{1}^{\frac{\mathcal{M}^{-1}(k)}{2}}w_1^{2a-\beta}dw_1\leq2^{\beta-1}R\mathcal{M}^{-1}(k)^{3-\beta}
        \frac{1}{1+2a-\beta}\left(\frac{\mathcal{M}^{-1}(k)}{2}\right)^{1+2a-\beta}, \numberthis \label{eq:second_term_diff}
\end{align*}}
\ch{since $\beta<a+3/2<1+2a$ i.e. $1+2a-\beta>0$, as $a>1/2$ (note that we need $a>1/2$ so that $a+3/2>2$).} Using \eqref{eq:first_term_diff} and \eqref{eq:second_term_diff}, we note that the term \eqref{eq:diff_of_two_terms} is at least
\begin{align*}
    &\frac{2{R}}{3+2a-2\beta}\mathcal{M}^{-1}(k)\left[\left(\frac{\mathcal{M}^{-1}(k)}{2} \right)^{3+2a-2\beta}-1\right]-\frac{2^{\beta-1}{R}}{1+2a-\beta}\mathcal{M}^{-1}(k)^{3-\beta}\left(\frac{\mathcal{M}^{-1}(k)}{2} \right)^{1+2a-\beta}\\&=\frac{2^{2\beta-2a-2}{R}}{3+2a-2\beta}\mathcal{M}^{-1}(k)\left[\mathcal{M}^{-1}(k)^{3+2a-2\beta}-2^{3+2a-2\beta}\right]-\frac{2^{2\beta-2a-2}{R}}{1+2a-\beta}\mathcal{M}^{-1}(k)^{4+2a-2\beta}\\
    &=2^{2\beta-2a-2}{R}\mathcal{M}^{-1}(k)^{4+2a-2\beta}\left[\frac{1}{3+2a-2\beta}-\frac{1}{1+2a-\beta} \right]+\Omega(\mathcal{M}^{-1}(k)).
\end{align*}
 Observing that $\left[{(3+2a-2\beta)^{-1}}-{(1+2a-\beta)^{-1}} \right]>0$ since $\beta>2$, and $4+2a-2\beta>1$ since $\beta<a+\frac{3}{2}$, the last expression is indeed $\Omega(\mathcal{M}^{-1}(k)^{4+2a-2\beta})$, which completes the proof for the case $1+2a-\beta>0$. }

\end{proof}

Note that combining Lemmas \ref{lem:scaling_T_w_UB} and \ref{lem:scaling_T_w_LB} we get the scaling of $T(w)$, as $w \to \infty$:
\begin{equation}\label{eq:scaling_T_in_m_inv}
    T(w)=\begin{cases}
        \Theta\left(w \right) &\mbox{if}\;\; \beta>(a+3/2) \vee 2, a \in (0,\infty)\\
        \Theta\left(w\log{w}\right) &\mbox{if}\;\; \beta=a+3/2, a\in(1/2,\infty)\\  
        \Theta\left(w^{4+2a-2\beta} \right) &\mbox{if}\;\; \beta \in (2,a+3/2), a \in(1/2,\infty).\end{cases}  
\end{equation}

{At this point we make the observation that the assumptions of Proposition \ref{prop:conc_weights} are indeed satisfied. Indeed, assumption (1.) is immediate from the choice of $g$ in \eqref{defn:interpolated_wt_fn}. It is clear from the expression \eqref{eq:M_w_integral_expression} that assumption (2.) is satisfied. Assumption (3.) follows from Proposition \ref{prop:scaling_Minv_k}. Finally, assumption (4.) follows from \eqref{eq:scaling_T_in_m_inv}. We make no more mention of this, and we apply Proposition \ref{prop:conc_weights} in the next section where we prove Theorem \ref{thm:scaling}, wherever necessary.}

\section{Proof of Theorem \ref{thm:scaling}}\label{sec:main_thm_proof}

In this section we prove Theorem \ref{thm:scaling}. Recall the definition of $T(w)$ from \eqref{eq:def_T_w}, and note that using the notation $\mathcal{S}_w(w_1,w_2)$ from \eqref{eq:spatial_int}, we have
\begin{equation}\label{eq:T_in_spatial}
    T(\mathcal{M}^{-1}(k))= \int_{1}^{\infty} \int_{1}^{\infty} \mathcal{S}_{\mathcal{M}^{-1}(k)}(w_1,w_2) f_W(w_1)dw_1 f_W(w_2)dw_2
\end{equation}

Observe that the sequence $(\mathrm{S}_k(a,\beta))_{k \ge 1}$ is fluctuation-robust, and thus so is $(k^2 \mathrm{S}_k(a,\beta))_{k \ge 1}$. Hence by Proposition~\ref{prop:conc_weights} it suffices to show that
\[
	\lim_{k \to \infty} \frac{T(\mathcal{M}^{-1}(k))}{k^2 \mathrm{S}_k(a,\beta)} = \Gamma(a,\alpha,\beta,d).
\]

Recall the definition of $\sigma$ from~\eqref{eq:def_scaling_T_w}. Following the strategy we laid out, we write
\[
	\frac{T(\mathcal{M}^{-1}(k))}{k^2 \mathrm{S}_k(a,\beta)} 
	= \frac{\sigma(\mathcal{M}^{-1}(k))}{k^2 \mathrm{S}_k(a,\beta)} \frac{T(\mathcal{M}^{-1}(k))}{\sigma(\mathcal{M}^{-1}(k))}.
\]
The first factor converges by Corollary~\ref{cor:limits_M_inv_k}. Comparing the constants in Corollary~\ref{cor:limits_M_inv_k} to $\Gamma(a,\alpha,\beta,d)$ in Theorem~\ref{thm:scaling}, it suffices to show that
\begin{equation}\label{eq:lim_t_over_sigma}
	\lim_{w \to \infty} \frac{T(w)}{\sigma(w)} = \begin{cases}
		\int_1^{\infty}\int_{1}^\infty I(w_1, w_2) f_W(w_1)dw_1 f_W(w_2)w_1
		&\text{if } \beta > \left(a + \frac{3}{2}\right) \vee 2, a > 0,\\
		\left(a + \frac{1}{2}\right) \int_{1}^\infty I(1, r) f_W(r)dr
		&\text{if } \beta = a + \frac{3}{2}, a > \frac{1}{2},\\
		\int_0^{\infty}\int_{0}^\infty \mathcal{S}_1(w_1, w_2) f_W(w_1)w_1 f_W(w_2)w_2
		&\text{if } \beta \in (2 , a + \frac{3}{2}), a > \frac{1}{2}.
	\end{cases}
\end{equation}
This is achieved in Lemmas~\ref{lem:conv_T_case_1}, \ref{lem:conv_T_case_2} and~\ref{lem:conv_T_case_3}.

The general strategy is to employ a specific change of variables $(x,y) = (\phi_1(x^\prime), \phi_2(y^\prime))$ and $(w_1,w_2) = (\psi_1(w_1^\prime), \psi_2(w_2^\prime))$ to the integral expression of $T(w)$ so that 
\[
	T(w) = k^2 \mathrm{S}_k(a,\beta) \int_1^{\infty}\int_1^{\infty} \hat{\mathcal{S}}_{w}(\psi_1(w_1^\prime), \psi_2(w_2^\prime)) d\psi_1(w_1^\prime) d\psi_2(w_2^\prime),
\]
where $\hat{\mathcal{S}}_{w}(\psi_1(w_1^\prime), \psi_2(w_2^\prime))$ is the integral $\mathcal{S}_{w}(w_1,w_2)$ after the change of variables and without the leading order scaling. We then analyze the new integral and use dominated convergence to prove that it converges to a constant as $k \to \infty$. This is also where we make use of the upper bounds established for $T(w)$ in Section~\ref{sec:T_bounds}.

The right choice of the change of variables depends on the regime of $a$ and $\beta$. In addition, for the third regime: $2 < \beta < (a+\frac{3}{2})$, the entire integration domain $w_1, w_2 \in (1,\infty)$ in $T(w)$ contributes to the constant. Instead in the other two regimes the main contribution comes from the integration over the domain $w_1, w_2 \le w$, which is summarized in the following result:

\begin{lemma}\label{lem:T_M_k_main_term}
Let $a > 0$ and $\beta > (a+\frac{3}{2}) \vee 2$ or $\beta = a + \frac{3}{2} > 2$. Then, as $w \to \infty$,
\[
	T(w) = (1+o(1)\int_1^w\int_1^{w} \mathcal{S}_{w}(w_1,w_2) 
	f_W(w_1)dw_1 f_W(w_2)dw_1.
\]
\end{lemma}

\begin{proof}
Lemmas~\ref{lem:scaling_T_w_UB} and \ref{lem:scaling_T_w_LB} together imply that in these regimes, $T(w)$ scales either as $w$ or $w\log(w)$. Hence, it suffices to show that
\[
	\left|T(w) - \int_1^w\int_1^{w} \mathcal{S}_{w}(w_1,w_2) f_W(w_1)dw_1 f_W(w_2)dw_2\right| = \smallO{w \log(w)}.
\]
Note that the LHS is nothing more than the integral over the domain $w_1 \vee w_2> w$. We must therefore show that
\[
	\int_1^{\infty}\int_{1}^{\infty} \mathcal{S}_{w}(w_1,w_2) \ind{w_1 \vee w_2> w} f_W(w_1)dw_1 f_W(w_2)dw_2 = \smallO{w \log(w)}.
\]


We first note that by symmetry the full integral is equal to twice the integral where $w_1 \le w_2$. We then split the inner integral over the two cases $1<w_1<w$ and $w <w_1 < w_2$ to obtain:
\begin{align*}
	&\hspace{-40pt}\int_1^{\infty}\int_{1}^{\infty} \mathcal{S}_{w}(w_1,w_2) \ind{w_1 \vee w_2> w} f_W(w_1)dw_1 f_W(w_2)dw_2\\
	&\le 2 \int_w^\infty \int_1^w \mathcal{S}_w(w_1,w_2) f_W(w_1)dw_1 f_W(w_2)dw_2 \numberthis \label{eq:int_large_weights_1} \\
	&\hspace{10pt}+ 2\int_w^\infty \int_w^{w_2} \mathcal{S}_w(w_1,w_2) f_W(w_1)dw_1 f_W(w_2)dw_2. \numberthis \label{eq:int_large_weights_2} 
\end{align*}

We first show that the second integral in~\eqref{eq:int_large_weights_2} is $\bigO{w^{4 + 2a-2\beta}}$ which is $o(w \log(w))$ for the two regimes of $a$ and $\beta$. Using the bound~\eqref{eq:UB_spatial_mixed_leq_1} we get
\begin{align*}
	&2\int_w^\infty \int_w^{w_2} \mathcal{S}_w(w_1,w_2) f_W(w_1)dw_1 f_W(w_2)dw_2
	\\&\le 2\left(\frac{\alpha \omega_d}{\alpha-1}\right)^2 \int_w^\infty \int_w^{w_2} g(w,w_1) g(w,w_2) f_W(w_1)dw_1 f_W(w_2)dw_2\\
	&= \bigO{1} \left(\int_w^\infty g(w,w_1) w_1^{-\beta} dw_1\right)^2 
		= \bigO{w^{4 + 2a-2\beta}}.
\end{align*}

We proceed with the first integral in~\eqref{eq:int_large_weights_1}. Note that $w_1 \le w_2$ which implies that $g(w,w_1) \le g(w,w_2)$. In addition, recall the function $h(w,w_1,w_2)$ from~\eqref{eq:notn_h_wst} and observe that because $w_1 \le w \le w_2$,
\[
	h(w,w_1,w_2) = [g(w,w_1)\wedge g(w,w_2)] g(w_1,w_2) = g(w,w_1)g(w_1,w_2) = w w_1^{2a} w_2.
\]
Next, by~\eqref{eq:min_UB_spatial},
\begin{align*}
    &2 \int_w^\infty \int_1^w \mathcal{S}_w(w_1,w_2) f_W(w_1)dw_1 f_W(w_2)dw_2
    \\&\le 2\left(\frac{\alpha \omega_d}{\alpha-1}\right)^2 w \int_w^\infty \int_1^w w_1^{2a} w_2 f_W(w_1)dw_1 f_W(w_2)dw_2 \\
    &= \bigO{1} w \int_w^\infty w_2^{1 - \beta} dw_2 \int_1^w w_1^{2a-\beta} dw_1 \\
    &= \begin{cases}
    		\bigO{w^{3-\beta}} &\text{if } \beta > 2a+1,\\
    		\bigO{w^{2-2a} \log w} &\text{if } \beta = 2a+1,\\
    		\bigO{w^{4+2a-2\beta}} &\text{if } \beta < 2a+1.
    	\end{cases}
\end{align*}
Since each term on the RHS is $o(w \log(w))$ when $\beta > (a+\frac{3}{2}) \vee 2$ and $a > 0$ or $\beta = a + \frac{3}{2} > 2$ with $a > \frac{1}{2}$ the result follows.
\end{proof}

\subsection{The case $\beta>(a+\frac{3}{2})\vee 2$ and $a \in [0,\infty)$}
By Lemma~\ref{lem:T_M_k_main_term} it suffices to show that
\begin{equation}\label{eq:conclusion_const_bgt1}
	\lim_{w \to \infty} \frac{\overline{T}(w)}{w} 
	= \int_1^{\infty}\int_1^\infty I(w_1,w_2) f_W(w_1)dw_1 f_W(w_2)dw_2,
\end{equation}
where 
\begin{align*}
    \overline{T}(w)
    &:=\int_{1}^{w} \int_{1}^{w} \mathcal{S}_{w}(w_1,w_2) f_W(w_1)dw_1 f_W(w_2)dw_2= 2 \int_{1}^{w} \int_{1}^{w_2} \mathcal{S}_{w}(w_1,w_2) f_W(w_1)dw_1 f_W(w_2)dw_2
\end{align*}
using the symmetry in $w_1$ and $w_2$. 

The following observation will prove useful throughout the whole section: for any $t,s > 0$ it holds that
\begin{equation}\label{eq:kappa_factor_extraction}
	\kappa(\|t^{1/d} s^{a/d}\mathbf{x}\|,tw_1,sw_2) = \kappa(\|\mathbf{x}\|,w_1,w_2).
\end{equation}

\begin{lemma}\label{lem:conv_T_case_1}
Let $\beta > (a+\frac{3}{2}) \vee 2$, with $a \ge 0$. Then, as $w \to \infty$,
\[
	\frac{\overline{T}(w)}{w} \to \int_1^\infty \int_1^\infty I(w_1,w_2) f_W({w}_1)d{w}_1 f_W({w}_2)d{w}_2.
\]
\end{lemma}

\begin{proof}
Recall that
\[
	\mathcal{S}_w(w_1,w_2) = \iint_{{(\mathbb{R}^d)^2}} \kappa(\|\mathbf{x}\|,w,w_1) \kappa(\|\mathbf{y}\|,w,w_2)
	\kappa(\|\mathbf{x}-\mathbf{y}\|,w_1,w_2) d\mathbf{x} d\mathbf{y}.
\]

We now apply the following change of variables:
\begin{equation}\label{eq:cov_bgt1}
    \mathbf{x} =w^{1/d} \mathbf{x}^\prime \quad \text{and} \quad \mathbf{y}^\prime = \mathbf{y} - \mathbf{x},
\end{equation}
to obtain
\begin{align*}
	\mathcal{S}_w(w_1,w_2)
	&= w \iint_{{(\mathbb{R}^d)^2}} \kappa(\|w^{1/d} \mathbf{x}^\prime\|,w,w_1) 
		\kappa(\|\mathbf{y}^\prime + w^{1/d} \mathbf{x}^\prime\|,w,w_2) \kappa(\|\mathbf{y}^\prime\|,w_1,w_2) {d\mathbf{x}^{\prime} d\mathbf{y}^{\prime}}\\
	&= w \iint_{{(\mathbb{R}^d)^2}} \kappa(\|\mathbf{x}^\prime\|,1,w_1)
		\kappa(\|\mathbf{x}^\prime + w^{-1/d}\mathbf{y}^\prime\|,1,w_2) \kappa(\|\mathbf{y}^\prime\|,w_1,w_2) d\mathbf{x}^\prime d\mathbf{y}^\prime,
\end{align*}
where we also use the scaling relation~\eqref{eq:kappa_factor_extraction} for the first two $\kappa$ terms in the integrand.

We now apply dominated convergence twice. For this we first note that the integrand above is dominated by $\kappa(\|\mathbf{x}^\prime\|,1,w_1) \kappa(\|\mathbf{y}^\prime\|,w_1,w_2)$. Hence, using straightforward calculations based on the definition of $\kappa$, we get
\begin{align*}
	\frac{\mathcal{S}_w(w_1,w_2)}{w} 
	&\le \iint_{{(\mathbb{R}^d)^2}} \kappa(\|\mathbf{x}\|,1,w_1) \kappa(\|\mathbf{y}\|,w_1,w_2) d\mathbf{x} d\mathbf{y}
		= \xi_\alpha^2 {w_1}(w_1 \vee w_2)(w_1 \wedge w_2)^a,
\end{align*}
which is finite for any finite $w_1,w_2$. Therefore since $\kappa(\|\mathbf{x} + w^{-1/d}\mathbf{y}\|,1,w_2) \to \kappa(\|\mathbf{x}\|,1,w_2)$ pointwise, it follows that
\[
	\frac{\mathcal{S}_w(w_1,w_2)}{w} \to \iint_{{(\mathbb{R}^d)^2}} \kappa(\|\mathbf{x}\|,1,w_1)
			\kappa(\|\mathbf{x}\|,1,w_2) \kappa(\|\mathbf{y}\|,w_1,w_2) d\mathbf{x} d\mathbf{y} = I(w_1,w_2),
\]
by dominated convergence.

For the last step we first observe that
\begin{align*}
	\frac{\overline{T}(w)}{w} 
	&= \int_1^\infty \int_1^\infty \frac{\mathcal{S}_w(w_1,w_2)}{w} \ind{w_1,w_2\le w} f_W(w_1)dw_1 f_W(w_2)dw_2\\
	&\le \xi_\alpha^2 \int_1^\infty \int_1^\infty 
		w_1^a(w_1 \vee w_2)(w_1 \wedge w_2)^a
		f_W(w_1)dw_1 f_W(w_2)dw_2\\
	&= 2 \xi_\alpha^2 (\beta-1)^2 \int_1^\infty \int_{w_1}^\infty w_1^{2a}w_2 (w_1w_2)^{-\beta}
		{dw_2 dw_1}.
\end{align*}
It is straightforward to see that the above integral is finite. Hence, since $\frac{\mathcal{S}_w(w_1,w_2)}{w} \ind{w_1,w_2\le w}$ converges pointwise to $I(w_1,w_2)$, which is integrable, we apply dominated convergence to conclude
\[
	\frac{\overline{T}(w)}{w} \to
	\int_1^\infty \int_1^\infty I(w_1,w_2) f_W(w_1)dw_1 f_W(w_2)dw_2,
\]
as required. This proves Theorem \ref{thm:scaling} for $\beta>(a+\frac{3}{2})\vee 2$ and $a \in (0,\infty)$.
\end{proof}


\subsection{The case $\beta=a+\frac{3}{2}$ and $a \in (\frac{1}{2},\infty)$}

The proof for this case follows a similar approach to that for $\beta > (a+\frac{3}{2})\vee 2$. Again, by Lemma~\ref{lem:T_M_k_main_term} we only need to consider $\overline{T}(w)$ instead of $T(w)$:

\begin{lemma}\label{lem:conv_T_case_2}
Let $\beta = a+\frac{3}{2}$, with $a > \frac{1}{2}$. Then, as $w \to \infty$,
\[
	\frac{\overline{T}(w)}{w \log(w)} \to {(\beta-1)^2\int_0^{\infty} I(1,r) r^{-\beta}dr.}
\]
\end{lemma}

\begin{proof}
Note that a similar change of variables as applied in Lemma~\ref{lem:conv_T_case_1} extracts a scaling factor $w$ outside the integral. However, for this case we need to extract an additional factor $\log(w)$. This is achieved by also scaling the variables $w_1$ and $w_2$. The proof thus proceeds as follows. We first apply a change of variable on $w_1$ and $w_2$ and then on $\mathbf{x}$ and $\mathbf{y}$. After this we get that $\overline{T}(w)$ is equal to $w \log(w)$ times some integral. We then use dominated convergence, in a similar way as in the proof of Lemma~\ref{lem:conv_T_case_1}, to show that this integral converges. 

Apply the change of variables $(w_1,w_2)\mapsto (u,r)$ given by
\[
	w_1 = w^u \quad \text{and} \quad w_2 = w^u r.
\]
Then $dw_1 = \log(w) w^udu$ and $dw_2 = w^udr$ so that{
\begin{align*}
	&\int_1^w \int_1^w \mathcal{S}_w(w_1,w_2) f_W(w_1)dw_1 f_W(w_2)dw_2
	\\&= \int_0^1 \int_{w^{-u}}^{w^{1-u}} \mathcal{S}_w(w^u,w^ur)w^{2u} \log(w) f_W(w^u) f_W(w^u r) dr du\\
	&= \int_0^1 \int_1^{w^{1-u}} \mathcal{S}_w(w^u,w^ur)w^{2u-2u\beta} \log(w) f_W(1) f_W(r) dr du\\&\hspace{10 pt}+(\beta-1)\int_0^1 \int_{w^{-u}}^{1} \mathcal{S}_w(w^u,w^ur)w^{2u-2u\beta} \log(w) f_W(1) r^{-\beta} dr du. \numberthis \label{eq:b_e_t1_rev_sum_of_tw_exp}
\end{align*}}

Next consider
\[
	\mathcal{S}_w(w^u,w^ur) = \iint_{{(\mathbb{R}^d)^2}} \kappa(\|\mathbf{x}\|,w,w^u) \kappa(\|\mathbf{y}\|, w, w^u r)
	\kappa(\|\mathbf{x}-\mathbf{y}\|,w^u,w^u r) d\mathbf{x} d\mathbf{y}.
\]
Here, we apply the change of variables
\[
	\mathbf{x} = w^{\frac{au+1}{d}} \mathbf{x}^\prime \quad \text{and} \quad \mathbf{y} = w^{\frac{(a+1)u}{d}} \mathbf{y}^\prime + w^{\frac{au+1}{d}} \mathbf{x}^\prime,
\]
to get
\begin{align*}
	\mathcal{S}_w(w^u,w^ur)&= w^{2au + u + 1} \iint_{{(\mathbb{R}^d)^2}} \kappa\left(\|w^{\frac{au+1}{d}} \mathbf{x}^\prime\|,w,w^u\right)
		\kappa\left(\|w^{\frac{au+1}{d}} (\mathbf{x}^\prime + w^{\frac{1-u}{d}}\mathbf{y}^\prime)\|,w^u,w^u r\right)
		\\&\hspace{200 pt}\times\kappa\left(\|w^{\frac{au+u}{d}} \mathbf{y}^\prime\|,w,w^u r\right) d\mathbf{x}^\prime d\mathbf{y}^\prime \\
	&= w^{2au + u + 1} \iint_{{(\mathbb{R}^d)^2}} \kappa\left(\|\mathbf{x}^\prime\|,1,1\right)
		\kappa\left(\|\mathbf{x}^\prime + w^{-\frac{1-u}{d}}\mathbf{y}^\prime\|,1,r\right)
		\kappa\left(\| \mathbf{y}^\prime\|,1, r\right) d\mathbf{x}^\prime d\mathbf{y}^\prime,
\end{align*}
where we have used~\eqref{eq:kappa_factor_extraction} twice to go from the second to the third line. Observe that the integrand is bounded by $\kappa\left(\|\mathbf{x}^\prime\|,1,1\right) \kappa\left(\|\mathbf{y}^\prime\|,1,r\right)$, which is integrable. Hence, by dominated convergence,
\[
	\frac{\mathcal{S}_w(w^u,w^ur)}{w^{2au+u+1}} 
	\to \iint_{\mathbb{R}^2} \kappa\left(\|\mathbf{x}^\prime\|,1,1\right)
			\kappa\left(\|\mathbf{x}^\prime\|,1,r\right)
			\kappa\left(\| \mathbf{y}^\prime\|,1, r\right) d\mathbf{x}^\prime d\mathbf{y}^\prime = I(1,r).
\]
{Recall the the split of $\overline{T}(w)$ in \eqref{eq:b_e_t1_rev_sum_of_tw_exp}. Consider the first term on the RHS of \eqref{eq:b_e_t1_rev_sum_of_tw_exp}. We have}
\begin{align*}
	&\int_0^1 \int_1^{w^{1-u}} \mathcal{S}_w(w^u,w^ur)w^{2u-2u\beta} \log(w) f_W(1) f_W(r) dr du\\
	&= w \log(w) f_W(1) \int_0^1 \int_1^{w^{1-u}} \frac{\mathcal{S}_w(w^u,w^ur)}{w^{2au +u +1}} w^{2au + 3u-2u\beta} 
		f_W(r) dr du \\
	&= w \log(w) f_W(1) \int_0^1 \int_1^{w^{1-u}} \frac{\mathcal{S}_w(w^u,w^ur)}{w^{2au +u +1}} f_W(r) dr du,
\end{align*}
where we used that $2a + 3 - 2\beta = 0$, since $\beta = a + \frac{3}{2}$.

We now apply dominated convergence. For this we recall that $\frac{\mathcal{S}_w(w^u,w^ur)}{w^{2au+u+1}} \to I(1,r)$ and that
\[
	I(1,r) \le \left(\int_{\mathbb{R}^d} \kappa(\|\mathbf{x}\|,1,r) d\mathbf{x}\right)^2
\]
where the RHS is integrable over $(1,\infty)$ with respect to $f_W(r)$. Therefore, by dominated convergence, {the first term on the RHS of \eqref{eq:b_e_t1_rev_sum_of_tw_exp} when divided by $w \log{(w)}$, converges as $w\to \infty$ to}
\[ 
	f_W(1) \int_0^1 \int_1^{\infty} I(1,r) f_W(r) dr du = (\beta-1) \int_1^{\infty} I(1,r) f_W(r)dr.
\]
{Now consider the second term on the RHS of \eqref{eq:b_e_t1_rev_sum_of_tw_exp}. Arguing similarly, we note that when divided by $w \log{(w)}$, this term becomes
\begin{align*}
    f_W(1)(\beta-1)\int_0^1 \int_{w^{-u}}^1 \frac{\mathcal{S}_w(w^u,w^ur)}{w^{2au+u+1}}r^{-\beta}drdu.
\end{align*}
We want to apply dominated convergence again. Here we use a different representation of the integral $\mathcal{S}_w(w^u,w^ur)$. We apply the change of variables
\[
	\mathbf{y} = w^{\frac{au+1}{d}} \mathbf{y}^\prime \quad \text{and} \quad \mathbf{x} = w^{\frac{(a+1)u}{d}} \mathbf{x}^\prime + w^{\frac{au+1}{d}} \mathbf{y}^\prime,
\]
to obtain in a similar fashion as before,
\begin{align*}
    &\frac{\mathcal{S}_w(w^u,w^ur)}{w^{2au+u+1}}\\&=\iint_{(\mathbb{R}^d)^2}\kappa(\|w^{\frac{(a+1)u}{d}} \mathbf{x}^\prime + w^{\frac{au+1}{d}} \mathbf{y}^\prime\|,w,w^u)\kappa(\|\mathbf{y}^\prime\|,1,r)\kappa(\|\mathbf{x}^\prime\|,1,r)d\mathbf{x}^\prime d\mathbf{y}^\prime. \numberthis \label{eq:rep_2_cov_b_e_t1}   
\end{align*}
In particular, we note that $\frac{\mathcal{S}_w(w^u,w^ur)}{w^{2au+u+1}}\ind{w^{-u}<r<1} \to I(1,r)\ind{0<r<1}$, as $w \to \infty$. Further using the representation \eqref{eq:rep_2_cov_b_e_t1}, we note that $\frac{\mathcal{S}_w(w^u,w^ur)}{w^{2au+u+1}}\ind{w^{-u}<r<1}$ is bounded from above by $\left(\int_{\mathbb{R}^d} \kappa(\|\mathbf{x}\|,1,r) d\mathbf{x}\right)^2\ind{0<r<1}=\overline{C}r^{2a}\ind{0<r<1}$, when $0<r<1$, using the definition of $\kappa$, for some constant $\overline{C}>0$. Using that $\beta=a+3/2$ and $a>1/2$, it is easily checked that the integral $(\beta-1)\int_0^1 \int_{0}^1 r^{2a}r^{-\beta}drdu$ is finite. Hence, by dominated convergence, the second term on the RHS of \eqref{eq:b_e_t1_rev_sum_of_tw_exp} when divided by $w\log{(w)}$, converges to
\begin{align*}
    f_W(1)(\beta-1) \int_0^1 \int_0^{1} I(1,r) r^{-\beta} dr du = (\beta-1)^2 \int_0^{1} I(1,r) r^{-\beta}dr.
\end{align*}
Overall, in this case, the quantity $\frac{\overline{T}(w)}{w\log{(w)}}$ converges to $(\beta-1)^2 \int_0^{\infty} I(1,r) r^{-\beta}dr$, as $w \to \infty$.
}This completes the proof of Theorem \ref{thm:scaling} for $\beta=a+\frac{3}{2}$ and $a\in(\frac{1}{2},\infty)$.
\end{proof}


\subsection{The case $\beta \in(2,a+\frac{3}{2})$ and $a \in (\frac{1}{2}, \infty)$}

%

\begin{lemma}\label{lem:conv_T_case_3}
Let $\beta \in (2,a+\frac{3}{2})$, with $a \in (\frac{1}{2},\infty)$. Then as $w \to \infty$,
\begin{align*}
    &\frac{T(w)}{w^{4+2a-2\beta}} \to  \int_{0}^{\infty}\int_{0}^{\infty} 
    \mathcal{S}_1(w_1,w_2) f_W(w_1)dw_1 f_W(w_2)dw_2.
\end{align*}
\end{lemma}
\begin{proof}
Recall the expression for $T(w)$ from \eqref{eq:def_T_w}. Furthermore, note that for any $t > 0$, $\mathbf{x} \in \mathbb{R}^d$ and $w_1,w_2 > 0$,
\[
	\kappa(\|t^{(1+a)/d} \mathbf{x}\|, tw_1, tw_2) = \kappa(\|\mathbf{x}\|,w_1,w_2).
\] 
We now apply the change of variables 
\begin{equation}\label{eq:cov_blt1}
    \mathbf{x}'=\frac{\mathbf{x}}{w^{(1+a)/d}}, \quad \mathbf{y}'=\frac{\mathbf{y}}{w^{(1+a)/d}}, \quad
    w_1'=\frac{w_1}{w}, \quad \text{and} \quad
    w_2'=\frac{w_2}{w}.
\end{equation}
Straightforward calculations imply that
\begin{align*}
    T(w)=w^{4+2a-2\beta}\int_{1/w}^{\infty}\int_{1/w}^{\infty} \mathcal{S}_1(w_1',w_2') f_W(w_1^\prime)w_1^\prime f_W(w_2^\prime)dw_2^\prime.
\end{align*}

To prove that the integral converges we first note that, by~\eqref{eq:scaling_T_in_m_inv}, 
\[
	\int_{1/w}^{\infty}\int_{1/w}^{\infty} \mathcal{S}_1(w_1',w_2') f_W(w_1^\prime)dw_1'f_W(w_2^\prime)w_2'
	= \frac{T(w)}{w^{4+2a-2\beta}} = \bigT{1}.
\]
Hence, by monotone convergence,
\[
	\int_{1/w}^{\infty}\int_{1/w}^{\infty} \mathcal{S}_1(w_1',w_2') f_W(w_1^\prime)dw_1' f_W(w_2^\prime)dw_2' \to 
	\int_{0}^{\infty}\int_{0}^{\infty} \mathcal{S}_1(w_1',w_2') f_W(w_1^\prime)dw_1' f_W(w_2^\prime)dw_2' < \infty,
\]
and thus
\begin{equation}\label{eq:conv_T_M_final_case}
	\lim_{w \to \infty} \frac{T(w)}{w^{4+2a-2\beta}}
	= \int_{0}^{\infty}\int_{0}^{\infty} \mathcal{S}_1(w_1,w_2) f_W(w_1)dw_1 f_W(w_2)dw_2
\end{equation}
This completes the proof of Theorem \ref{thm:scaling} for $\beta \in(2,a+\frac{3}{2})$ and $a\in(\frac{1}{2},\infty)$.
\end{proof}

\section{Discussion}\label{sec:disc}
This section consists of discussions of our proof, possible extensions, and simulations.
\subsection{Geometry of typical triangles}
\begin{figure}[H]
    \centering
    \resizebox{0.3\textwidth}{!}{
        \begin{tikzpicture}

\def\triangleheight{5};
\def\mathleftupper{(\|\mathbf{x}\|,{w}_1)};
\def\mathrightupper{(\|\mathbf{y}\|,{w}_2)};
\def\mathlower{(X_0=\mathbf{0},W_0)};
\def\mathbrace{\|\mathbf{x}-\mathbf{y}\|};

\def\trianglewidth{4.2};

\coordinate (lower) at (\trianglewidth/2,0);
\coordinate (upperleft) at (0,\triangleheight);
\coordinate (upperright) at (\trianglewidth,\triangleheight);

\draw [decorate,decoration = {brace,raise=3pt},ultra thick] (upperleft) -- node[above=6pt,align=center] {$\mathbrace$} (upperright);
\draw[thick] (upperright) node[right,align=center] {$\mathrightupper$} -- (lower) node[below,align=center] {$\mathlower$} -- (upperleft) node[left,align=center] {$\mathleftupper$} -- cycle;
\end{tikzpicture}
    }
    \caption{The triangle between $0$ and its two randomly sampled neighbors $(\mathbf{x},w_1)$ and $(\mathbf{y},w_2)$.}
    \label{fig:typ_tri}
\end{figure}
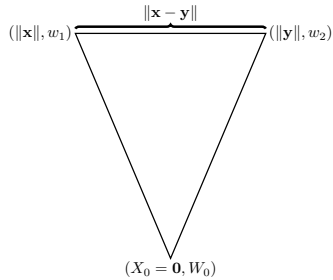

Note that because of the Poisson point process representation of the neighbors of $0$ in $\mathbb{G}^{\infty}$ as obtained in Lemma \ref{lem:nbr_PP}, the fraction $T(w)/\mathcal{M}(w)^2$ can be seen as the probability that two randomly sampled neighbors of $0$ have an edge between them, given that $W_0=w$. Consequently, since the density $f_k(w)$ behaves as a Dirac mass at $\mathcal{M}^{-1}(k)$ thanks to Proposition \ref{prop:conc_weights}, the fraction $T(\mathcal{M}^{-1}(k))/k^2$ can be interpreted as the probability that two randomly sampled neighbors of $0$ are also neighbors, given $d_0(\mathbb{G}^{\infty})=k$. The change of variables that we encounter in Section \ref{sec:main_thm_proof} gives us interesting insights about the locations and weights $(\mathbf{x},w_1),(\mathbf{y},w_2)\in \mathbb{R}^d \times \mathbb{R}_+$ of two such randomly sampled neighbors. This further gives us insight into, what a randomly sampled triangle between $0$ and two such neighbors look like. Figure~\ref{fig:typ_tri} shows an abstract depiction of such a typical triangle. The quantity $\|\mathbf{x}-\mathbf{y}\|$ is the spatial distance from $\mathbf{x}$ to $\mathbf{y}$ in $\mathbb{R}^d$. We will next explain how this distance and the weights of each node in a typical triangle behaves in the three different regimes of $\beta$ and $a$.

To begin, recall Corollary \ref{cor:limits_M_inv_k} and \eqref{eq:lim_t_over_sigma}. Our main idea is to analyse the change of variables encountered in Lemmas \ref{lem:conv_T_case_1}, \ref{lem:conv_T_case_2} and~\ref{lem:conv_T_case_3}, to get the convergence of $T(w)/\sigma(w)$, for $w=\mathcal{M}^{-1}(k)$. 

\subsubsection{Case $\beta>(a+\frac{3}{2})\vee 2$ and $a\in [0,\infty)$} 
Let us carefully dissect the variable change in \eqref{eq:cov_bgt1}. Clearly, when $w=\mathcal{M}^{-1}(k)$, by \eqref{eq:cov_bgt1}, $\|\mathbf{x}\|$ scales like $\mathcal{M}^{-1}(k)^{1/d}$, which scales like $k^{1/d}$ (recall Proposition \ref{prop:scaling_Minv_k}). On the other hand, we did not need to scale $\mathbf{y}-\mathbf{x}=\mathbf{y}'$, i.e. $\|\mathbf{x}-\mathbf{y}\|$ remains bounded, which implies that $\mathbf{y}$ stays in a ball of finite radius about $\mathbf{x}$. Thus, $\|\mathbf{y}\|$ is also of the order $k^{1/d}$, and is of constant order spatial distance away from $\mathbf{x}$. Finally, note that we did not have to rescale the respective weights $w_1$ and $w_2$ corresponding to $x$ and $y$, which means that they are of constant order.

\subsubsection{Case $\beta=a+\frac{3}{2}$ and $a \in (\frac{1}{2},\infty)$.} Arguing as in the previous case, and observing that $\mathcal{M}^{-1}(k)$ scales like $k$ using Proposition \ref{prop:scaling_Minv_k}, we note from the proof of Lemma \ref{lem:conv_T_case_2} that when $w=\mathcal{M}^{-1}(k)$, the relevant change of variables yield $\mathbf{x}=C_1\mathbf{x}'k^{1/d}=C_2\overline{\mathbf{x}}k^{(1+au')/d}$, for some constants $C_1, C_2>0$ and $u'\in (0,1)$. Consequently, $\|\mathbf{x}\|=\Theta(k^{(1+au')/d})$. Next note that $\mathbf{y}-\mathbf{x}=\mathbf{y}'=C_3\overline{\mathbf{y}}k^{(au'+u')/d}$, where $C_3>0$ is some constant. This implies that $\|\mathbf{x}-\mathbf{y}\|=\Theta(k^{(au'+u')/d})$. Note that the distance $\|\mathbf{x}-\mathbf{y}\|$ between $\mathbf{x}$ and $\mathbf{y}$ is of smaller order than $\|\mathbf{x}\|$, since $u' \in (0,1)$. This means that $\|\mathbf{y}\|$ is also of the order $\|\mathbf{x}\|$. Finally, note that for this case we had to rescale the weights as well. In particular, we changed variables as $w_1=e^u=C_4k^{u'}$, for some constant $C_4>0$, which means that the weight $w_1$ corresponding to $\mathbf{x}$ is $\Theta(k^{u'})$, with $u'\in (0,1)$. We also had to do a variable change as $\mathbf{w}={w_2}/{C_4k^{u'}}={w_2}/{w_1}$, i.e. the weight $w_2$ corresponding to $\mathbf{y}$ is of the same order as $w_1$, the weight of $\mathbf{x}$.

\subsubsection{Case $\beta \in (2,a+\frac{3}{2})$ and $a \in (\frac{1}{2},\infty)$. }
Note that for this case, for $w=\mathcal{M}^{-1}(k)$, the variable change in \eqref{eq:cov_blt1} suggests that both $\mathbf{x}$ and $\mathbf{y}$ are $\Theta((\mathcal{M}^{-1}(k))^{(1+a)/d})$, while the respective weights $w_1$ and $w_2$ are both $\Theta(\mathcal{M}^{-1}(k))$. Finally, note that for this case we did not have to rescale $\mathbf{x}-\mathbf{y}$, which means that the distance $\|\mathbf{x}-\mathbf{y}\|$ between $\mathbf{x}$ and $\mathbf{y}$ is also of the order $(\mathcal{M}^{-1}(k))^{(1+a)/d}$. Depending on how $\mathcal{M}^{-1}(k)$ scales in $k$ (recall Proposition \ref{prop:scaling_Minv_k}), the orders of the different quantities $\|\mathbf{x}\|,\|\mathbf{y}\|,w_1,w_2,\|\mathbf{x}-\mathbf{y}\|$ differ in different regimes of $\beta$ and $a$. The relationship and the behaviour is summarized in Table \ref{tab:abc}.

\begin{table}[h]
\centering
\begin{tabular}{l | l | l |l}
Relation between $\beta$ and $a$ & $W_0$ & $(\|\mathbf{x}\|,w_1)$ and $(\|\mathbf{y}\|,w_2)$ &$\|\mathbf{x}-\mathbf{y}\|$ \\
\hline \hline
$\beta>(a+\frac{3}{2})\vee2,a\in [0,\infty)$ & $\Theta(k)$ & $(\Theta(k^{1/d}),\Theta(1))$ & $\Theta(1)$ \\ \hline
$\beta=a+\frac{3}{2}, a\in (\frac{1}{2},\infty)$ & $\Theta(k)$ & $\left(\Theta\left(k^{\frac{1+au}{d}} \right),\Theta\left(k^u\right)\right)$& $\Theta\left(k^{\frac{au+u}{d}}\right)$\\
 & & for $u \in (0,1)$& for $u \in (0,1)$\\ \hline
$\beta \in(a+1,a+\frac{3}{2}), a \in (\frac{1}{2},\infty)$ & $\Theta(k)$ & $\left(\Theta\left(k^{\frac{1+a}{d}} \right),\Theta\left(k\right)\right)$& $\Theta\left(k^{\frac{1+a}{d}} \right)$\\ \hline
$\beta=a+1, a \in (1,\infty)$ & $\Theta\left(\frac{k}{\log{k}}\right)$ & $\left(\Theta\left(\left(\frac{k}{\log{k}}\right)^{\frac{1+a}{d}} \right),\Theta\left(\frac{k}{\log{k}}\right)\right)$& $\Theta\left(\left(\frac{k}{\log{k}}\right)^{\frac{1+a}{d}} \right)$\\ \hline
$\beta \in (2,a+1), a\in (1,\infty)$&$\Theta\left(k^{\frac{1}{2+a-\beta}}\right)$ & $\left(\Theta\left(k^{\frac{1+a}{d(2+a-\beta)}} \right),\Theta\left(k^{\frac{1}{2+a-\beta}}\right)\right)$& $\Theta\left(k^{\frac{1+a}{d(2+a-\beta)}}\right)$ \\ \hline
\end{tabular}
\caption{Summary of the behavior of the typical triangle between $0$ and its two randomly sampled neighbors $(\mathbf{x},w_1)$, $(\mathbf{y},w_2)$.}
\label{tab:abc}
\end{table}

\subsection{The Boolean model}
In this section, we discuss a popular related model from the literature. Consider a homogeneous Poisson point process $\Gamma$ on $\mathbb{R}^d$, and a non-negative random variable $W$. Form the open ball $B(\mathbf{x},W_{\mathbf{x}})$ around each point $\mathbf{x} \in \Gamma$, where the collection $(W_{\mathbf{x}})_{\mathbf{x} \in \Gamma}$ is a collection of i.i.d.\ copies of $W$. The random graph formed by placing an edge between distinct pairs $\mathbf{x}, \mathbf{y} \in \Gamma$ whenever $B(\mathbf{x},W_{\mathbf{x}})\cap B(\mathbf{y},W_{\mathbf{y}}) \neq \emptyset$ is called the \emph{Boolean} model \cite{Subcrit_boolean,meester_roy_1996}. Let $W$ have distribution \eqref{eq:Pareto_weights}. Clearly, this random graph model is a SIRG, with connection function \eqref{eq:conn_fn_form} with $\alpha=\infty$ and $g(s,t)=(s+t)^d$. One can also consider the corresponding $\alpha \in(1,\infty)$ version. For this model, $\CExp{d_0(\mathbb{G}^{\infty})}{W_0=w}=\mathcal{M}(w)$ is integrable with respect to the density \eqref{eq:Pareto_weights} when $\beta>d$. 

A related model was studied in \cite{WDRCM_Rec_trans}, where the authors consider a weight function of the form $g^{\mathrm{sum}}(s,t)=\beta^{-1}(s^{-\gamma/d}+t^{-\gamma/d})^{-d}$, with $s$ and $t$ now uniformly distributed on $(0,1)$, and $\beta \in (0,\infty)$, $\gamma \in [0,1)$ parameters, see \cite[(1.5)]{WDRCM_Rec_trans}. A particular case of this model can be considered, when one takes $\beta=1$, the profile function $\rho$ as considered in \cite[(1.2)]{WDRCM_Rec_trans} to be $\kappa_{\alpha}(\|\mathbf{x}-\mathbf{y}\|,s',t')$ as in \eqref{eq:conn_fn_form}, where $s'=s^{-\gamma}$, and $t'=t^{-\gamma}$, with the weight function $g(s',t')=(s'+t')^d$. Note that $s'$ and $t'$ now have density \eqref{eq:Pareto_weights} with $\beta-1=1/\gamma$. Again, as in the last paragraph, $\CExp{d_0(\mathbb{G}^{\infty})}{W_0=w}=\mathcal{M}(w)$ is integrable with respect to the density \eqref{eq:Pareto_weights} when $\beta>d$, i.e. $\gamma<1/(d-1)$.

These models does not directly come as a particular case of the weight function \eqref{defn:interpolated_wt_fn} that we consider. However, following the same proof strategy and calculations as in the proof of Theorem \ref{thm:scaling}, it can be shown that the corresponding clustering function $\gamma(k)$ {has an inverse linear scaling:} $\lim_{k \to \infty}k\gamma(k)=\Gamma \in (0,\infty)$. The calculations are very similar to the $a=0$ and $\beta>d$ case for the weight function \eqref{defn:interpolated_wt_fn} that we consider, and so we omit these calculations. {In particular, we state the following theorem without proof:

\begin{theorem}
Consider the infinite SIRG $\mathbb{G}^{\infty}(\kappa,W)$, where $W$ has distribution \eqref{eq:Pareto_weights} and $\kappa=\kappa_{\alpha}$ is of the form \eqref{eq:conn_fn_form}, with $g(s,t)=(s+t)^d$. Assume $\beta>d$. Then, with $\gamma(k)$ denoting the clustering function of this graph,
\[\lim_{k \to \infty}k\gamma(k)=\Gamma,\]
where (recall \eqref{eq:def_xi} and the definitions below it)
\[\Gamma=\frac{\beta-\alpha-1}{\xi_{\alpha}}\int_{1}^{\infty}\int_{1}^{\infty}I(w_1,w_2)f_W(w_1)dw_1 f_W(w_2)dw_2.\]

\end{theorem} As remarked before, we need the condition $\beta>d$ to ensure that $\Exp{d_0(\mathbb{G}^{\infty})}<\infty$.}

\subsection{On finite models}\label{ssec:finite_models} 
Note that Theorem \ref{thm:scaling} gives results on how $\gamma(k)$ scales. However $\gamma(k)$ is the $n \to \infty$ limit of the sequence of clustering functions $\left(\mathrm{CC}_{\mathbb{G}^{(n)}}(k)\right)_{n \geq 1}$ of the finite models $\mathbb{G}^{(n)}$. What happens if we let $k=k_n$ depend on $n$, with $k_n \to \infty$, for the graph $\mathbb{G}^{(n)}$? That is, let us consider the sequence of random variables $\mathrm{CC}_{\mathbb{G}^{(n)}}(k_n)$, and ask the question: do results similar as in the theorems above hold with high probability as $n \to \infty$? We expect if $k_n \to \infty$ sufficiently slowly, then we have a positive answer to the question. We formulate this as a conjecture in this section. We focus on the finite SIRG $\mathbb{G}^{(n)}$, where the locations are i.i.d.\ on $\left[-{n^{1/d}}/{2},{n^{1/d}}/{2} \right]^d$. As remarked earlier, the finite model locally converges to the infinite model $\mathbb{G}^{\infty}$ (see \cite[Theorem 1.2]{LWC_SIRGs_2020}). For the finite model $\mathbb{G}^{(n)}$, recall the clustering function $\mathrm{CC}_{\mathbb{G}^{(n)}}(k)$ from \eqref{eq:cLust_fn_SIRG_finite}. Recall for fixed $k$, $\mathrm{CC}_{\mathbb{G}^{(n)}}(k) \plim \gamma(k)$. For $k_n \to \infty$, denote the sequence $\phi(k_n) \to \infty$ by
\begin{align*}
    \phi(k_n):=\begin{cases}
        \frac{k_n}{\log{k_n}} &\mbox{if}\;\; \beta=a+1, a\in(1,\infty),\\
        k_n &\mbox{otherwise,}
    \end{cases}
\end{align*}
and denote the sequence $\psi(n)\to \infty$ as 
\begin{align*}
    \psi(n)=\begin{cases}{n^{\frac{1}{\beta-1}}} &\mbox{if}\;\; \beta\geq (a+1)\vee 2, a\in[0,\infty),\\
    {n^{\frac{2+a-\beta}{\beta-1}}} &\mbox{otherwise.}\end{cases}
\end{align*}

We conjecture the following limit laws for finite models:
\begin{conjecture}[Scaling for finite graphs]\label{conj:finite_models}
Consider the finite SIRG $\mathbb{G}^{(n)}$.
\begin{itemize}
    \item[a.] Let $k_n \to \infty$ satisfy $\phi(k_n)=\smallO{\psi(n)}$. Then, as $n \to \infty$, $\frac{\mathrm{CC}_{\mathbb{G}^{(n)}}(k_n)}{\gamma(k_n)}\plim 1$. 
\item[b.] Let $k_n \to \infty$ satisfy $\phi(k_n)=(1+\smallO{1})\psi(n)$, as $n \to \infty$.
Then, as $n \to \infty$, $\frac{\mathrm{CC}_{\mathbb{G}^{(n)}}(k_n)}{\gamma(k_n)}\dlim L_0\text{Poi}(\lambda)$, for some $L_0,\lambda>0$.
\item[c.] Let $k_n \to \infty$ satisfy $\phi(k_n)=\omega\left({\psi(n)}\right)$.
Then as $n\to \infty$, $\frac{\mathrm{CC}_{\mathbb{G}^{(n)}}(k_n)}{\gamma(k_n)}\plim 0$.
\end{itemize}
\end{conjecture}

\paragraph{Heuristic arguments for Conjecture \ref{conj:finite_models}.} Let $W$ be a random variable with density \eqref{eq:Pareto_weights}. Observe that by \eqref{eq:deg_law_given_wt}, because $\mathbb{G}^{(n)}$ converges locally to $\mathbb{G}^{\infty}$, the expected number of vertices with degree at least $k_n$ in $\mathbb{G}^{(n)}$ $\sim n\Prob{\mathrm{Poi}(\mathcal{M}(W))>k_n}\sim n\Prob{\mathcal{M}(W)>k_n}=n\Prob{W>\mathcal{M}^{-1}(k_n)}$ as $W \to \infty$. Then \eqref{eq:Pareto_weights} and Proposition \ref{prop:scaling_Minv_k} ensure that when $k_n$ satisfies $\phi(k_n)=\smallO{\psi(n)}$, this expected value does not decay to $0$. In other words, \emph{as soon as it is possible to choose vertices of degree $k_n$ in the graph $\mathbb{G}^{(n)}$, the ratio of $\mathrm{CC}_{\mathbb{G}^{(n)}}(k_n)$ and $\gamma(k_n)$ should converge in probability to $1$, i.e., they are asymptotically the same, as claimed by Conjecture \ref{conj:finite_models} (a).} {There are two main difficult steps in proving Conjecture \ref{conj:finite_models} (a). Firstly, to rigorously prove that as soon as $\phi(k_n)=\smallO{\psi(n)}$, the expected number of vertices in $\mathbb{G}^{(n)}$ with degree $k_n$, does not decay to $0$. Secondly, to prove $\mathrm{CC}_{\mathbb{G}^{(n)}}(k_n)/\gamma(k_n)$ converges to $1$, as soon as $\phi(k_n)=\smallO{\psi(n)}$.} Parts (b) and (c) of Conjecture \ref{conj:finite_models} are to be understood as quantifying the finite-size effects in our models. That is, when $\phi(k_n)=\omega(\psi(n))$ as in Conjecture \ref{conj:finite_models} (c), there are not many vertices of degree $k_n$ in $\mathbb{G}^{(n)}$, and even when there are, they simply do not have enough edges among their neighbors to contribute to the clustering in comparison to their degrees. Conjecture \ref{conj:finite_models} (b) claims that at the boundary of the finite-size effects, which we call the \emph{finite-size threshold}, the limit is random, and in fact Poissonian. {Parts (b) and (c) require a careful and tedious analysis of boundary effects around the finite-size threshold, which constitutes its main difficulty.} A result of similar flavor for Hyperbolic Random Graphs, which corresponds to taking $d=1$ and $\alpha=\infty$ in our model, was proved in \cite[Theorem 1.5]{Clustering_HRGs_2020}. {This proof relied heavily upon the particular choice of the model and its parameters, and does not appear to be easily generalizable to our setting. }

As we have seen in Theorem \ref{thm:scaling}, there are two phase transitions in the scaling of $\gamma(k)$, so three different phases. In Figure \ref{fig:finite_graph_sims}, we see simulations of the clustering as a function of the degree in the finite graph $\mathbb{G}^{(n)}$, in these three different phases of $\beta$ and $a$. Here for $k \in \{d_v(\mathbb{G}^{(n)}):v\in V(\mathbb{G}^{(n)})\}$, the log-log plot of $k \to \mathrm{CC}_{\mathbb{G}^{(n)}}(k)$ is shown. Since in all these three phases, $\gamma(k)=\Theta(k^{-p})$ as $k \to \infty$, for some $p \geq 0$, we see a straight line in the log-log plot, before the decay becomes much faster due to finite-size effects. The (approximate) finite-size threshold is indicated by the green dashed vertical line. In the plots of Figure \ref{fig:finite_graph_sims}, the number of vertices $n$ of $\mathbb{G}^{(n)}$ is fixed to be $n=22000$, $a=2$, $\alpha=1$ and $d=2$ are fixed. $\beta$ is varied to obtain random graphs in the three different regimes: $\beta>a+\frac{3}{2}$, $\beta \in (a+1,a+\frac{3}{2})$ and $\beta\in (\frac{a+3}{2},a+1)$. Here for the third phase, we do not go all the way to $\beta=2$, since, in the regime $\beta \in (2,\frac{a+3}{2})$, $\Exp{d_0(\mathbb{G}^{\infty})}=\infty$, which we recall from Section \ref{sec:m_inv_scaling}.

\begin{figure}[H]
     \centering
     \begin{subfigure}[b]{1\textwidth}
         \centering
   \includegraphics[width=15.7cm]{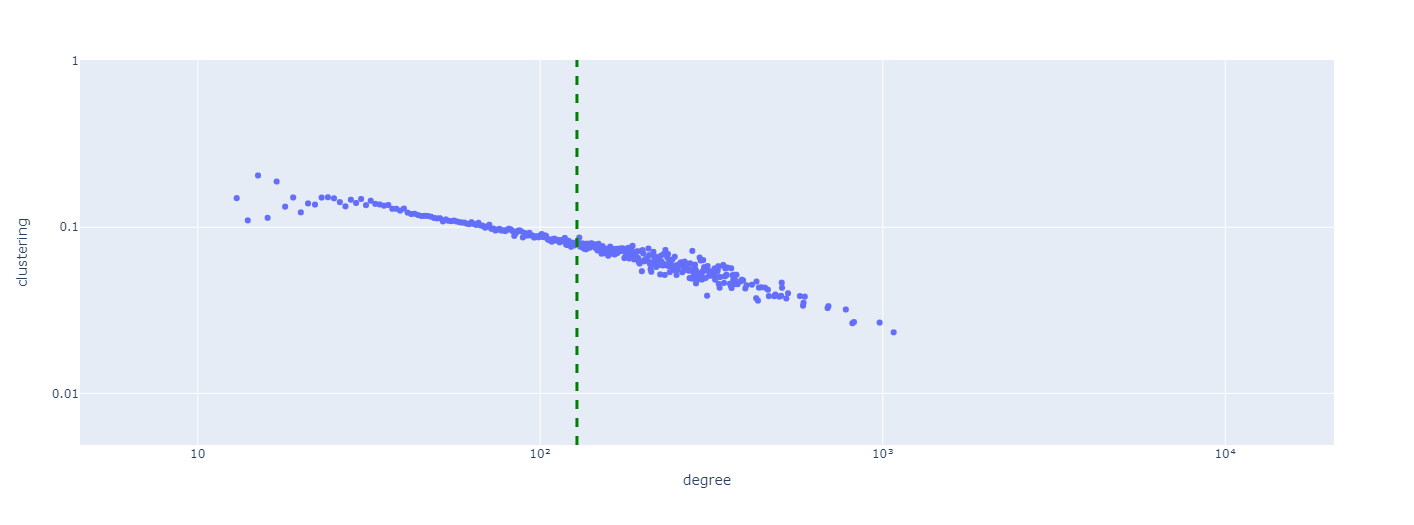}
  \caption{Clustering function in the regime $\beta>a+3/2$. Here $n=22000$, $a=2$ and $\beta=4$. 
  }
    \label{fig:sim_b_g_t1}
     \end{subfigure}
     \hfill
     \begin{subfigure}[b]{1\textwidth}
         \centering
    \includegraphics[width=15.7cm]{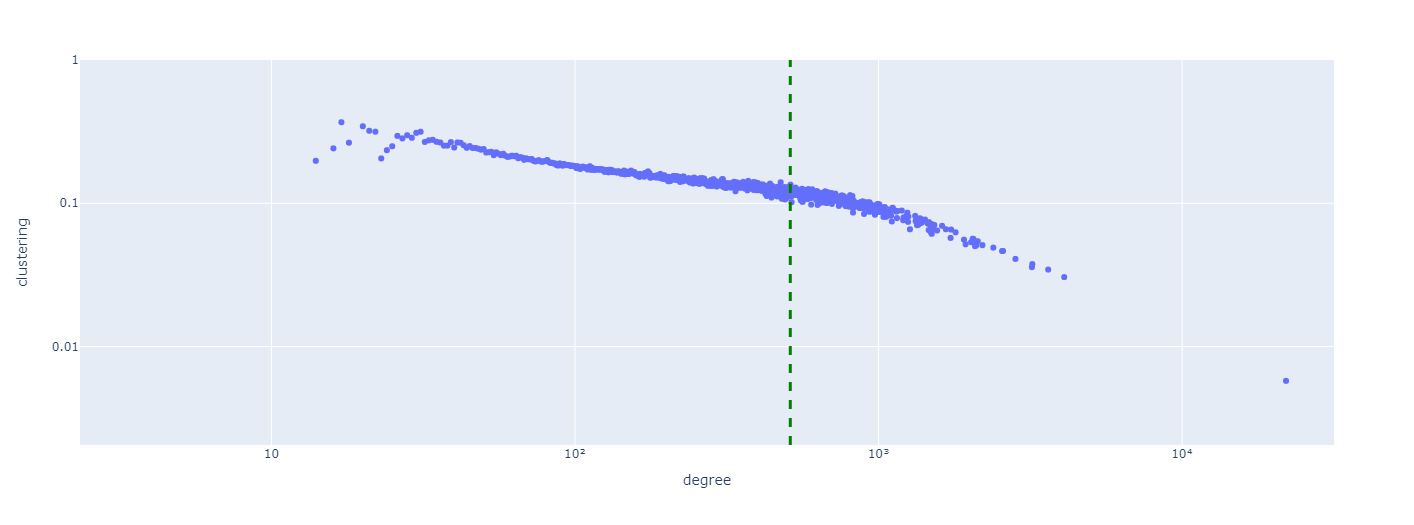}
  \caption{Clustering function in the regime $\beta \in (a+1,a+\frac{3}{2})$. Here $n=22000$, $a=2$ and $\beta=3.01$. }
    \label{fig:sim_b_g_t2}
     \end{subfigure}
     \end{figure}
     \begin{figure}[H]
     \ContinuedFloat
     \begin{subfigure}[b]{1\textwidth}
         \centering
    \includegraphics[width=15.7cm]{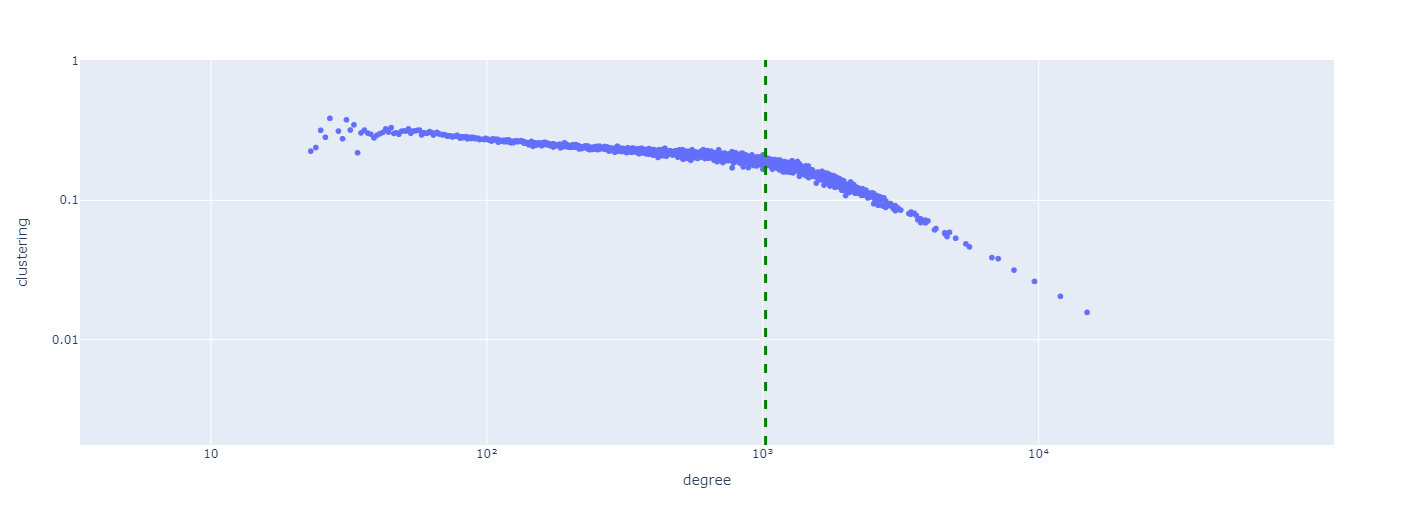}
  \caption{Clustering function in the regime $\beta \in (2,a+1)$. Here $n=22000$, $a=2$ and $\beta=2.6$. }
    \label{fig:sim_b_l_t2}
     \end{subfigure}
        \caption{Plots of clustering against degree in the log-log scale. The vertical green dashed line indicates the approximate finite-size thresholds.}
        \label{fig:finite_graph_sims}
\end{figure}
As can be seen in Figure \ref{fig:finite_graph_sims}, as the weight tail $\beta$ becomes heavier, the decay of the clustering function, as the degree becomes large, becomes slower and slower. The outliers in the simulations for extremely high-degree vertices (extreme right of the plots) show low clustering because of finite-size effects. Note that the clustering contribution from extremely low-degree vertices (extreme left of the plots) are not of interest to us, since we are interested in the clustering behaviour as the degree diverges. 

\begin{figure}[H]
     \centering
     \begin{subfigure}[b]{0.30\textwidth}
         \centering
   \includegraphics[width=4.5cm]{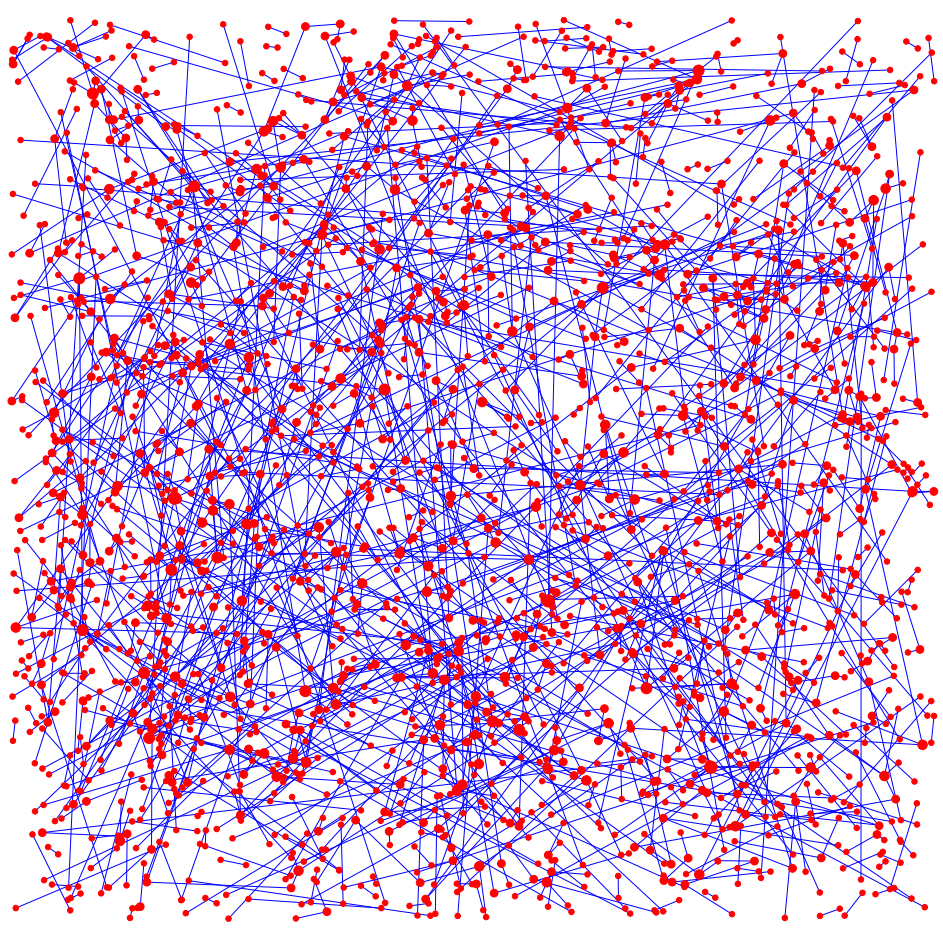}
  \caption{The graph $\mathbb{G}^{(n)}$ in the regime $\beta>a+3/2$. Here $n=8000$, $a=2$, $\beta=4$.}
    \label{fig:sim_graph_b_g_t1}
     \end{subfigure}
     \hfill
     \begin{subfigure}[b]{0.30\textwidth}
         \centering
   \includegraphics[width=4.5cm]{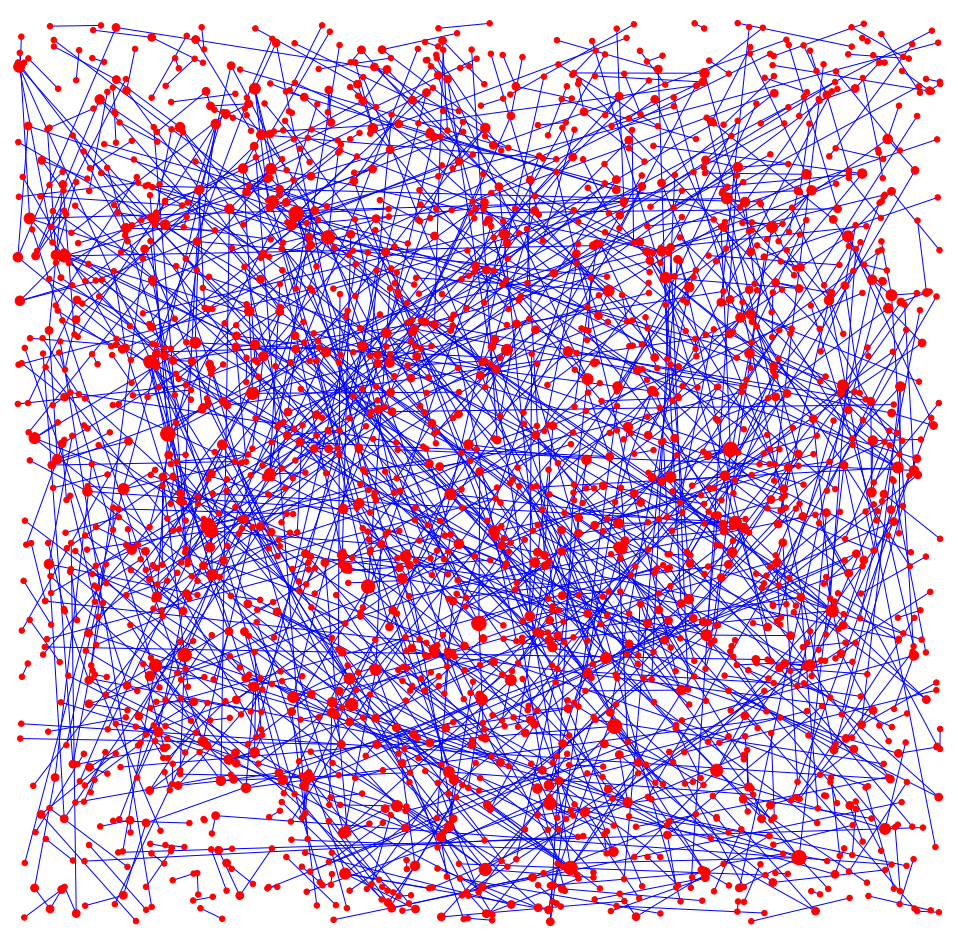}
  \caption{The graph $\mathbb{G}^{(n)}$ in the regime $a+3/2>\beta>a+1$. Here $n=8000$, $a=2$, $\beta=3.01$. }
    \label{fig:sim_graph_b_l_t1}
     \end{subfigure}
     \hfill
     \begin{subfigure}[b]{0.30\textwidth}
         \centering
   \includegraphics[width=4.5cm]{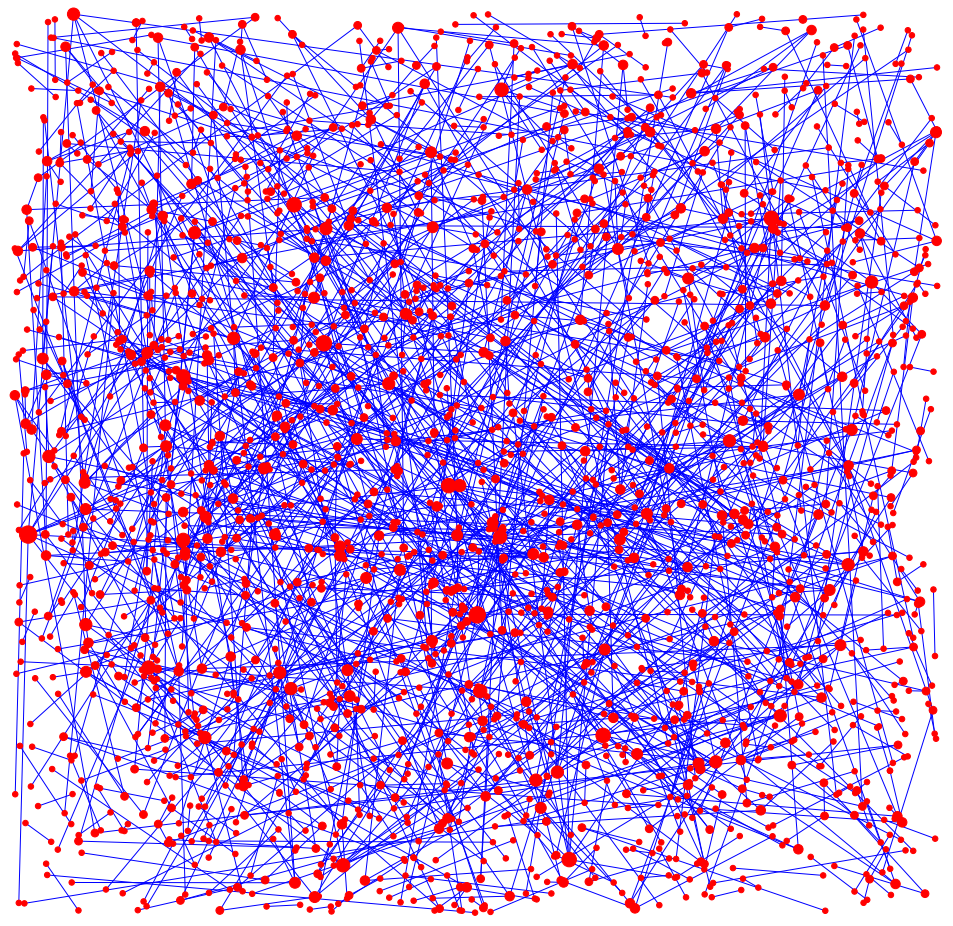}
  \caption{The graph $\mathbb{G}^{(n)}$ in the regime $a+1>\beta>2$. Here $n=8000$, $a=2$, $\beta=2.6$. }
    \label{fig:sim_graph_b_l_t2}
     \end{subfigure}
     
        \caption{In Figures \ref{fig:sim_graph_b_l_t1} and \ref{fig:sim_graph_b_l_t2}, bond percolation is applied so that the average degree is approximately the same as the graph in Figure \ref{fig:sim_graph_b_g_t1}. }
        \label{fig:finite_graphs}
\end{figure}

In Figure \ref{fig:finite_graphs}, for visualization, we have included simulations of the graph $\mathbb{G}^{(n)}$ in the three different regimes $\beta>a+\frac{3}{2}$, $\beta \in (a+1,a+\frac{3}{2})$ and $\beta \in (2,a+1)$. In all these simulations, $n=8000$ and $\alpha=1$ are fixed. Here larger sized vertices have more degree. A uniform bond percolation parameter of $p=1/7000$ is first applied on all the three graphs for better visibility of edges. Further, bond percolation is applied to the graphs in the regimes $\beta \in (a+1,a+\frac{3}{2})$ and $\beta \in (2,a+1)$ so that they have (approximately) the same average degree, as the graph in the regime $\beta>a+\frac{3}{2}$. In Figure \ref{fig:sim_graph_b_l_t2}, we see the graph $\mathbb{G}^{(n)}$ in the regime $\beta \in (2,a+1)$. Here $\beta$ is also chosen such that $\beta>(a+3)/2$, so that $\Exp{d_0(\mathbb{G}^{\infty})}<\infty$. However, the clustering $\gamma(k)$ of the corresponding infinite model converges to a constant, see Theorem \ref{thm:scaling}. As a consequence, the graph admits vertices of very high degrees, and these vertices have neighborhoods that are almost cliques. Could these vertices be part of a large dense subgraph? A future research direction for us is to study other properties of the graphs in this \emph{exotic} regime.

\appendix

\section{Markings and thinnings: proof of Lemma \ref{lem:nbr_PP}}\label{sec:app}

Before proving Lemma \ref{lem:nbr_PP}, we first recall the concepts of marking and thinning of a point process. For more on what follows, we refer the reader to \cite[Chapter 5]{Last_Penrose_LPP}. 

\begin{definition}[Probability kernel] 
For measurable spaces $(\mathbb{X},\mathscr{X})$ and $(\mathbb{Y},\mathscr{Y})$, a probability kernel $K$ from $\mathbb{X}$ to $\mathbb{Y}$ is a map $K:\mathbb{X} \times \mathscr{Y} \to [0,1]$, such that for every fixed $x \in \mathbb{X}$, $K(x,\cdot)$ is a probability measure on $\mathbb{Y}$, and for every fixed $A \in \mathscr{Y}$, $K(\cdot, A)$ is an $\mathscr{X}$-measurable function on $\mathbb{X}$.
\end{definition}

\begin{definition}[Proper point process]
A point process $\eta$ on a measurable space $(\mathbb{X},\mathscr{X})$ is called proper, if there exists $\mathbb{X}$-valued random variables $X_1,X_2,\dots$ and an $\mathbb{N}_0:=\mathbb{N}\cup\{0\}$-valued random variable $\mathscr{N}$, such that almost surely,
\begin{align*}
    \eta=\sum_{i=1}^{\mathscr{N}}\delta_{X_i}.
\end{align*}
\end{definition}
It is a well-known fact that any Poisson point process is proper, see \cite[Corollaries 3.7 and 6.5]{Last_Penrose_LPP}. 

\begin{definition}[Marking]
Let $\eta=\sum_{i=1}^{\mathscr{N}}\delta_{X_i}$ be a proper point process on $\mathbb{X}$. Let $K$ be a probability kernel from $\mathbb{X}$ to $\mathbb{Y}$ for some measurable space $(\mathbb{Y}, \mathscr{Y})$. Let $Y_1,Y_2,\dots$ be $\mathbb{Y}$-valued random variables and {let} $(Y_i)_{i \leq m}$, given $\mathscr{N}=m$ and $(X_i)_{i \leq m}$, be independent random variables with distribution $K(X_i,\cdot)$, $i \leq m$. Then the $\mathbb{X} \times \mathbb{Y}$-valued point process 
\begin{align*}
    \zeta=\sum_{i=1}^{\mathscr{N}}\delta_{X_i,Y_i}
\end{align*} 
is called a \emph{$K$-marking} of $\eta$.
\end{definition}
If there is a probability measure $\mathbb{Q}$ on $\mathbb{Y}$ such that $K(x,\cdot)=\mathbb{Q}(\cdot)$ for all $x \in \mathbb{X}$, then $\zeta$ is called an \emph{independent $\mathbb{Q}$-marking} of $\eta$.
The following is the main theorem on markings of Poisson point processes (see \cite[Theorem 5.6]{Last_Penrose_LPP}) which we state without proof. Recall an $s$-finite measure is a countable sum of finite measures:
\begin{theorem}[Marking theorem]
Let $\zeta$ be a $K$-marking of a proper Poisson point process $\eta$ with $s$-finite intensity measure $\lambda$. Then $\zeta$ is a Poisson point process on $\mathbb{X} \times \mathbb{Y}$ with intensity measure $\lambda \otimes K$, where, for any $C \in \mathscr{X} \times \mathscr{Y}$, 
\begin{align*}
    \lambda \otimes K(C):= \int_{\mathbb{X} \times \mathbb{Y}} \ind{(x,y) \in C} K(x,dy)\lambda(dx).
\end{align*}
\end{theorem}

\begin{definition}[Thinning]
Let $p:\mathbb{X} \to [0,1]$ be measurable, and consider the probability kernel $K$ from $\mathbb{X}$ to $[0,1]$ defined by 
\begin{align*}
    K_p(x,\cdot)=(1-p(x))\delta_0+p(x)\delta_1. \numberthis \label{def:thin_marking}
\end{align*}
Then if $\zeta$ is a $K_p$-marking of a proper point process $\eta$, the point process $\zeta(\cdot \times \{1\})$ on $\mathbb{X}$ is called a $p$-thinning of $\eta$. 
\end{definition}

The following is the main theorem on thinnings of Poisson point process (see \cite[Theorem 5.8]{Last_Penrose_LPP}), which we state without proof:

\begin{theorem}[Thinning theorem]
Let $\zeta$ be a $K_p$-marking of a proper Poisson process $\eta$, where $K_p$ is as in \eqref{def:thin_marking}. Then $\zeta_i:=\zeta(\cdot \times \{i\})$ for $i=1,0$ are independent Poisson point processes on $\mathbb{X}$.  
\end{theorem}

Note that since the intensity measure of $\zeta$ is $\lambda \otimes K_p$, the intensity measures $\zeta_0$ and $\zeta_1$ of respectively $\zeta(\cdot \times \{0\})$ and $\zeta(\cdot \times \{1\})$ are 
\begin{align*}
    \zeta_0(dx)=(1-p(x))\lambda(dx),\;\;\; \zeta_1(dx):=p(x)\lambda(dx). \numberthis \label{eq:int_meas_thin}
\end{align*}

We now have the necessary ingredients to give the proof of Lemma \ref{lem:nbr_PP}.
\begin{proof}[Proof of Lemma \ref{lem:nbr_PP}]

Let $\eta$ denote a Poisson point process on $\mathbb{U}$ with intensity measure given by $\lambda_d(d\mathbf{x}) \times f_W(x)dx$, where {we} recall that $\lambda_d$ denotes the Lebesgue measure on $\mathbb{R}^d$. Let $\eta_0$ be the point process $\eta + \delta_{(\mathbf{0},w)}$, i.e. $\eta_0$ is $\eta$ with the point $(\mathbf{0},w) \in \mathbb{U}$ added - a Palm version of $\eta$. Observe that the collection $\{(X_i,W_i)\}_{i \in V(\mathbb{G}^{\infty})}$ of tuples of the locations and weights of the vertices of $\mathbb{G}^{\infty}$, are distributed as the atoms of $\eta_0$. Thus, we view the atoms of $\eta_0$ as the vertices of $\mathbb{G}^{\infty}$.

For $w>0$ fixed, let $p_w:\mathbb{U} \to [0,1]$ be defined by 
\begin{align*}
   p_w(\mathbf{p})= p_w(\mathbf{x},x):= \kappa(\|\mathbf{x}\|,w,x).
\end{align*}
That is, given $\mathbf{p}=(\mathbf{x},x) \in \eta$, $p_w(\mathbf{p})$ is the probability that there is an edge between the vertices $\mathbf{p}$ and $(\mathbf{0},w)$ in $\mathbb{G}^{\infty}$.

Consider the probability kernel $K_{p_w}$ from $\mathbb{U}$ to $[0,1]$ defined by
\begin{align*}
    K_{p_w}(\mathbf{p},\cdot):=(1-p(\mathbf{p}))\delta_0(\cdot)+p(\mathbf{p})\delta_1(\cdot).
\end{align*}

Let $\eta^{(p_w)}$ be the $K_{p_w}$-marking of $\eta$, and consider $\eta^{(p_w)}_1:=\eta^{(p_w)}(\cdot \times \{1\})$, the $p_w$-thinning of $\eta$. Observe that the neighbors of $(\mathbf{0},w)$ in $\mathbb{G}^{\infty}$ {have} the same distribution as the atoms of $\eta^{(p_w)}_1$. From \eqref{eq:int_meas_thin}, the intensity measure of the Poisson point process $\eta^{(p_w)}_1$ is
\begin{align*}
    \kappa(\|\mathbf{x}\|,w,x)\lambda_d(d\mathbf{x}) \times f_W(x)dx.
\end{align*} 
In other words, the neighbors of $(\mathbf{0},w)$ are distributed as a Poisson point process on $\mathbb{U}$, with intensity function $\rho_w(\mathbf{p})$ as in \eqref{eq:density_bias}. This completes the proof.\end{proof}

\begin{remark}
Note that our proof shows that the point process of the locations and weights of the non-neighbors of $0$, given it has weight $W_0=w$, also form a Poisson point process, namely $\eta^{p_w}_0:=\eta^{(p_w)}(\cdot \times \{0\})$.
\end{remark}

\section{Inverse scaling of almost power functions: proof of Lemma \ref{lem:log_poly_inv_scaling}}
We prove the scaling of $f^{-1}(k)$ by contradiction. To this end we assume that there exists some $\delta>0$, and a sequence of reals $(k_n)_{n \geq 1}$, with $k_n \to \infty$ as $n \to \infty$, such that, for all $n \geq 1$,
\begin{align*}
    \left|\frac{f^{-1}(k_n)}{\left(\frac{a^{b}k}{c(\log{k})^{b}}\right)^{1/a}}- 1 \right|>\delta.
\end{align*}
The above inequality implies that along a subsequence
\begin{align*}
    f^{-1}(k_n)>(1+\delta)\left(\frac{a^{b}k}{c(\log{k})^{b}}\right)^{1/a} := y_n^+,
\end{align*}
or
\begin{align*}
f^{-1}(k_n)<(1-\delta)\left(\frac{a^{b}k}{c(\log{k})^{b}}\right)^{1/a} := y_n^-.
\end{align*}
Without loss of generality, let us assume the first case. The second case can be dealt with in a similar fashion (reversing all inequalities). 

Since $f$ is strictly increasing on $(t,\infty)$, and $f^{-1}(k_n) \to \infty$, for all large $n\geq 1$,
\begin{equation}\label{eq:app_1}
    k_n > f^{-1}(y_n^+).
\end{equation}
Since $f(k)=ck^{a}(\log{k})^{b}+\smallO{ck^{a}(\log{k})^{b}}$ as $k \to \infty$, and
\begin{align*}
    c (y_n^+)^{a} (\log y_n^+)^b &= c\left( (1+\delta)\left(\frac{a^{b}k}{c(\log{k})^{b}}\right)^{1/a}\right)^{a} \left(\log{\left( (1+\delta)\left(\frac{a^{b}k}{c(\log{k})^{b}}\right)^{1/a}\right)} \right)^{b}\\
    &=(1+\delta)(k_n+\smallO{k_n}),
\end{align*}
as $k_n \to \infty$, from \eqref{eq:app_1}, we have
\begin{align*}
    k_n > k_n (1 + \delta)(1 + o(1)),
\end{align*}
which is a contradiction. This completes the proof.

\paragraph{Acknowledgements.}The work of RvdH is supported in part by the Netherlands Organisation for Scientific Research (NWO) through Gravitation-grant {\sc NETWORKS}-024.002.003. NM thanks Martijn Gösgens for help with the pictures, and Joost Jorritsma for help with the simulations, based on the code \cite{JorritsmaLapinskasSpreadingGirgs}. {We thank an anonymous reviewer for a meticulous reading of the submitted version, whose observations greatly improved the paper.}

\paragraph{Availability of data and materials.} The code for the simulations that appear in Section \ref{sec:disc} have been built upon the code \cite{JorritsmaLapinskasSpreadingGirgs}. The exact code for the simulations in Section \ref{sec:disc} are available from the corresponding author on reasonable request.

\paragraph{Conflict of interest.} The authors declare that they have no conflict of interest.

\bibliographystyle{abbrvnat}
\bibliography{ref}

\end{document}